\documentclass[a4paper,11pt]{article}
\usepackage{amsfonts,amssymb,amsmath,epsf,eepic,graphics,graphicx}
\usepackage{color}
\usepackage{bbm}
\usepackage{enumerate}
\usepackage{ mathrsfs }
\usepackage[english]{babel}
\usepackage[hidelinks]{hyperref}

 \makeatletter
 \@addtoreset{enunciato}{section}
 \makeatother

 \newcounter{enunciato}[section]

 \newtheorem{ittheorem}{Theorem}
 \newtheorem{itlemma}{Lemma}
 \newtheorem{itproposition}{Proposition}
 \newtheorem{itdefinition}{Definition}
 \newtheorem{itconjecture}{Conjecture}
 \newtheorem{itassumption}{Assumption}
 \newtheorem{itcorollary}{Corollary}

 \newenvironment{remark}{\addtocounter{enunciato}{1} \noindent \textbf{Remark \thesection.\arabic{enunciato}}\,\,}
 {\smallskip }

 \newenvironment{theorem}{\addtocounter{enunciato}{1}
 \begin{ittheorem}}{\end{ittheorem}}

 \newenvironment{lemma}{\addtocounter{enunciato}{1}
 \begin{itlemma}}{\end{itlemma}}

 \newenvironment{proposition}{\addtocounter{enunciato}{1}
 \begin{itproposition}}{\end{itproposition}}

 \newenvironment{definition}{\addtocounter{enunciato}{1}
 \begin{itdefinition}}{\end{itdefinition}}

 \newenvironment{conjecture}{\addtocounter{enunciato}{1}
 \begin{itconjecture}}{\end{itconjecture}}

\newenvironment{assumption}{\addtocounter{enunciato}{1}
\begin{itassumption}}{\end{itassumption}}

\newenvironment{corollary}{\addtocounter{enunciato}{1}
\begin{itcorollary}}{\end{itcorollary}}

 \newcommand{\be}[1]{\begin{equation}\label{#1}}
 \newcommand{\ee}{\end{equation}}

 \newcommand{\bea}[1]{\begin{eqnarray}\label{#1}}
 \newcommand{\eea}{\end{eqnarray}}

 \newcommand{\bl}[1]{\begin{lemma}\label{#1}}
 \newcommand{\el}{\end{lemma}}

 \newcommand{\bt}[1]{\begin{theorem}\label{#1}}
 \newcommand{\et}{\end{theorem}}

 \newcommand{\bd}[1]{\begin{definition}\label{#1}}
 \newcommand{\ed}{\end{definition}}

 \newcommand{\bp}[1]{\begin{proposition}\label{#1}}
 \newcommand{\ep}{\end{proposition}}

 \newcommand{\bc}[1]{\begin{corollary}\label{#1}}
 \newcommand{\ec}{\end{corollary}}

 \newcommand{\bcj}[1]{\begin{conjecture}\label{#1}}
 \newcommand{\ecj}{\end{conjecture}}

\newcommand{\bass}[1]{\begin{assumption}\label{#1}}
\newcommand{\eass}{\end{assumption}}

 \newcommand{\bpr}{\begin{proof}}
 \newcommand{\epr}{\end{proof}}

 \newcommand{\bprl}[1]{\begin{proofof}{\it\ref{#1}}.\,\,}
 \newcommand{\eprl}{\end{proofof}}
 
  \newcommand{\bprp}[1]{\begin{proofofp}{\it\ref{#1}}.\,\,}
 \newcommand{\eprp}{\end{proofofp}}

 \newcommand{\bi}{\begin{itemize}}
 \newcommand{\ei}{\end{itemize}}

 \newcommand{\ben}{\begin{enumerate}}
 \newcommand{\een}{\end{enumerate}}

\newenvironment{proof}{\noindent {\em Proof}.\,\,}
{\hspace*{\fill}$\square$\medskip}
\newenvironment{proofof}{\noindent {\em Proof of Lemma\,\,}}
{\hspace*{\fill}$\halmos$\medskip}
\newenvironment{proofofp}{\noindent {\em Proof of Proposition\,\,}}
{\hspace*{\fill}$\halmos$\medskip}
\newcommand{\halmos}{\rule{1ex}{1.4ex}}

\newcommand\restr[2]{{
		\left.\kern-\nulldelimiterspace 
		#1 
		\vphantom{\big|} 
		\right|_{#2} 
}}

\newcommand\smallrestr[2]{{
		\left.\kern-\nulldelimiterspace 
		#1 
		\vphantom{|} 
		\right|_{#2} 
}}

\def \Z {{\mathbb Z}}
\def \R {{\mathbb R}}
\def \N {{\mathbb N}}

\def \ra {\rightarrow}
\def \ba {\begin{array}}
\def \ea {\end{array}}

\def \P {{\mathbb P}}
\def \E {{\mathbb E}}
\def \cX {{\mathcal X}}
\def \cL {{\mathcal L}}

\def \cC {{\mathcal C}}
\def \cP {{\mathcal P}}

\def \L {\Lambda}

\def\b {\beta}

\def\s {\sigma}

\def\L {\Lambda}
\def\G {\Gamma}
\def\e {\varepsilon}
\def\1 {\mathbbm{1}}

\def\ho {h_\mathrm{2}}
\def\he {h_\mathrm{1}}
\def\Lo {\L_\mathrm{2}}
\def\Le {\L_\mathrm{1}}

\definecolor{darkgreen}{rgb}{0,.6,0}
\definecolor{darkagenta}{rgb}{.5,0,.5}
\definecolor{darkred}{rgb}{1,0,0}
\definecolor{darkblue}{rgb}{0,0,.4}
\definecolor{black}{rgb}{0,0,0}
\definecolor{gray}{rgb}{.4,.4,.4}
\definecolor{white}{rgb}{0.99,0.99,0.99}
\definecolor{geel1}{rgb}{1,1,0.3}

\numberwithin{equation}{section}
\begin{document}


\title{On the metastability in three \\
 modifications of the   Ising model}

\author{
\renewcommand{\thefootnote}{\arabic{footnote}}
K.\ Bashiri
\footnotemark[1]
}

\footnotetext[1]{
Institut f\"ur Angewandte Mathematik,
Rheinische Friedrich-Wilhelms-Universit\"at, 
Endenicher Allee 60, 53115 Bonn, Germany
}

\date{\today}
\maketitle


\begin{abstract}

We consider three extensions of the standard 2D Ising model with Glauber dynamics on a finite torus at low temperature. 
The first model (see Chapter \ref{S2}) is an anisotropic version, where the interaction energy takes different values on vertical and on horizontal bonds.
The second model (Chapter \ref{S4}) adds next-nearest-neighbor attraction to the standard Ising model. 
And the third model  (Chapter \ref{pm}) associates different alternating signs for the magnetic fields on even and odd rows. 
All these models have already been studied, and results concerning metastability have been established using the so-called \textit{path-wise approach} (see \cite{KO92},\cite{KO94},\cite{NO96}).
In this text, we extend these earlier results, and apply the \emph{potential-theoretic approach} to metastability to obtain more precise asymptotic information on the transition time from the metastable phase to the stable phase.
\let\thefootnote\svthefootnote
\let\thefootnote\relax\footnotetext{
	{\it Date.} 
	May 1, 2017.}
\let\thefootnote\svthefootnote
\let\thefootnote\relax
\let\thefootnote\relax
\let\thefootnote\svthefootnote
\let\thefootnote\relax\footnotetext{
	{\it Preprint.} 
	https://arxiv.org/abs/1705.07012}
\let\thefootnote\svthefootnote

\let\thefootnote\relax\footnotetext{
The research in this paper is partially supported by the German Research Foundation in the Collaborative Research Center 1060 "The
Mathematics of Emergent Effects", and the Bonn International Graduate School in Mathematics (BIGS) in the Hausdorff
Center for Mathematics (HCM)}
\let\thefootnote\svthefootnote
\end{abstract}


\section*{Introduction}
%
%
%

In many physical, biological or chemical evolutions, one can observe a phenomenon called \textit{metastability}.
If the states of the system are associated to an energy functional, this phenomenon can be described as follows.
For a relatively long time, the state of the system resides around a local minimum of the energy landscape, which is not the global minimum. 
This state is called the \textit{metastable state}.
However, under thermal fluctuations and after many unsuccessful attempts the system can finally free itself from this valley in the energy landscape and it manages the crossover to the global minimum, which is called the \textit{stable state}.
Often, this crossover is triggered by reaching a \textit{critical state} in the system.
An example is   \textit{over-saturated water vapor}, where below critical temperature, the formation of a critical droplet is needed to achieve the transition from the gas-phase to the liquid-phase.
An analogue situation holds for  \textit{over-cooled liquids} and for
 \textit{magnetic hysteresis}.

Through the last decades many mathematical models have been considered to study this phenomenon. 
One is mostly interested in 
\begin{enumerate}[i)]
	\setlength\itemsep{-0.3em} 
	\item the average transition time from the metastable to the stable state,
	\item  the estimate in probability and the exponential distribution of this transition time,
	\item the  typical paths for the transition from the metastable to the stable state, and 
	\item  in showing that the critical states has to be passed in order to make this transition.\label{S1.1}
\end{enumerate}

Mainly two methods have been crystallized to be very fruitful to tackle these problems. 
The first one is the \textit{path-wise approach}, initiated by Cassandro, Galves, Olivieri and Vares in \cite{CGOV84}. 
Motivated by the \textit{Freidlin and Wentzell theory}, one uses large deviation estimates on the path space to identify the most likely paths of the system for the transition from the metastable to the stable state.
The advantage is a very detailed description of the tube of typical paths for the transition, but at the same time the average transition time can only be computed up to a multiplicative factor of the order $e^{\e \beta}$ as $\beta \ra \infty$, where $\beta$ is the inverse temperature and $\e > 0 $ is independent of $\beta$ and can be chosen arbitrary small. 
For an extensive treatment on the path-wise approach to metastability, the reader is referred to \cite{OV04}, \cite{MNOS}, \cite{CNS} or \cite{FMNSS}. 

The second method is the \textit{potential-theoretic approach to metastability}, which was initiated by Bovier, Eckhoff, Klein and Gayrard in \cite{BEGK01} and \cite{BEGK02}.
Here, one uses \textit{potential theory} to rewrite the average transition time in terms of quantities from \textit{electric networks}, namely \textit{capacities}.
Now, using variational principles for the capacity, the average transition time can be computed up to a multiplicative error that tends to one as $\beta \ra \infty$, which provides a sharp estimate.
This method is also the basis of this text, and in Chapter \ref{ModelChap} we will shortly review a general recipe to obtain metastability results for a stochastic process on a finite graph, whose dynamics are given through a metropolis algorithm.
This general recipe is based on the paper \cite{BM02} by Bovier and Manzo.
For a more detailed overview on the potential-theoretic approach to metastability, we refer to the 2015 monograph \cite{BdH15} by Bovier and den Hollander (especially Chapter 16).

Probably the easiest application of these methods, where one can rigorously investigate metastable  behavior, is the two-dimensional standard Ising model on a finite torus in the low temperature regime. 
Neves and Schonmann applied 1991 in \cite{NS91} the path-wise approach to this model, which was later rewritten in the Chapters 7.1--7.5 of \cite{OV04}.
And in the year 2002 the potential-theoretic approach was used in \cite{BM02} by Bovier and Manzo (see also Chapter 17 of \cite{BdH15}).
Moreover, several other settings and regimes in the Ising model have been considered   as well. 
For example,  the Ising model on $ \Z^d $  was considered in \cite{DS} ($ d=2 $) and in  \cite{CF} ($ d\geq 3 $), and the regime, where the magnetic field tends to zero was treated in \cite{S}.
Of course, the metastable analysis goes far beyond the Ising model, and for numerous other models, such as coupled diffusion process or lattice gas models, a metastable behavior was studied rigorously.
For an extensive historical review on this we refer to  \cite{BdH15},  \cite{OV04}, \cite{MNOS} and references therein.
 
In this paper we study three modifications of the Ising model. 
Roughly speaking,  the crucial difference between all three models and  the standard Ising model is the fact that we lose the applicability of  \emph{isoperimetrical inequalities}.
Namely, in the Ising case, for a given number of \emph{up-spins}, the configurations with minimal energy  are those droplets of up-spins whose  shape  is given by a  square (or a quasi-square) with a possible bar of up-spins attached to one of its sides.
Here we do not have this property.
Instead we need to look at the \emph{stability} of certain classes of configurations separately in order to specify the metastable and the critical state rigorously. 
The path-wise approach  has already been applied to  these models in  \cite{KO92}, \cite{KO94} and \cite{NO96}, respectively.
In the Chapters 7.7--7.10 of \cite{OV04}, a brief overview on these three papers is given. 
Here we complement these results and apply the potential-theoretic approach. 

In the following chapter we  introduce these three models,  formulate the main results, and provide a more detailed comparison of the results in this paper and the results from  \cite{KO92}, \cite{KO94} and \cite{NO96}.
The proofs are moved to the Chapters \ref{S2}, \ref{S4} and \ref{pm}, respectively.

\section{Main results}
\label{ModelChap}

This chapter is organized as follows.
In Section \ref{S1.2}  we introduce the setting and the results of Chapter 16 in \cite{BdH15}.\footnote{
	The setting in \cite{BdH15} is more general, but to keep it as simple as possible, we restrict to our situation of a dynamical spin-flip model on the two-dimensional lattice.}
More precisely, we define a dynamical spin-flip model on the two-dimensional lattice, which is  driven by a general energy function.
Furthermore, we provide definitions concerning the geometrical properties of the energy landscape. 
At the end of that section, we state the so-called \emph{metastability theorems}, which are the key results on which we rely in this paper.

In Section \ref{resultsec} we introduce the three modifications of the Ising model that we are studying in this paper and state the main results. 
In fact, these models are given by the abstract set-up from Section \ref{S1.2} but with a specific energy function.
In that section we also  compare our results and our approach with those from the papers \cite{KO92}, \cite{KO94} and \cite{NO96}.

In Section \ref{SecGeo} we list some further definitions that will be used in this paper.


\subsection{The abstract set-up  and the metastability theorems from \cite{BdH15}}
\label{S1.2}
%
%
%

\begin{figure}[htbp]
 	\centering
 	\includegraphics[scale=0.25]{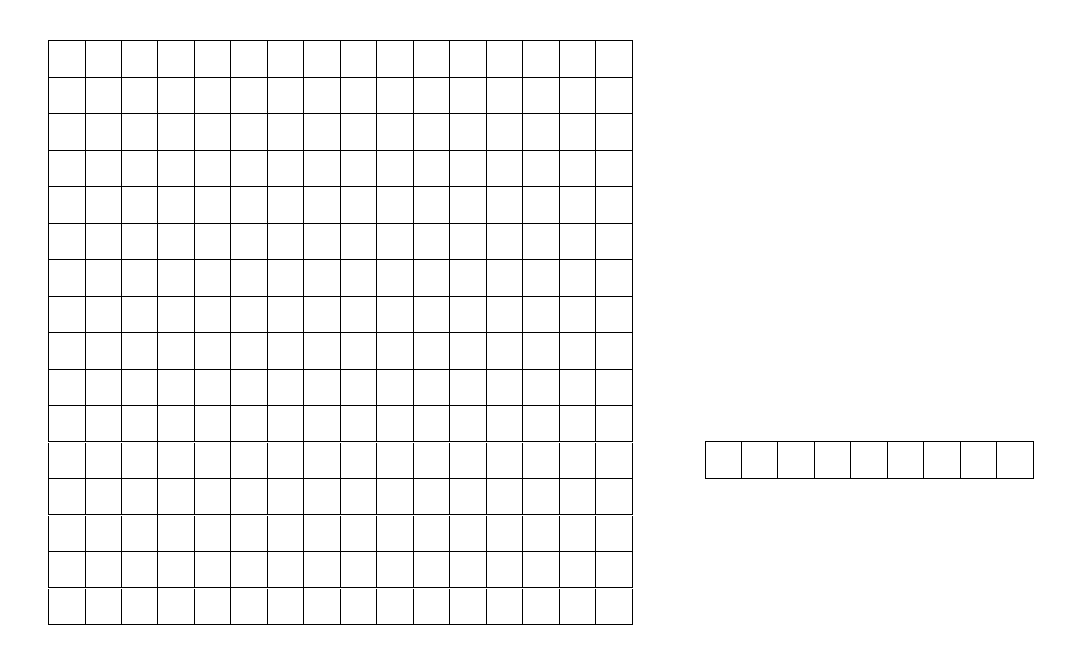}	
 	\vspace{-0.5cm}
 	\caption{{\footnotesize Geometric representation of a configuration that assigns to each site in $\L$ the spin-value $-1$  except on a square of size $16 \times 16$ and a rectangle of size $1 \times 8$}}
 	\label{Example}
 \end{figure}

 Let $\L \subset \Z^2$ be a finite, square box with periodic boundary conditions, centered at the origin. 
 $S=\{-1,1\}^\L$ will be called the \textit{configuration space}.
 An element $\s \in S$ is called \textit{configuration}, and at each \textit{site} $x\in\L$, $\s(x)\in \{-1,1\}$ is called the \textit{spin-value} at $x$.
 By abuse of notation, we often identify each configuration $\s \in S$ with the sites that have spin value $+1$, i.e.
 \be{confgeo}
 \s \equiv  \{ x\in \L \ | \  \s(x)=+1 \}.
 \ee
 Moreover, we represent $\s$ geometrically by identifying  each $x \in \s$ with $  \s(x)=+1 $ with a closed unit square centered at $ x $.
 See Figure \ref{Example} for an example.
 
 The \emph{energy} of the system is given by a \textit{Hamiltonian} $\mathrm{H}:S\ra \R$.
 If $\b >0$ is the \textit{inverse temperature}, the  \textit{Gibbs measure} associated with $\mathrm{H}$ and $\b$ is given by
 \begin{align}
 \mu_\b(\s)	= \frac{1}{Z_\b}e^{-\b \, \mathrm{H}(\s)}, \qquad \text{for	} \s \in S,
 \end{align}
 where $ Z_\b $ is a normalization constant called \textit{partition function}.
 
 For $ \s \in S $ and $x\in\L$ we define $\s^x \in S$ by
 \begin{align}
 \s^x(y)=
 \begin{cases}
 \s(y)			&:y\neq x, \\
 -\s(x)			&:y=x.
 \end{cases}
 \end{align}
 For all $ \s,\s'\in S $, we say that $ \s $ and $\s' $ \emph{communicate} and write $\s \sim \s'$ if there exists  $x\in\L$ such that $\s^x =\s'$.
 This induces a graph structure on $S$ by defining an edge between each $\s,\s' \in S$ whenever  $\s \sim \s'$.
 
 The dynamics of the system is given by the continuous time Markov Chain $(\s_t)_{t\geq 0}$ on $S$, whose generator $ \cL_\b $ is given by
 \begin{align}
 (\cL_\b f) (\s) = \sum_{x \in \L}  c_\b (\s,\s^x ) (f(\s^x) - f(\s)),
 \end{align} 
 where $f:S\ra \R$ is a function and
 \begin{align}
 c_\b (\s,\s' )=
 \begin{cases}
 e^{-\b \max  \{0,\mathrm{H}( \s') - \mathrm{H}(\s) \} } 			&:\s \sim \s' , \\
 0		&:\text{else}.
 \end{cases}
 \end{align}
 Notice that for $\b=\infty$ only moves to configurations with lower or equal energy are permitted. 
 Moreover, one can immediately see that the following \textit{detailed balance condition} holds:
 \be{reversible}
 \mu_\b(\s) c_\b(\s,\s') = \mu_\b(\s') c_\b(\s',\s) \qquad \forall\,\s,\s' \in \cX_\b^{(n_\b)}. 
 \ee
 Hence, the dynamics is \textit{reversible} with respect to the Gibbs measure. 
 The law of $(\s_t)_{t\geq 0}$ given that $\s_0= \s \in S$ will be denoted by $\P_\s$, and for a set $A\subset S$, we denote its \emph{first hitting time after the starting configuration has been left} by $\tau_\mathrm{A}$, i.e.
 \begin{align}
 \tau_\mathrm{A}= \inf \{	t> 0|\ \s_t \in A, \, \exists\, 0<s<t : \s_s \neq \s_0	\}.
 \end{align}
 If $A = \{\s\}$ for some $\s\in S$, then we write $ \tau_\s = \tau_\mathrm{A} $.
 \bd{geom}
 \begin{itemize}
 	\item[i)] Let $ \s,\s' \in S $. The \emph{communication height} between $ \s $ and $ \s' $ is defined by
 	\begin{align}
 	\Phi(\s,\s' ) = \min_{\gamma:\s \ra \s'} \max_{\eta \in \gamma} \mathrm{H}(\eta),
 	\end{align}
 	where the minimum is taken over all finite paths $ \gamma$ of allowed moves  in $S$ going from $\s $ to $ \s' $.
 	
 	\item[ii)] Let $ \s,\s' \in S $. A finite path $ \gamma:\s \ra \s'  $ is called \emph{optimal path} between $ \s $ and $ \s' $ if
 	\begin{align}
 	\Phi(\s,\s' ) =  \max_{\eta \in \gamma} \mathrm{H}(\eta).
 	\end{align}
 	The set of all optimal paths between $ \s $ and $ \s' $ is denoted by $ (\s \ra \s' )_{\mathrm{opt}} $.
 	
 	\item[iii)] Let $ \s \in S $. The \emph{stability level} of $ \s $  is defined by
 	\begin{align}
 	V_\s =  \min_{\eta \in S: \mathrm{H}(\eta) < \mathrm{H}(\s)} \Phi(\s,\eta) - \mathrm{H}(\s).
 	\end{align}
 	Moreover, for $ V \in \R $, we define 
 	\begin{align}\label{StabilityEq}
 	S_V =  \{ \s \in S\,|\, V_\s > V\},
 	\end{align}
 	which is the set of all configurations, whose stability level is greater than $ V $.
 	\item[iv)] The set of \emph{stable configurations} in $S$ is defined by:
 	\begin{align}
 	S_{\mathrm{stab}}= \{ \s \in S\  |\  \mathrm{H}(\s) = \min_{\eta \in S} \mathrm{H}(\eta)				\}.
 	\end{align}
 	
 	\item[v)] The set of \emph{metastable configurations} in $S$ is defined by:
 	\begin{align}
 	S_{\mathrm{meta}}= \{ \s \in S\  |\  V_\s = \max_{\eta \in S\setminus S_{\mathrm{stab}} } V_\eta  \,  	\}.
 	\end{align}
 	
 	\item[vi)] Let $ (m,s) \in S_{\mathrm{meta}} \times S_{\mathrm{stab}}$.  The \emph{energy barrier} $ \G^\star(m,s) $ between $ m $ and $ s $ is defined by
 	 \begin{align}
 	\G^\star(m,s)  = \Phi(m,s) - \mathrm{H}(m). 
 	\end{align}
 \end{itemize}
 \ed
 
 Note that by Theorem 2.4 in \cite{CN} we have that $ \G^\star (m,s) =  \max_{\eta \in S\setminus S_{\mathrm{stab}} } V_\eta =: \G^\star $ for all $ (m,s) \in S_{\mathrm{meta}} \times S_{\mathrm{stab}}$.
In the following definition we  introduce the notion of a  \emph{critical configuration}.
This quantifies the idea of the critical state from the introduction. 

 \begin{definition}[Definition 16.3 in \cite{BdH15}]\label{Crit}
 Let $ (m,s) \in S_{\mathrm{meta}} \times S_{\mathrm{stab}}$. 
 Then $(\mathcal{P}^\star(m,s), \cC^\star (m,s) )$ is defined as the maximal subset of $ S\times S $ such that 
 \begin{enumerate}[1.)]
 	\item $ \forall \s \in \cP^\star (m,s) \ \exists \s' \in \cC^\star (m,s) \  : \  \s \sim \s' $, and \\
 	$ \forall \s' \in \cC^\star (m,s) \  \exists \s \in \cP^\star (m,s) \ : \  \s \sim \s' $,	
 	\item 	$ \forall \s \in \cP^\star (m,s) \   : \ \Phi(m,\s) < \Phi(\s,s) $,
 	\item	$ \forall \s' \in \cC^\star (m,s) \  \exists \gamma: \s' \ra s  \, : \,   \max_{\eta \in \gamma } \mathrm{H}(\eta ) - \ \mathrm{H}(m) \leq \G^\star \, , \, \Phi(m,\eta) \geq \Phi(\eta,s) \,  \forall \eta \in \gamma $.
 \end{enumerate}
 We call $ \cP^\star (m,s) $ the \emph{set of protocritical configurations} and $ \cC^\star (m,s) $ the \emph{set of critical configurations}.
  \end{definition}

 The results from  Section \ref{resultsec}  will be based on the following metastability theorems (see Theorem 16.4 -- 16.6 in \cite{BdH15}).
 These will hold subject to the  hypothesis
 \begin{itemize}
 	\item[(H1)] $ S_{\mathrm{meta}}= \{ m 	\} $ and $ S_{\mathrm{stab}}= \{ s \} $,
 \end{itemize}
 where $ m,s \in S $.
 One challenge in the Chapters \ref{S2}--\ref{pm}  is to verify this hypothesis for the three specific models.
 Under (H1), it would not lead to confusions if we abbreviate $\cP^\star= \cP^\star (m,s) $ and  $ \cC^\star= \cC^\star (m,s) $.

 \begin{theorem}[Theorem 16.4 in \cite{BdH15}]\label{T1}
 	Consider $(\mathcal{P}^\star ,\cC^\star)  $ from Definition \ref{Crit}.
 	Suppose  (H1).
 	Then, 
 	\begin{itemize}
 		\item[a)]	$ \lim_{\b \ra \infty} \P_m [	\tau_{\cC^\star} < \tau_s \ | \ \tau_s < \tau_m	] = 1 $, and
 		\item[b)]	if, moreover, the following assumption holds
 		\bi
 		\item[\emph{(H2)}] \qquad \qquad	$  \s'\ra |\{ \s \in \cP^\star \  : \  \s \sim \s' \}| $ is constant on $ \cC^\star$,
 		\ei
 		then for all $ \chi	 \in \cC^\star $, it holds: 
 		$
 		\lim_{\b \ra \infty} \P_m [	\s_{\tau_{\cC^\star}} = \chi	] = \frac{1}{|\cC^\star|}.
 		$
 	\end{itemize}
 \end{theorem}
 Theorem \ref{T1} says that the set of critical configurations has to be reached in order to cross over from the metastable to the stable configuration.
 If the additional assumption holds, then part b) of Theorem \ref{T1} says that the entrance into $\cC^\star$ is uniformly distributed on $\cC^\star$. 
 
 \begin{theorem}[Theorem 16.6 in \cite{BdH15}]\label{T2}
 	Subject to (H1), it holds that 
 	\begin{itemize}
 		\item[a)]	$ \lim_{\b \ra \infty} \lambda_\b \E_m [	\tau_s 	] = 1 $, where $ \lambda_\b $ is the second largest eigenvalue of $- \cL_\b $, and
 		\item[b)]	
 		$ \lim_{\b \ra \infty} \P_m [	\tau_s > t \cdot \E_m [	\tau_s 	]	] = e ^ {-t} $ for all $ t \geq 0 $.
 	\end{itemize}
 \end{theorem}
 Theorem \ref{T2} represents the average transition time of the system in terms of the spectrum of its generator and part b) yields the asymptotic exponential distribution of $  \tau_s $.
 
 \begin{theorem}[Theorem 16.5 and Lemma 16.17 in \cite{BdH15}]\label{T3}
 	Suppose (H1). Then, 
 	\begin{itemize}
 		\item[a)]	
 		there exists a constant $ K \in (0, \infty)$ such that 
 		$ \lim_{\b \ra \infty} e^{-\b \G^\star} \E_m [	\tau_s 	] = K $, and
 		\item[b)]	define 
 		\begin{itemize}
 			\item $ S^\star \subset S  $ be the subgraph obtained by removing all vertices $\eta$ with $ H(\eta) > \G^\star + H(m) $ and all edges incident to these vertices,
 			\item $ S^{\star \star} \subset S^\star  $ be the subgraph obtained by removing all vertices $\eta$ with $ H(\eta) = \G^\star + H(m) $ and all edges incident to these vertices,
 			\item $ S_m = \{  \eta \in S\  |\   \Phi(m,\eta) < \Phi(\eta,s) = \G^\star + H(m)   \} $,
 			\item $ S_s = \{  \eta \in S\  |\   \Phi(\eta,s)< \Phi(m,\eta) = \G^\star + H(m)   \} $,
 			\item 
 			$ S_1, \dots , S_I \subset S^{\star \star} $  be such that $ S^{\star \star} \setminus (S_m  \cup S_s   ) = \cup_{i=1}^I S_i $ and each  $S_i$ is a maximal set of communicating configurations,
 		\end{itemize} 
 		then,
 		\begin{align}\label{EqPrefactor}
 		\frac{1}{K} = \min_{C_1,\dots,C_I \in [0,1]}
 		\min_{\substack{ h:S^\star \ra [0,1] \\
 				\smallrestr{h}{S_m} =1, \smallrestr{h}{S_s} =0, \smallrestr{h}{S_i} =C_i \, \forall i }}
 		\frac{1}{2} \sum_{\eta,\eta' \in S^\star } \mathbbm{1}_{\{	\eta \sim \eta'	 \}} [h(\eta) - h(\eta')]^2.
 		\end{align}
 		Note that the first minimum runs over all constants $ C_1,\dots,C_I \in [0,1] $.
 	\end{itemize}
 \end{theorem}
 Theorem \ref{T3} yields the precise asymptotics of the average transition time  and provides a variational formula to compute the pre-factor. 
Note that the sum in   \eqref{EqPrefactor} is taken, in particular, also  over elements from $ \cup_{i=1}^I S_i $. 
However, in the three models that we consider these terms are  negligible.
Indeed, for the lower bound we  just use that all terms in the sum in \eqref{EqPrefactor} are non-negative.
For the upper bound we restrict the minimum to a suitable class of functions $ h $, where a transition $  \eta \sim \eta' $ with $ h(\eta ) \neq h(\eta' )  $ is not possible if $ \eta \in S^{\star \star }\setminus \cP^\star $ or $\eta' \in S^\star \setminus \cC^\star$.
For more details we refer to Section \ref{Step6A}.
 
 \subsection{The models and the main results}
 \label{resultsec}

The first goal of this paper is to verify Theorem \ref{T1}, \ref{T2} and \ref{T3} for three specific models that will be introduced below in the Subsections \ref{resultsecA}--\ref{resultsecpm}.
Thus we have to show that the conditions of those theorems are satisfied.
In all three models we will have  that 
\begin{align}
m= \boxminus , \qquad \text{and} \qquad
s= \boxplus,
\end{align}
where $ \boxminus \in S $ is the configuration, where all spin values are $ -1 $ and $ \boxplus $ is the configuration with all spin values being $ +1 $.
The second goal  of this paper is to compute for each case the precise value of the pre-factor $ K $ from Theorem \ref{T3}.
Hence, for each model we have to 
\begin{itemize} 
	\setlength\itemsep{-0.3em} 
	\item compute $  \Phi(\boxminus,\boxplus) - \mathrm{H}_\mathrm{A}(\boxminus)$,
	\item identify the sets $ \cP^ \star $ and	$ \cC^ \star $,
	\item verify hypothesis (H1) and if possible hypothesis (H2), and 
	\item compute $ K $.\
\end{itemize}

These tasks are treated in the Chapters \ref{S2}--\ref{pm}.
In this section we only introduce the models and formulate the results. 
However, before we introduce the models, we  comment on the results that were already obtained in \cite{KO92}, \cite{KO94} and \cite{NO96}.
In all these papers  
Theorem \ref{T1} a) has already been established for the respective models; see\ Theorem 1 in \cite{KO92}, Theorem 1 in \cite{KO94} and Theorem 1 i) in \cite{NO96}. 
Moreover, an estimate in probability of $ 	\tau_\boxplus	  $ is proven  and  the typical paths for the transition from $ \boxminus $ to $ \boxplus $ are identified; see Theorem 2 and 3 in \cite{KO92}, Theorem 2 and 3 in \cite{KO94}, and  Theorem 1 ii), Section 4 and 5 in \cite{NO96}.
Then, using standard techniques (cf.\ Theorem 6.30  and (6.171) in \cite{OV04}),
one can use these results to obtain the value of $  \Phi(\boxminus,\boxplus) - \mathrm{H}_\mathrm{A}(\boxminus)$, Theorem \ref{T2} b) and the exponential asymptotics of $ \E_\boxminus [	\tau_\boxplus	]  $ (i.e.\ the asymptotics without the pre-factor $ K $). 
Hence, we provide here a new approach to prove Theorem \ref{T1} a),  Theorem \ref{T2} b) and to compute the  value of $  \Phi(\boxminus,\boxplus) - \mathrm{H}_\mathrm{A}(\boxminus)$. Moreover, we prove a sharper estimate for $ \E_\boxminus [	\tau_\boxplus	]  $.

 \subsubsection{Anisotropic Ising model}
  \label{resultsecA}
  
The first model we study is the same model as in \cite{KO92}, where the interaction between  neighboring spins is anisotropic in the sense that the attraction on horizontal bonds is stronger than on vertical bonds.
More precisely, the Hamiltonian here is explicitly given by
 \be{HAIntro}
 \mathrm{H}_\mathrm{A}(\s) = - \frac{J_H}{2} \sum_{(x,y) \in \L_H^\star} \s(x)\s(y)
 - \frac{J_V}{2} \sum_{(x,y) \in \L_V^\star} \s(x)\s(y)
 - \frac{h}{2} \sum_{x \in \L} \s(x),
 \ee
 where $\s \in S$, $ J_H>J_V>0,\ h>0 $, $ \L_H^\star $ is the set of \textit{unordered horizontal nearest-neighbor bonds} in $ \L $ and $ \L_V^\star$ is the set of \textit{unordered vertical nearest-neighbor bonds} in $ \L $.
 Here and in the following the subscript $ \mathrm{A} $ is added to remind that we are in the anisotropic case. 
 The critical length in this model is given by
 \begin{align}\label{CritlengthAIntro}
 L_V^\star = \left\lceil   \frac{2J_V}{h} \right\rceil.
 \end{align}
 
 We now formulate the main result for this model.
 For a more precise formulation and the proof we refer to Chapter \ref{S2}. 
 \bt{ThmAIntro}
 Under Assumption \ref{AssA}, the pair $ (\boxminus,\boxplus) $ satisfies (H1) and (H2) so that Theorems \ref{T1}--\ref{T3} hold for the anisotropic Ising model. 
 Moreover, $ \cP^ \star  $ and $   \cC^\star  $ are given by the set of configurations given in Definition \ref{defiA} and Theorem \ref{ThmA},  
\begin{align}
\begin{split}
  \Phi(\boxminus,\boxplus) - \mathrm{H}_\mathrm{A}(\boxminus) &= 2 L_V^\star (J_H+J_V) - h (1+(L_V^\star-1)L_V^\star)    , \text{ and } \\
K^{-1}&=   \frac{4(2L_V^\star-1)}{3}|\L|  .  
\end{split}
\end{align}
 \et

 \subsubsection{Ising model with next-nearest-neighbour attraction}
  \label{resultsecNN}
  
  In the second model we consider in this paper, we allow next-nearest-neighbor attraction, i.e.\ two spins that have euclidean distance of $\sqrt{2}$ feel an interaction energy, which is strictly less than the interaction energy between nearest-neighbor bonds.
  This has the physical intuition that  next-nearest-neighbor attraction is seen as a perturbation of nearest-neighbor attraction.
  An interesting fact is that the local minima of the energy landscape are given by droplets of \emph{octagonal shape}.
  For the path-wise approach to this model we refer to \cite{KO94}. 
  
Here the Hamiltonian is given by
 \be{HNNIntro}
 \mathrm{H}_\mathrm{NN}(\s) = - \frac{\tilde{J}}{2} \sum_{(x,y) \in \L^\star} \s(x)\s(y)
 - \frac{K}{2} \sum_{(x,y) \in \L^{\star\star}} \s(x)\s(y)
 - \frac{h}{2} \sum_{x \in \L} \s(x),
 \ee
 where $\s \in S$, $ \tilde{J}>K,\ h>0 $, $ \L^\star $ is the set of  \textit{unordered nearest-neighbor bonds} in $ \L $ and $ \L^{\star\star}$ is the set of \textit{unordered next-nearest-neighbor bonds} in $ \L $, i.e 
 \begin{align}  \label{resultsecNNEq}
 \L^{\star\star} = \big\{	\{x,y\} \in \L^2\ \big|\ |x-y|=\sqrt{2} 	\big\}. 
 \end{align}
 Also here the subscript $ \mathrm{NN} $ is added to remind that we are in the case with next-nearest-neighbor attraction. 
 Set $  J= \tilde{J}+2K $.
 The critical lengths in this model will be given by
 \begin{align}\label{CritlengthsNNIntro}
 \ell^\star = \left\lceil  \frac{2K}{h} \right\rceil  \qquad \text{and} \qquad D^\star = \left\lceil   \frac{2J}{h} \right\rceil  \qquad \text{and} \qquad L^\star = D^\star - 2(\ell^\star -1).
 \end{align}
We now formulate the main result for this model. A more precise formulation and the proof are given in  Chapter \ref{S4}.  
 \bt{ThmNNIntro}
 Under Assumption \ref{AssNN}, the pair $ (\boxminus,\boxplus) $ satisfies (H1) and (H2) so that Theorems \ref{T1}--\ref{T3} hold for the Ising model with next-nearest-neighbor attraction.
 Moreover, $  \cP^ \star $ and  $\cC^\star  $ are given by the set of configurations given in Definition \ref{defiNN} and Theorem \ref{ThmNN},  
\begin{align}
\begin{split}
  \Phi(\boxminus,\boxplus) - \mathrm{H}_\mathrm{NN}(\boxminus) &=\mathrm{H}_\mathrm{NN}(Q(D^\star - 1,D^\star))  +2J-4K-h  , \text{ and }\\
  K^{-1}&=   \frac{4(2L^\star-5)}{3}|\L|  .
\end{split}
\end{align}
 \et
 \subsubsection{Ising model with alternating magnetic field}
 \label{resultsecpm}
 
In the third modification of  the standard Ising model, the magnetic field is allowed to take alternating signs and absolute values on even and on odd rows. 
 The path-wise approach has been applied to this model in \cite{NO96}.
The Hamiltonian here is given by
\be{HpmIntro}
\mathrm{H}_{\pm}(\s) = - \frac{J}{2} \sum_{(x,y) \in \L^\star} \s(x)\s(y)
+ \frac{\ho }{2} \sum_{x \in \L_\mathrm{2}} \s(x)
- \frac{h_\mathrm{1}}{2} \sum_{x \in \L_\mathrm{1}} \s(x),
\ee
where $\s \in S$, $ J ,\ho,\he >0 $,  $ \Lo = \{	 (x_1,x_2)\in \L \, | \, x_2 \text{ is odd}	\} $ are the \emph{odd rows} in $ \L $,  $ \Le = \L\setminus\Lo $ are the \emph{even rows} and   $ \L^\star $ is the set of \textit{unordered nearest-neighbor bonds} in $ \L $.
The critical lengths in this model will be given by
\begin{align}\label{CritlengthspmIntro}
l_b^\star = \left\lceil  \frac{\mu}{\e} \right\rceil  
\qquad \text{and} \qquad l_h^\star = 2 l_b^\star - 1,
\end{align}
where
\begin{align}
\begin{split}
& \e = \he - \ho , \quad \text{and} \\
&\mu = 2J - \ho.
\end{split}
\end{align}
$ l_b^\star $ will be the length of the basis of the critical droplet, and $ l_h^\star $ will be its height. 
We now state the main result for this model. More details and the proof are given in Chapter \ref{pm}.
\bt{ThmpmIntro}
Under Assumption \ref{Asspm}, the pair $ (\boxminus,\boxplus) $ satisfies (H1) so that Theorem \ref{T1} a), Theorem \ref{T2} and Theorem \ref{T3} hold for the Ising model with alternating magnetic field.
Moreover, 
$  \cP^ \star $ and  $\cC^\star  $ are given by the set of configurations given in Definition \ref{defipm} and Theorem \ref{Thmpm},
\begin{align}
\begin{split}
   \Phi(\boxminus,\boxplus) - \mathrm{H}_{\pm}(\boxminus) &=  4J \, l_b^\star + \mu (l_b^\star - 1) - \e ( l_b^\star  (l_b^\star - 1) + 1 )    , \text{ and }\\
K^{-1}&=  
\frac{14\,( l_b^\star-1)}{3}|\L| .
\end{split}
\end{align}
\et

%

\subsection{Further definitions}\label{SecGeo}

We conclude this chapter with some definitions hat are used in all three situations in Chapters \ref{S2}--\ref{pm}.

 \begin{itemize}
	\item 
	For $x \in \R$, $\lceil x \rceil$ denotes the smallest integer greater than $x$.
	\item
	For $l_1,l_2 \in \N$, $R(l_1 \times l_2)$ denotes the set of all configurations consisting of a single rectangle with horizontal length $l_1$ and vertical length $l_2$ somewhere on the torus $\L$.
	An element $ \s \in R(l_1 \times l_2) $ is called \emph{rectangle} and will often be denoted by $  l_1 \times l_2$, since usually we can ignore the position of the rectangle in the torus.
	For this reason, by abuse of notation, we often identify the whole set $R(l_1 \times l_2)$ with $l_1 \times l_2$.
	We also define $  R  (l_1 , l_2) =R(l_1 \times l_2)\cup R(l_2 \times l_1)$. 
	If $ |l_1 - l_2|= 1 $ or $ |l_1 - l_2|= 0 $, then $l_1 \times l_2$ is called \emph{quasi-square} or \emph{square}, respectively.
	$ 1 \times l_2 $ is called \emph{vertical bar} or \emph{column} and $ l_1 \times 1 $ is called \emph{horizontal bar} or \emph{row}.
	\item
	For a rectangle $ R \in S $, we denote by  $ P_HR\in \N$ its \emph{horizontal length}, and by $ P_VR \in \N$ its \emph{vertical length}.
	\item
	For $ \s \in S $, let $ |\s|  $ be the \emph{area} of $\s$, i.e.\ its number of $(+1)$--spins.
	Further, $ \partial (\s) $ is the Euclidean boundary of $ \s $ in its geometric representation and $ |\partial (\s)| $ denotes the \emph{perimeter}, i.e.\ the length of $ \partial (\s) $. 
	\item 
	Let $ \s \in S $.
	We say that  $ \s  $  is \emph{connected} if $\s \setminus \partial (\s) $  is connected in the Euclidean space $ \R^2 $.
	\item
	Let $ \s \in S $. 
	A \emph{cluster} of $ \s $ is a maximally connected component of $ \s $.
	\item 
	Two droplets on the torus are called \emph{isolated} if their Euclidean distance is greater or equal to $ \sqrt{2} $.	
	\item
	Let $ \s \in S $ and $ x\in \L $ be such that $ \s(x)=+1 $. Then $x$ is called \emph{protuberance} if $ \sum_{y \in \L: |y-x|=1} \s(y) = -2 $.
	\item
	Let $ \s \in S $ be connected and $ l$  be either a vertical bar or a horizontal bar. Then $l$ is called  \emph{attached} to $ \s $ if  for all $ x \in l $ there exists $ y \in \L\setminus l $ such that  $ |y-x|=1$ and $ z \in \s $ such that  $ |z-x|=1$. 
	\item 
	If $ \s \in S $ consists of a single, connected droplet, then $ R(\s) $ is the smallest rectangle that contains $ \s $.
	\item
	A \emph{row} or a \emph{column} of  a connected configuration $ \s\in S $ is defined as the intersection of a row or a column of $ \L $ with $ \s $.
	\item 
	$ \s \in S $ is called a  \emph{local minimum} of  $ \mathrm{H} $ if $ \mathrm{H}(\s^x) > \mathrm{H}(\s) $ for all $ x \in \L $. 
	\item
	For $ A \subset S$, let $ \partial ^+ A= \{	 \s \in S \setminus A \, | \, \exists \s' \in S : \s \sim \s'	\}$ denote the \emph{outer boundary} $A$.
	We also define $ A^+ = A \cup  \partial ^+ A $.
	Moreover, if $ \eta \in S $, then $ A \sim \eta \subset S$ is defined by $A \sim \eta= \{	 \s \in A \, | \, \s \sim \eta	\}$.
\end{itemize}

\section{Anisotropic Ising model}
\label{S2}

Recall the  setting from Section \ref{resultsecA} and that the
Hamiltonian is given  by
\be{HA}
\mathrm{H}_\mathrm{A}(\s) = - \frac{J_H}{2} \sum_{(x,y) \in \L_H^\star} \s(x)\s(y)
- \frac{J_V}{2} \sum_{(x,y) \in \L_V^\star} \s(x)\s(y)
- \frac{h}{2} \sum_{x \in \L} \s(x),
\ee
where $\s \in S$, $ J_H,J_V,h>0 $, $ \L_H^\star $ is the set of \textit{unordered horizontal nearest-neighbor bonds} in $ \L $ and $ \L_V^\star$ is the set of \textit{unordered vertical nearest-neighbor bonds} in $ \L $.
 
Using the geometric representation of $\s$, one can rewrite $ \mathrm{H}_\mathrm{A}(\s) $ as
\be{HgeomA}
\mathrm{H}_\mathrm{A}(\s) = 
\mathrm{H}_\mathrm{A}(\boxminus) 
-h |\s|
+J_H |\partial_V(\s)|
+J_V |\partial_H(\s)|,
\ee
where $  |\partial_V(\s)|  $ is the length of the vertical part of $\partial (\s)$ and $  |\partial_H(\s) | $ is the length of the horizontal part of $\partial (\s)$. In Figure \ref{Example} we observe that  $ | \partial_V(\s)| = 34  $  and $  |\partial_H(\s)|= 40  $.

Recall that the critical length in this model is given by
\begin{align}\label{CritlengthA}
L_V^\star = \left\lceil   \frac{2J_V}{h} \right\rceil.
\end{align}
We make the following assumptions in this chapter.
\bass{AssA}
\begin{enumerate}[a)]
	\item $ J_H > J_V $,
	\item $ 2J_V > h $,
	\item $ \frac{2J_V	}{h} \notin \N $, 
	\item $ |\L| > \left( \max\{\frac{2 J_H}{hL_V^\star - 2J_V}, \frac{2J_H (L_V^\star - 1)}{ 2J_V - h(L_V^\star -1)} + L_V^\star \}  \right)^2 $. 
\end{enumerate}
\eass
By symmetry, Assumption \ref{AssA} a)  could be chosen the other way around. 
Assumption \ref{AssA} b) implies that the dynamics prefers aligned  neighboring spins to $(+1)$--spins.
This is essential to obtain the metastable  behavior of the system. 
Indeed, if $ 2J_V \leq h $, then $ L_V^\star = 1 $ and therefore, each configuration with a single $(+1)$--spin somewhere in $ \L $ is a critical configuration of the system.
It  follows from Assumption \ref{AssA} c) that
\begin{align}\label{InequA}
(L_V^\star -1)h < 2J_V < L_V^\star  h.
\end{align}
In Section \ref{Step2A} and Section \ref{Step4A} the importance of \eqref{InequA} will become clear. 
 Assumption \ref{AssA} d) is made to avoid certain degenerate situations. 
 For instance, if $ |\L| $ is small enough, all optimal paths between $ \boxminus $ and $ \boxplus $ contain a configuration, which consists of a single rectangle, where one side wraps around the torus and the other side is of length strictly smaller than $ L_V^\star -1$. 
For more details see \eqref{EqAssA01} or the proof of Lemma \ref{Gate1A}.
 Moreover, d) ensures that  the torus is large enough to contain at least a critical droplet. 

Recall the definition of $ R(l_1,l_2) $ from Section \ref{SecGeo}.
Before stating the main result of this chapter, we need the following definition. 

\begin{figure}[htbp]
	\centering
	\includegraphics[scale=0.6]{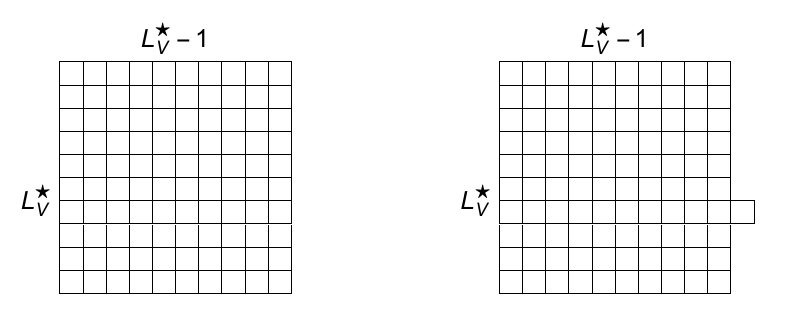}	
	\vspace{-0.7cm}
	\caption{{\footnotesize Configurations in $ \cP^ \star$ and $ \cC^ \star$}}
	\label{PandC_A}
\end{figure}

\bd{defiA}
We denote by $R (L_V^\star  -1 , L_V^\star) ^\mathrm{1pr}$  the set of all configurations consisting only of a rectangle from $ R (L_V^\star  -1 , L_V^\star) $ and with an additional protuberance attached to one of its longer sides.
The right droplet in Figure \ref{PandC_A} provides an example.

Moreover, we denote by 
$R (L_V^\star  -1 , L_V^\star)^\mathrm{2pr}$  the set of all configurations that are obtained from a configuration in $R (L_V^\star  -1 , L_V^\star)^\mathrm{1pr}$ by adding a second $ (+1) $--spin, which is attached to the rectangle and adjacent to the protuberance. 
\ed

We now formulate the main result of this chapter.
\bt{ThmA}
Under Assumption \ref{AssA}, the pair $ (\boxminus,\boxplus) $ satisfies (H1) and (H2) so that Theorems \ref{T1}--\ref{T3} hold for the anisotropic Ising model. 
Moreover, 
\bi
\item	$ \cP^ \star =R (L_V^\star  -1 , L_V^\star) $,
\item	$ \cC^ \star =R (L_V^\star  -1 , L_V^\star) ^\mathrm{1pr} $,
\item $  \Phi(\boxminus,\boxplus) - \mathrm{H}_\mathrm{A}(\boxminus) = 2 L_V^\star (J_H+J_V) - h (1+(L_V^\star-1)L_V^\star)  =:\G^\star_\mathrm{A} =:  \mathrm{E}_\mathrm{A}^\star  - \mathrm{H}_\mathrm{A}(\boxminus) $,
\item $K^{-1}=   \frac{4(2L_V^\star-1)}{3}|\L|  $.  
\ei
\et
\bpr
The proof is divided into the Sections \ref{Step1A}--\ref{Step6A}. 
\epr

\subsection{Proof of $ \Phi(\boxminus,\boxplus) - \mathrm{H}_\mathrm{A}(\boxminus) \leq  \G^\star_\mathrm{A} $}

\label{Step1A}
It will be enough to construct a path $ \gamma_\mathrm{A}=(\gamma_\mathrm{A}(n))_{n\geq0}: \boxminus \ra \boxplus$ such that 
\be{RefPathCondA}
 \max_{\eta \in  \gamma_\mathrm{A} } \mathrm{H}_\mathrm{A}(\eta) \leq  \mathrm{H}_\mathrm{A}(\boxminus) +  \G^\star_\mathrm{A} =  \mathrm{E}_\mathrm{A}^\star .
 \ee 
This path will be called \emph{reference path}.
\medskip

\emph{Construction of $  \gamma_\mathrm{A} $.}
%
Let $  \gamma_\mathrm{A}(0)= \boxminus $. 
In the first step an arbitrary $(-1)$--spin is flipped.
Then $  \gamma_\mathrm{A} $  first passes through a sequence of squares and quasi-squares as follows.
If at some step $ i $, $ \gamma_\mathrm{A}(i)$ is a square, then a protuberance is added above the droplet.
Afterwards, this row is filled by successively flipping in this row adjacent $(-1)$--spins until the droplet has the shape of a quasi-square.
Next, a protuberance is added on the right of the droplet. 
Similarly as before,  successively, adjacent $(-1)$--spins are flipped in this column until the droplet has the shape of a square again. 
%
%
%
This procedure is stopped, when $ R (( L_V^\star  -1) \times L_V^\star) $ is reached.

Now a protuberance is added on the right of the droplet and this column is filled until $ R (( L_V^\star  -1) \times (L_V^\star+1)) $
is reached. 
This adding structure is repeated until the droplet winds around the torus.
Next, a protuberance is added above the droplet and the corresponding row is filled until this row also winds around the torus. 
This is repeated until $ \boxplus $ is reached. 
\medskip

\emph{Inequality \eqref{RefPathCondA} holds.}
Let $ k^\star  $ be such that $  \gamma_\mathrm{A}(k^\star) \in R (( L_V^\star  -1) \times L_V^\star) $.
Then $ \mathrm{H}_\mathrm{A}( \gamma_\mathrm{A}(k^\star)) = \mathrm{E}^\star_\mathrm{A} - 2 J_V + h <  \mathrm{E}^\star_\mathrm{A}$. 
If we go backwards in the path from that point on, then we will have to cut the top row of $ R (( L_V^\star  -1) \times L_V^\star) $, which has the length $ L_V^\star  -1 $.
This is an increase of the energy in each step by $ h  $ for $ (L_V^\star  -2) $ times until the top row turns into  a protuberance.
At this point the energy equals to
\begin{align}
 \mathrm{H}_\mathrm{A}( \gamma_\mathrm{A}(k^\star-(L_V^\star  -2))) = \mathrm{E}^\star_\mathrm{A} - 2 J_V + (L_V^\star  -1)h  < \mathrm{E}^\star_\mathrm{A}
\end{align}
by \eqref{InequA}.
Cutting the last protuberance decreases the energy   by $ 2 J_H -h $. 
By the same arguments, if we keep on going backwards in the path of $  \gamma_\mathrm{A} $, we will always stay below $ \mathrm{E}^\star_\mathrm{A} $, since the size of the above and right bars of the droplets will be at most $ L_V^\star  -  1$.
Hence, we get that 
\be{RefPathFirstPartA} 
\max_{i=1,\dots,k^\star}\mathrm{H}_\mathrm{A}( \gamma_\mathrm{A}(i)) < \mathrm{E}^\star_\mathrm{A} .
\ee

We now consider the remaining path of $  \gamma_\mathrm{A} $ after the step $ k^\star +2 $.
It holds that $ \mathrm{H}_\mathrm{A}( \gamma_\mathrm{A}(k^\star+2)) = \mathrm{E}^\star_\mathrm{A}- h <  \mathrm{E}^\star_\mathrm{A}$. 
While filling the right column, the energy decreases by $ h $ at every step.
After  the right column is filled, a protuberance is added on the right side and the energy increases by $ 2J_V-h $.
Again by  \eqref{InequA}, we get that
\begin{align}
\mathrm{H}_\mathrm{A}( \gamma_\mathrm{A}(k^\star+(L_V^\star  +1))) = \mathrm{E}^\star_\mathrm{A} + 2 J_V - L_V^\star  h  < \mathrm{E}^\star_\mathrm{A}.
\end{align}
Repeating this until the droplet wraps around the torus, the following energy level is reached
\begin{align}
\mathrm{E}^\star_\mathrm{A} - (h L_V^\star - 2 J_V) (\sqrt{|\L|}- L_V^\star) - h(L_V^\star-1)-2J_V L_V^\star.
\end{align}
Now we add a protuberance above the droplet and the energy increases by $ 2J_H - h  $. 
Assumption \ref{AssA} d)  and \eqref{InequA} imply that
 \begin{align}\label{EqAssA01}
\mathrm{E}^\star_\mathrm{A} &- (h L_V^\star - 2 J_V) (\sqrt{|\L|}- L_V^\star) - h(L_V^\star-1)- 2J_V L_V^\star + 2J_H - h \nonumber \\
&\leq \mathrm{E}^\star_\mathrm{A} + (h L_V^\star - 2 J_V) L_V^\star - h L_V^\star-2J_V L_V^\star   \\
&< \mathrm{E}^\star_\mathrm{A} -2J_V L_V^\star.\nonumber 
\end{align}
Filling this row, decreases the energy by $ (\sqrt{|\L|}- 1)h +2J_V$.
In the same way, one can show that the remaining part of the path stays below $ \mathrm{E}^\star_\mathrm{A} $. Combining this with \eqref{RefPathFirstPartA} and the fact that $ \mathrm{H}_\mathrm{A}( \gamma_\mathrm{A}(k^\star+1)) = \mathrm{E}^\star_\mathrm{A} $, we infer \eqref{RefPathCondA}.

\subsection{Proof of $ \Phi(\boxminus,\boxplus) - \mathrm{H}_\mathrm{A}(\boxminus) \geq  \G^\star _\mathrm{A}$}
 \label{Sec2Geq}
\label{Step2A}
It suffices to show that every optimal path from $ \boxminus $ to $ \boxplus $ has to pass through $ R (L_V^\star  -1 , L_V^\star) ^\mathrm{1pr} $.
We first list a few observations. Recall the definition of \emph{local minimum} from Section \ref{SecGeo}. 
\bl{LocMinA}
Let $ \s \in S $ be a local minimum of $ \mathrm{H}_\mathrm{A} $.
Then $ \s $ is a union of isolated rectangles. 
\el
\bpr
Suppose that $\s$ has a connected component $ \s_1 $ that is not a rectangle. 
Consider a connected component $ \gamma_1 $ of $ R(\s_1) \cap (\Z^2 \setminus \s_1) $.
Let $ l_1 $ be the maximal component of the boundary of $ \gamma_1 $ that does not belong to the boundary of $ R(\s_1) $. An example would be:
\begin{center}
	\includegraphics[scale=0.4]{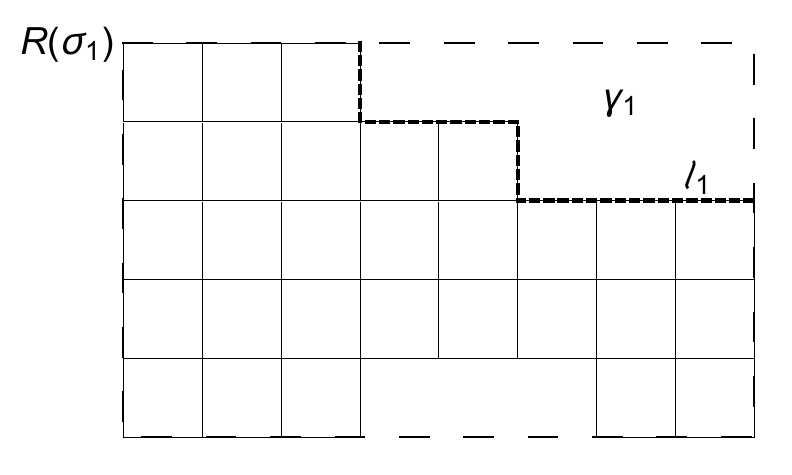}
\end{center}
Then, since $ \s_1 $  is connected and $ l_1 $ lies inside $ R(\s_1) $, $ l_1 $ has both a horizontal part and a vertical part. 
Let $ x \in \gamma_1 $ be a site, whose boundary intersects both a horizontal part and a vertical part of $ l_1 $.
In particular, $ \s(x)=-1 $ and $ x $ has at least two nearest-neighbor $ (+1) $--spins.
%
It is easy to see that $ \s^x $ has strictly lower energy than $ \s $.
\epr
%
\bc{MinRectA}
Assume that $ \s \in S $ consists  of a unique cluster. 
Then
\begin{align}
 \mathrm{H}_\mathrm{A}(\s) \geq  \mathrm{H}_\mathrm{A}(R(\s)),
\end{align}
and equality holds if and only if $ \s=R(\s) $.
\ec

We first show that every optimal path has to cross $ R (L_V^\star  -1 , L_V^\star)$.
\bl{Gate1A}
Let $ \gamma \in (\boxminus , \boxplus )_{\mathrm{opt}} $. 
Then $ \gamma $  has to cross $ R (L_V^\star  -1 , L_V^\star)$.
\el
\bpr
Assume the contrary, i.e.\ $ \gamma \cap R (L_V^\star  -1 , L_V^\star) = \emptyset$.
Let us first assume that throughout its whole path $ \gamma $ consists of a unique cluster. 
On its way to $ \boxplus $, $ \gamma  $ has to cross a configuration, whose rectangular envelope has both  horizontal and vertical length greater or equal to $ L_V^\star $.
Let 
\begin{align}\label{EqBartA}
 \bar{t} = \min\{ l\geq 0 \, | \, P_HR(\gamma(l))  , P_VR(\gamma(l)) \geq  	L_V^\star  	\} .
\end{align}
Since $ \gamma $ is assumed to consist of a unique cluster, we have that either $ P_HR(\gamma(\bar{t}-1)) =  	L_V^\star  -1$ holds or $ P_VR(\gamma(\bar{t}-1)) =  	L_V^\star -1 $.
In the following we analyze both cases and show that the assumption $ \gamma \cap R (L_V^\star  -1 , L_V^\star) = \emptyset$ leads to a contradiction.
\medskip

  \emph{Case 1}. [$ P_VR(\gamma(\bar{t}-1)) =  	L_V^\star -1 $].\\
From the definition of $ \bar{t} $, it is clear that $ R(\gamma(\bar{t}-1))\in R((L_V^\star +m)\times (L_V^\star -1 )) $ for some $ m\geq  0 $.

  \emph{Case 1.1}. [$ m=0$].\\
By hypothesis, $ \gamma  $ does not cross $ R ( L_V^\star \times ( L_V^\star  -1) ) $.
Hence,  Corollary \ref{MinRectA} yields that 
\begin{align}
\mathrm{H}_\mathrm{A}(\gamma(\bar{t}-1)) > \mathrm{H}_\mathrm{A} (L_V^\star \times ( L_V^\star  -1) ) =\mathrm{H}_\mathrm{A}(\boxminus) +  \G^\star _\mathrm{A} -2J_H +h=   \mathrm{E}_\mathrm{A}^\star -2J_H +h.
\end{align}
The minimal increase of energy to enlarge  the vertical length of the rectangular envelope of a configuration  is $ 2J_H -h  $.
Hence,
\begin{align}
\mathrm{H}_\mathrm{A}(\gamma(\bar{t})) \geq \mathrm{H}_\mathrm{A}(\gamma(\bar{t}-1)) + 2J_H -h>    \mathrm{E}_\mathrm{A}^\star.
\end{align}
This contradicts $ \gamma \in (\boxminus , \boxplus )_{\mathrm{opt}} $, since we already know from Section \ref{Step1A} that $ \Phi(\boxminus,\boxplus) \leq  E^\star _\mathrm{A}$.

  \emph{Case 1.2}. [$ m\in [1,\sqrt{|\L|} - L_V^\star)$].\\
Again, by Corollary \ref{MinRectA}  we have that 
\begin{align}
\mathrm{H}_\mathrm{A}(\gamma(\bar{t}-1)) &\geq \mathrm{H}_\mathrm{A} (  (L_V^\star +m) \times ( L_V^\star  -1)) \nonumber \\
&=\mathrm{H}_\mathrm{A} ( L_V^\star \times ( L_V^\star  -1)) + m (2J_V- h (L_V^\star -1))\\
&>  \mathrm{E}_\mathrm{A}^\star -2J_H +h,\nonumber 
\end{align}
where we used inequality \eqref{InequA}  in the last step.
As before, this leads to a contradiction, since
\begin{align}
\mathrm{H}_\mathrm{A}(\gamma(\bar{t})) \geq \mathrm{H}_\mathrm{A}(\gamma(\bar{t}-1)) + 2J_H -h>    \mathrm{E}_\mathrm{A}^\star.
\end{align}

  \emph{Case 1.3}. [$ m=\sqrt{|\L|} - L_V^\star$].\\
In this case, $ \gamma(\bar{t}-1) $ wraps around the torus. Using Assumption \ref{AssA} d), we infer that 
\begin{align}
\mathrm{H}_\mathrm{A}&(\gamma(\bar{t}-1)) \geq \mathrm{H}_\mathrm{A} (  \sqrt{|\L|} \times ( L_V^\star  -1)) \nonumber  \\
&=\mathrm{H}_\mathrm{A} (L_V^\star \times ( L_V^\star  -1) ) +(\sqrt{|\L|} - L_V^\star ) (2J_V - h(L_V^\star -1)) - 2J_V(L_V^\star -1)\\
&> \mathrm{H}_\mathrm{A} (L_V^\star \times ( L_V^\star  -1) )
=  \mathrm{E}_\mathrm{A}^\star - 2J_H +h. \nonumber 
\end{align}
Finally, 
\begin{align}
\mathrm{H}_\mathrm{A}(\gamma(\bar{t})) \geq \mathrm{H}_\mathrm{A}(\gamma(\bar{t}-1)) + 2J_H -h>    \mathrm{E}_\mathrm{A}^\star,
\end{align}
which is a contradiction.
\medskip

  \emph{Case 2}. [$ P_HR(\gamma(\bar{t}-1)) =  	L_V^\star -1 $].\\
Here we have that $ R(\gamma(\bar{t}-1))\in R( (L_V^\star -1 )\times (L_V^\star +m')) $ for some $ m'\geq  0 $.

  \emph{Case 2.1}. [$ m'=0$].\\
Since $ \gamma  $ does not cross $ R ( ( L_V^\star  -1) \times  L_V^\star) $, we have by Corollary \ref{MinRectA} that 
\begin{align}
\mathrm{H}_\mathrm{A}(\gamma(\bar{t}-1)) > \mathrm{H}_\mathrm{A} (( L_V^\star  -1) \times  L_V^\star) =   \mathrm{E}_\mathrm{A}^\star -2J_V +h.
\end{align}
The minimal increase of energy to enlarge  the horizontal length of the rectangular envelope of a configuration  is $ 2J_V -h  $.
Hence,
\begin{align}
\mathrm{H}_\mathrm{A}(\gamma(\bar{t})) \geq \mathrm{H}_\mathrm{A}(\gamma(\bar{t}-1)) + 2J_V -h>    \mathrm{E}_\mathrm{A}^\star.
\end{align}
As before, this contradicts $ \gamma \in (\boxminus , \boxplus )_{\mathrm{opt}} $.

 \emph{Case 2.2}. [$ m'\in [1,\sqrt{|\L|} - L_V^\star)$].\\
This case also leads to a contradiction, since
\begin{align}
\mathrm{H}_\mathrm{A}(\gamma(\bar{t})) &\geq \mathrm{H}_\mathrm{A}(\gamma(\bar{t}-1)) + 2J_V -h \geq \mathrm{H}_\mathrm{A} (    ( L_V^\star  -1)\times (L_V^\star +m')) + 2J_V -h \nonumber \\
&=\mathrm{H}_\mathrm{A} ( ( L_V^\star  -1)\times L_V^\star ) + m' (2J_H- h (L_V^\star -1)) + 2J_V -h\\
&>  \mathrm{E}_\mathrm{A}^\star,\nonumber 
\end{align}
where we have used inequality \eqref{InequA} and Assumption \ref{AssA} a) in the last step.

 \emph{Case 2.3}. [$ m'=\sqrt{|\L|} - L_V^\star$].\\
Using Assumption \ref{AssA} d), we infer that 
\begin{align}
\mathrm{H}_\mathrm{A}&(\gamma(\bar{t}-1)) \geq \mathrm{H}_\mathrm{A} (   ( L_V^\star  -1)  \times \sqrt{|\L|} ) \nonumber \\
&=\mathrm{H}_\mathrm{A} (( L_V^\star  -1) \times  L_V^\star)  +(\sqrt{|\L|} - L_V^\star ) (2J_H - h(L_V^\star -1)) - 2J_H(L_V^\star -1)\\
&> \mathrm{H}_\mathrm{A} (( L_V^\star  -1) \times  L_V^\star)
=  \mathrm{E}_\mathrm{A}^\star - 2J_V +h.\nonumber 
\end{align}
Finally, 
\begin{align}
\mathrm{H}_\mathrm{A}(\gamma(\bar{t})) \geq \mathrm{H}_\mathrm{A}(\gamma(\bar{t}-1)) + 2J_V -h>    \mathrm{E}_\mathrm{A}^\star,
\end{align}
which is a contradiction.
\medskip

Now suppose that $ \gamma 	 $ can consist of several clusters, i.e.\ 
at each step $ j \in \N $, $ \gamma(j) $ consists of $ n_j \in \N $ clusters, which are denoted by $ \gamma^1(j), \dots, \gamma^{n_j}(j) 	 $.
The proof follows from similar arguments as in the first part of the proof of this lemma. 
Thus, we 
only provide the main arguments and omit the details.

Using formula \eqref{HgeomA} and Corollary \ref{MinRectA}, we infer that for all $ j \in \N $, 
\begin{align}\label{TwoCluster00Eq}
\begin{split}
\mathrm{H}_\mathrm{A}( \,\gamma(j) \,)&=  \sum_{k=1}^{n_j} \mathrm{H}_\mathrm{A}( \,\gamma^k(j) \,) - (n_j-1) \, 
\mathrm{H}_\mathrm{A}(\boxminus)  \\
&\geq 
\sum_{k=1}^{n_j} \mathrm{H}_\mathrm{A}( \,R(\gamma^k(j) )\,) - (n_j-1) \, 
\mathrm{H}_\mathrm{A}(\boxminus).
\end{split}
\end{align}
For all $ j \in \N $ and $ k\leq n_j $, set $ \ell_V^k(j) = P_VR(\gamma^k(j))$ and $\ell_H^k(j)= P_HR(\gamma^k(j)) $.
Then, 
\begin{align}\label{TwoCluster01Eq}
\begin{split}
\sum_{k=1}^{n_j} &\mathrm{H}_\mathrm{A}( \,R(\gamma^k(j) )\,) - (n_j-1) \, 
\mathrm{H}_\mathrm{A}(\boxminus)\\
&=\mathrm{H}_\mathrm{A}(\boxminus) + 2 J_H \sum_{k=1}^{n_j}   \ell_V^k(j)  + 2J_V \sum_{k=1}^{n_j}    \ell_H^k(j) - h \sum_{k=1}^{n_j} \ell_V^k(j)  \ell_H^k(j).
\end{split}
\end{align}
Set $ \ell_V(j) = \sum_{k=1}^{n_j} \ell_V^k(j)$ and $\ell_H(j)=\sum_{k=1}^{n_j} \ell_H^k(j) $ and 
define 
\begin{align} \label{TwoCluster0Eq}
\tilde{t} = \min \left\{	j\in \N \ | \  \ell_H (j), \ell_V(j) \geq L_V^\star		\right\}. 
\end{align}
We have  that either $  \ell_V (\tilde{t}-1) =  	L_V^\star  -1$ holds or $  \ell_H (\tilde{t}-1) =  	L_V^\star -1 $.
We only treat the case when $  \ell_H (\tilde{t}-1) =  	L_V^\star  -1$, since the other case is 
a straightforward combination of Case 1 above and the following arguments. 

By the definition of $ \tilde{t} $, we have that $  \ell_V (\tilde{t}-1) =  	L_V^\star  + \tilde{m}$ for some $ \tilde{m}\geq -1 $. (If $ n_{\tilde{t}-1} \geq n_{\tilde{t}}  $, then $ \tilde{m}\geq 0 $, and if $ n_{\tilde{t}-1} < n_{\tilde{t}}  $, then $ \tilde{m}\geq -1 $).
From \eqref{TwoCluster00Eq} and \eqref{TwoCluster01Eq}, we infer that 
\begin{align}\label{TwoCluster02Eq}
\nonumber
\mathrm{H}_\mathrm{A}( \,\gamma(\tilde{t}-1) \,)&\geq 
\mathrm{H}_\mathrm{A}(\boxminus) + 2 J_H \,  (	L_V^\star  + \tilde{m})  + 2J_V     (L_V^\star  -1) - h \sum_{k=1}^{n_{\tilde{t}-1}} \ell_V^k(\tilde{t}-1)  \ell_H^k(\tilde{t}-1), \text{ and } \\
\mathrm{H}_\mathrm{A}( \,\gamma(\tilde{t}) \,)&\geq 
\mathrm{H}_\mathrm{A}( \,\gamma(\tilde{t}-1) \,)  + 2J_V -h\\
&\geq 
\mathrm{H}_\mathrm{A}(\boxminus) + 2 J_H  \,  (	L_V^\star  + \tilde{m})  + 2J_V     L_V^\star - h \left(\sum_{k=1}^{n_{\tilde{t}-1}} \ell_V^k(\tilde{t}-1)  \ell_H^k(\tilde{t}-1) + 1\right).
\nonumber
\end{align}
Notice the following estimate 
\begin{align}\label{TwoCluster03Eq}
\sum_{k=1}^{n_{\tilde{t}-1}} \ell_V^k(\tilde{t}-1)  \ell_H^k(\tilde{t}-1)  \leq   \sum_{p=1}^{n_{\tilde{t}-1}} \ell_V^p(\tilde{t}-1) \sum_{k=1}^{n_{\tilde{t}-1}}  \ell_H^k(\tilde{t}-1)= (	L_V^\star  + \tilde{m})    ( L_V^\star-1),
\end{align} 
where the inequality is strict whenever $ n_{\tilde{t}-1} > 1 $.
We now have to show  that all possible values for $ \tilde{m} $ and the  hypothesis that $ \gamma \cap R (L_V^\star  -1 , L_V^\star) = \emptyset$ yield to the fact that $ \mathrm{H}_\mathrm{A}( \,\gamma(\tilde{t}) \,) > \mathrm{E}_\mathrm{A}^\star $, which is a contradiction.
However, using  \eqref{TwoCluster02Eq} and      \eqref{TwoCluster03Eq}, we can proceed
 as in Case 2.1--Case 2.3 above.
The details are straightforward adaptations   and are therefore omitted. 
This concludes the proof of this lemma. 
\epr

The following lemma  concludes the proof of $ \Phi(\boxminus,\boxplus) - \mathrm{H}_\mathrm{A}(\boxminus) \geq  \G^\star _\mathrm{A}$.
\bl{Gate2A}
Let $ \gamma\in (\boxminus , \boxplus )_{\mathrm{opt}}  $.
In order to cross a  configuration whose rectangular envelope has both vertical and horizontal length greater or equal to $ L_V^\star $, $ \gamma $ has to pass through $ R (L_V^\star  -1 , L_V^\star)$ and $ R (L_V^\star  -1 , L_V^\star)^\mathrm{1pr} $.
In particular, each optimal path between $ \boxminus $ and  $ \boxplus $ has to cross $ R (L_V^\star  -1 , L_V^\star)^\mathrm{1pr} $.
\el
\bpr
Consider the time step $ \bar{t} $ defined in \eqref{EqBartA}. 
In the proof of Lemma \ref{Gate1A} we have seen that necessarily $ \gamma(\bar{t}-1) \in R (L_V^\star  -1 , L_V^\star)$. 
Note that $ \min\{P_V(\gamma(\bar{t}-1)),P_H(\gamma (\bar{t}-1))\}= L_V^\star-1 $ and that 
$ P_VR(\gamma(\bar{t})),P_HR(\gamma(\bar{t})) \geq L_V^\star  $.
Therefore, 
$ \gamma(\bar{t}) $ must be obtained from $ \gamma(\bar{t}-1) $
by adding a protuberance at a longer side of the rectangle $ \gamma(\bar{t}-1) $. 
This implies that $ \gamma(\bar{t}) $ needs to belong to $ R (L_V^\star  -1 , L_V^\star)^\mathrm{1pr} $.
\epr

\subsection{Identification of $ \cP^ \star $ and	$ \cC^ \star $}

\label{Step3A}

From Section \ref{Step1A}, we get that $ R (L_V^\star  -1 , L_V^\star) \subset  \cP^ \star $.
Now let $ \s \in \cP^ \star  $ and $ x \in \L $ be such that $ \s^x \in \cC^ \star $.
If follows from the definition of $ \cP^ \star $ and $ \cC^ \star $ that there exists $ \gamma \in (\boxminus,\boxplus)_{\mathrm{opt}}$ and $ \ell \in \N $ such that 
\begin{enumerate}[(i)]
		\item $ \gamma(\ell)= \s $ and $ \gamma(\ell +1)=\s^x$,
		\item $ \mathrm{H}_\mathrm{A}(\gamma(k)) < \mathrm{E}^\star_\mathrm{A}  $ for all $  k \in \{ 0,\dots,\ell	 	\}  $,
		\item $ \Phi(\boxminus,\gamma(k) ) \geq  \Phi(\gamma(k) , \boxplus )$  for all $  k \geq \ell + 1  $.
\end{enumerate}
By Lemma \ref{Gate2A}, (ii) implies that  $\min(P_HR(\s),P_VR(\s))  \leq  L_V^\star  -1 $, since otherwise the energy level $ \mathrm{E}^\star_\mathrm{A} $ would have been reached.
There are two possible cases.

  \emph{Case 1}. [$ P_HR(\s^x), P_VR(\s^x)  \geq  L_V^\star $].\\
Lemma  \ref{Gate2A} implies that we necessarily have that $ \s \in R (L_V^\star  -1 , L_V^\star) $ and  $ \s^x \in R (L_V^\star  -1 , L_V^\star)^\mathrm{1pr} $.

 \emph{Case 2}. [$ \min( P_HR(\s^x),P_VR(\s^x))  \leq  L_V^\star -1$].\\
Also by  Lemma \ref{Gate2A}, there must exist some $ k^\star \geq \ell + 2  $ such that $ \gamma(k^\star ) \in R (L_V^\star  -1 , L_V^\star) $. 
But this contradicts (iii), since $ \Phi(\boxminus,\gamma(k^\star ) ) < \Phi(\gamma(k^\star) , \boxplus )=  \mathrm{E}_\mathrm{A}^\star $.

Hence, only Case 1 can hold true.
We conclude that $ \cP^ \star=R (L_V^\star  -1 , L_V^\star) $ and	$ \cC^ \star= R (L_V^\star  -1 , L_V^\star)^\mathrm{1pr} $.

\subsection{Verification of (H1)}

\label{Step4A}
Obviously, $ S_{\mathrm{stab}}= \{ \boxplus	\} $, since $ \boxplus $ minimizes all three sums in \eqref{HA}.
It remains to show that $ S_{\mathrm{meta}}= \{ \boxminus	\}. $

Let $ \s \in S\setminus\{\boxminus,\boxplus  \} $.
We have to show that $V_\s < \G^\star_\mathrm{A}$, i.e.\  there exists $ \s' \in S $ such that $ \mathrm{H}_\mathrm{A}(\s' ) < \mathrm{H}_\mathrm{A}(\s)$ and $ \Phi ( \s,\s') - \mathrm{H}_\mathrm{A}(\s) < \G_\mathrm{A}^\star $. 
There are four possible cases.

 \emph{Case 1}. [$ \s $ contains a cluster, which is not a rectangle].\\
Lemma \ref{LocMinA} implies that $ \s $ is not a local minimum, i.e.\ there exists $ x \in \L $ such that $ \mathrm{H}_\mathrm{A}(\s^x ) < \mathrm{H}_\mathrm{A}(\s)$.
Moreover, $ \Phi ( \s,\s^x) - \mathrm{H}_\mathrm{A}(\s) =0 < \G_\mathrm{A}^\star $.

 \emph{Case 2}. [$ \s $ contains a cluster, which is a  rectangle $  R= l_1 \times l_2 $ with $ l_2 \geq L_V^\star $ and $ l_1<\sqrt{|\L|} $].\\
Let $ \s' $ be obtained from $ \s $ by attaching on the right of $ R $ a new column of length $ l_2 $.
Then,
\begin{align}
	\mathrm{H}_\mathrm{A}(\s' ) &\leq  \mathrm{H}_\mathrm{A}(\s ) + 2J_V - l_2 h \leq \mathrm{H}_\mathrm{A}(\s ) + 2J_V - L_V^\star  h< \mathrm{H}_\mathrm{A}(\s), \quad \text{and} \nonumber \\
	\Phi ( \s,\s') - \mathrm{H}_\mathrm{A}(\s) &= 2J_V - h < \G_\mathrm{A}^\star.
\end{align}
 
 \emph{Case 3}. [$ \s $ contains a cluster, which is a  rectangle $  R= l_1 \times l_2 $ with $ l_2 < L_V^\star $ and $ l_1<\sqrt{|\L|} $].\\
Let $ \s' $ be obtained from $ \s $ by cutting the right column of $ R $.
Then,
\begin{align}
\mathrm{H}_\mathrm{A}(\s' ) &=  \mathrm{H}_\mathrm{A}(\s ) - 2J_V + l_2 h \leq \mathrm{H}_\mathrm{A}(\s ) - 2J_V +  (L_V^\star-1) h < \mathrm{H}_\mathrm{A}(\s), \quad \text{and}\nonumber  \\
\Phi ( \s,\s') - \mathrm{H}_\mathrm{A}(\s) &= (l_2-1)h < \G_\mathrm{A}^\star.
\end{align}

 \emph{Case 4}. [$ \s $ contains a cluster, which is a  rectangle $  R= l_1 \times l_2 $ with $ l_1=\sqrt{|\L|} $].\\
Let $ \s' $ be obtained from $ \s $ by attaching above $ R $ a row that also wraps around the torus.
Then, by Assumption \ref{AssA} d),
\begin{align}
\mathrm{H}_\mathrm{A}(\s' ) &=  \mathrm{H}_\mathrm{A}(\s ) + 2J_H - l_1 h < \mathrm{H}_\mathrm{A}(\s), \quad \text{and} \nonumber \\
\Phi ( \s,\s') - \mathrm{H}_\mathrm{A}(\s) &= 2J_H - h < \G_\mathrm{A}^\star.
\end{align}

We conclude that $ S_{\mathrm{meta}}= \{ \boxminus	\}. $

\subsection{Verification of (H2)}

\label{Step5A}

Obviously, $ |\{ \s \in \cP^\star \,  | \,  \s \sim \s' \}| = 1 $ for all $ \s' \in \cC^\star$.
Therefore, (H2) holds.

\subsection{Computation of $K$}

\label{Step6A}

The starting point for the computation of $ K$ is the variational formula  \eqref{EqPrefactor}. 
Recall the definitions of $ \partial^+ A $ and  $ A^+  $ for a subset $ A \subset S$ (see Section \ref{SecGeo}).

\emph{Lower bound.}
Since the sum in \eqref{EqPrefactor} has only non-negative summands, we can bound $ K^{-1} $ from below by
	\begin{align}
\frac{1}{K} \geq \min_{C_1,\dots,C_I \in [0,1]}
\min_{\substack{ h:S^\star \ra [0,1] \\
		\smallrestr{h}{S_\boxminus} =1, \smallrestr{h}{S_\boxplus} =0, \smallrestr{h}{S_i} =C_i \, \forall i }}
\frac{1}{2} \sum_{\eta,\eta' \in (\cC^\star)^+ } \mathbbm{1}_{\{	\eta \sim \eta'	 \}} [h(\eta) - h(\eta')]^2.
\end{align}
Obviously,  $ \partial^+ \cC^\star \cap S^\star = R (L_V^\star  -1 , L_V^\star)  \cup R (L_V^\star  -1 , L_V^\star)^\mathrm{2pr}  $.
Moreover, similar computations as in Section \ref{Step1A} show that $  R (L_V^\star  -1 , L_V^\star)  \subset  S_\boxminus  $ and $ R (L_V^\star  -1 , L_V^\star)^\mathrm{2pr} \subset S_\boxplus  $. This leads to
\begin{align}
\frac{1}{K} &\geq 
\min_{ h:\cC^\star \ra [0,1]}
 \sum_{\eta \in \cC^\star } \left( \sum_{\eta' \in R (L_V^\star  -1 , L_V^\star), \eta'\sim \eta  }	[1- h(\eta) ]^2	+ \sum_{\eta' \in R (L_V^\star  -1 , L_V^\star)^\mathrm{2pr}, \eta'\sim \eta  }	h(\eta)^2			\right) \nonumber  \\
 &=
 \sum_{\eta \in \cC^\star }  \min_{ h \in [0,1]} \Big( |R (L_V^\star  -1 , L_V^\star)\sim \eta| \ 	[1- h ]^2	+ 
 |R (L_V^\star  -1 , L_V^\star)^\mathrm{2pr}\sim \eta| \  h^2			\Big) \\
 &=
 \sum_{\eta \in \cC^\star } \frac
 {|R (L_V^\star  -1 , L_V^\star)\sim \eta| 	\cdot 
 	|R (L_V^\star  -1 , L_V^\star)^\mathrm{2pr}\sim \eta| }
 {|R (L_V^\star  -1 , L_V^\star)\sim \eta| 	+ 
 	|R (L_V^\star  -1 , L_V^\star)^\mathrm{2pr}\sim \eta| }.\nonumber 
\end{align}
For all $ \eta \in \cC^\star  $ we have that $ |R (L_V^\star  -1 , L_V^\star)\sim \eta| = 1 $.
If the protuberance in $ \eta  $ is attached at a corner of $ ( L_V^\star  -1) \times L_V^\star $, then $ |R (L_V^\star  -1 , L_V^\star)^\mathrm{2pr}\sim \eta|=1 $, otherwise $ |R (L_V^\star  -1 , L_V^\star)^\mathrm{2pr}\sim \eta|=2 $.
Taking into account that there are $ |\L| $ possible locations for each shape of a critical droplet and 2 possible rotations, we obtain that
\begin{align}
\frac{1}{K} &\geq
 \Big( 2(L_V^\star-2) \frac{2} {3 } + 4 \frac {1 } {2 } \Big) 2|\L|= \frac{4(2L_V^\star-1)}{3}|\L|.
\end{align}

\emph{Upper bound.}
Define 
\begin{align}
\begin{split}
\mathscr{S}^- =  \{\s\in S^\star\, &| \,  
\min(P_HR(\eta),P_VR(\eta))  \leq  L_V^\star  -1 \text{ for all clusters $ \eta $ of $ \s $} \}, \text{ and} \\
\mathscr{S}^+ =  \{\s\in S^\star\, &| \, \text{ there exists a cluster $ \eta  $ of $ \s $ such that } 
P_HR(\eta),P_VR(\eta)  \geq  L_V^\star  \} .
\end{split}
\end{align} 
Note that $ S^\star = \mathscr{S}^-\cup \mathscr{S}^+$, $ \mathscr{S}^- \cap \mathscr{S}^+ = \emptyset $, $ \cP^\star \subset \mathscr{S}^- $ and $ \cC^\star \subset \mathscr{S}^+ $.
Using the same arguments as in Lemma \ref{Gate1A},
we can show the following fact for transitions between $ \mathscr{S}^- $ and $ \mathscr{S}^+ $.
\bl{UpperBoundA}
Let $ \s \in \mathscr{S}^- $ and $ \s' \in \mathscr{S}^+ $.
Then $ \s \sim \s' $ if and only if  $ \s \in \cP^\star$ and $ \s' \in \cC^\star$.
\el
\bpr
We omit the details of this proof, since they walk along the same lines as  the proof of Lemma \ref{Gate1A}.
\epr

Recall that the sets $ S_1,\dots,S_I $ are assumed to be maximal sets of communicating configurations. 
Hence, for all $ i = 1,\dots,I $, we have that either $ S_i \subset \mathscr{S}^- $ or $ S_i \subset \mathscr{S}^+ $, since $ S^\star = \mathscr{S}^-\cup \mathscr{S}^+$ and $ \mathscr{S}^- \cap \mathscr{S}^+ = \emptyset $.
For the same reason and by Section \ref{Step1A}, we have that $ S_\boxminus \subset \mathscr{S}^- $ and $ S_\boxplus \subset \mathscr{S}^+ $.
Therefore, we can estimate $ K^{-1} $ from above by restricting the minimum in \eqref{EqPrefactor} only to those functions $ h:S^\star \ra [0,1] $ such that 
\begin{align}\label{EqStep6FunctionA}
\begin{split}
&h(\eta) =1 \text{ for all } \eta \in \mathscr{S}^{-}, \text{ and }\\  
&h(\eta) =0 \text{ for all } \eta \in \mathscr{S}^{+}\setminus \cC^\star.
\end{split}
\end{align}  
The restriction to such functions is allowed, since we can choose $ C_i = 1 $ for all $ i\in \{1,\dots,I\} $ such that $ S_i \subset \mathscr{S}^- $ and $ C_i = 0 $ for all $ i\in \{1,\dots,I\} $ such that $ S_i \subset \mathscr{S}^+ $.
Thus,
\begin{align}
\begin{split}
\frac{1}{K} &\leq 
\min_{\substack{ h:S^\star \ra [0,1] \\
		\smallrestr{h}{\mathscr{S}^-} =1, \smallrestr{h}{\mathscr{S}^+\setminus \cC^\star } =0}}
\frac{1}{2} \sum_{\eta,\eta' \in S^\star } \mathbbm{1}_{\{	\eta \sim \eta'	 \}} [h(\eta) - h(\eta')]^2 \\
&=\min_{\substack{ h:S^\star \ra [0,1] \\
		\smallrestr{h}{\mathscr{S}^-} =1, \smallrestr{h}{\mathscr{S}^+\setminus \cC^\star } =0}}
\frac{1}{2} \sum_{\eta \in \mathscr{S}^-,\eta' \in \mathscr{S}^-} \mathbbm{1}_{\{	\eta \sim \eta'	 \}} [h(\eta) - h(\eta')]^2.
\end{split}
\end{align}
Using Lemma \ref{UpperBoundA}, the right-hand side is equal to 
\begin{align}
\begin{split}
&\ 
\min_{\substack{ h:(\cC^\star)^+ \ra [0,1] \\
		\smallrestr{h}{\mathscr{S}^-\cap \partial^+\cC^\star} =1, \smallrestr{h}{\mathscr{S}^+\cap \partial^+\cC^\star} =0}}
\frac{1}{2} \sum_{\eta,\eta' \in (\cC^\star)^+ } \mathbbm{1}_{\{	\eta \sim \eta'	 \}} [h(\eta) - h(\eta')]^2 \\
&=
\min_{ h:\cC^\star \ra [0,1]}
\sum_{\eta \in \cC^\star } \left( \sum_{\eta' \in R (L_V^\star  -1 , L_V^\star), \eta'\sim \eta  }	[1- h(\eta) ]^2	+ \sum_{\eta' \in R (L_V^\star  -1 , L_V^\star)^\mathrm{2pr}, \eta'\sim \eta  }	h(\eta)^2			\right) \\
&=
\frac{4(2L_V^\star-1)}{3}|\L|.
\end{split}
\end{align}
This concludes the proof.

\section{Ising model with next-nearest-neighbor attraction}
\label{S4}

In this chapter  (cf.\ Section \ref{resultsecNN})  the Hamiltonian is given by
\be{HNN}
\mathrm{H}_\mathrm{NN}(\s) = - \frac{\tilde{J}}{2} \sum_{(x,y) \in \L^\star} \s(x)\s(y)
- \frac{K}{2} \sum_{(x,y) \in \L^{\star\star}} \s(x)\s(y)
- \frac{h}{2} \sum_{x \in \L} \s(x),
\ee
where $\s \in S$, $ \tilde{J},K,h>0 $, $ \L^\star $ is the set of  \textit{unordered nearest-neighbor bonds} in $ \L $ and $ \L^{\star\star}$ is the set of \textit{unordered next-nearest-neighbor bonds} in $ \L $ (cf.\ \eqref{resultsecNNEq}).
We can rewrite $ \mathrm{H}_\mathrm{NN}(\s) $ as
\be{HgeomNN}
\mathrm{H}_\mathrm{NN}(\s) = 
\mathrm{H}_\mathrm{NN}(\boxminus) 
-h |\s|
+J |\partial(\s)|
-K |A(\s)|,
\ee
where $  J= \tilde{J}+2K $  and $  |A(\s)| $ is the number of corners (or right angles) of $\s$. 
Indeed, a unit segment of $ \partial (\s) $ breaks two next-nearest-neighbor-bonds
However, at each corner the same broken next-nearest-neighbor-bond is counted twice.
This explains the term $ -K |A(\s)| $ in \eqref{HgeomNN}.
Moreover,  in the  situation 
\begin{center}
	\includegraphics[height=3cm]{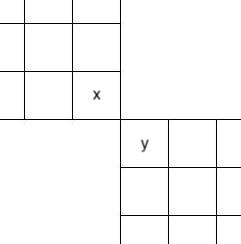}
\end{center}
we count four corners, since the bond between $ x $ and $ y $ is not broken, but we have counted it as such due to the four  unit segments surrounding this bond.

Recall that the critical lengths in this model are given by
\begin{align}\label{CritlengthsNN}
\ell^\star = \left\lceil  \frac{2K}{h} \right\rceil  \qquad \text{and} \qquad D^\star = \left\lceil   \frac{2J}{h} \right\rceil  \qquad \text{and} \qquad L^\star = D^\star - 2(\ell^\star -1).
\end{align}
We make the following assumptions  for this chapter.
\bass{AssNN}
\begin{enumerate}[a)]
	\item $ K > h $,
	\item $ \tilde{J}\geq 2K+ h $,
	\item $ \frac{2J	}{h} \notin \N $, $ \frac{2K	}{h} \notin \N $,
	 \item $ |\L| > \left( \frac{2J (D^\star - 1)}{ 2J - h(D^\star -1)} + D^\star  \right)^2 $.  
\end{enumerate}
\eass
Similarly as in Chapter 2,  a) and b) induce a hierarchy in the sense that for the system it is most important to align nearest-neighbors, then next-nearest-neighbors and then to align the spin values with the sign of the magnetic field. 
As in Chapter \ref{S2}, this assumption is essential to obtain the metastable  behavior of the system. 
Moreover, a) respects a hypothesis made in \cite{KO94} (but not every hypothesis in there). 
Assumption c) is made for non-degeneracy reasons.
Assumption d) implies that it is not profitable to  enlarge a droplet such that one side is sub-critical and the other side wraps around the torus.
This will become clear later in Lemma \ref{Gate1NN}. 
Moreover, d) ensures that  the torus is large enough to contain at least a critical droplet. 
It immediately follows from Assumption \ref{AssNN} c) that
\begin{align}\label{InequNN}
(\ell^\star -1)h < 2K < \ell ^\star  h  \qquad \text{and} \qquad (D^\star -1)h < 2J < D^\star  h.
\end{align}

\begin{figure}[htbp]
	\centering
	\includegraphics[scale=0.6]{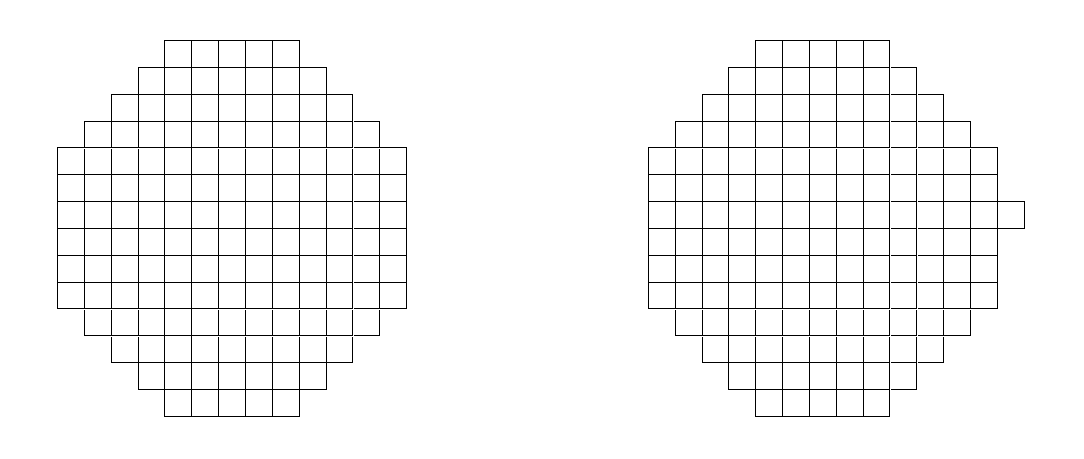}	
	\vspace{-0.5cm}
	\caption{{\footnotesize Configurations in $ \cP^ \star$ and $ \cC^ \star$ with $ \ell^\star = 5$ and $ D^\star = 14 $}}
	\label{PandC_NN}
\end{figure}

We need a few definitions that are mostly carried over from \cite{KO94}. 
\bd{defiNN}
\begin{itemize}
	\setlength\itemsep{-0.3em} 
\item
$ A\subset \Z^2 $ is called an \emph{oblique bar} if $ A=\{x_1, \dots, x_n\} $ for some $ n\in \N $ and it holds that either $ x_i= x_{i-1} + (1,1)^T  $  or $ x_i= x_{i-1} + (1,-1)^T  $ for all $2\leq i\leq n $.
\item
We say that $ \s \in S $ is an \emph{octagon of side lengths $ D_n,D_w\in \N\cap [1,\sqrt{|\L|}-1] $ and oblique edge lengths} $ \ell_{ne},\ell_{nw},\ell_{sw},\ell_{se} \in \N  $ and write $ \s \in Q(D_n,D_w; \ell_{ne},\ell_{nw},\ell_{sw},\ell_{se}) $ if the geometric representation of $ \s  $ has the following form (cf. Scheme 2.2 in \cite{KO94}).
$ \s $ is connected and inscribed in a rectangle from $R( D_n , D_w )$. 
Moreover, $ \s $ has four straight edges with endpoints $ a_i,b_i $, $ i=1,\dots,4 $ and four oblique edges that have a local staircase structure with endpoints $( b_1,a_2) $, $( b_2,a_3) $, $( b_3,a_4) $, $( b_4,a_1) $.
The \emph{lengths of its oblique edges} are defined by 
\begin{align}
&\ell_{ne} =	1+\frac 1{\sqrt{2}} |b_1-a_2|	,  &\ell_{se} =1+\frac 1{\sqrt{2}} |b_2-a_3|, \\
&\ell_{sw} =1+\frac 1{\sqrt{2}} |b_3-a_4|, &\ell_{nw} = 1+\frac 1{\sqrt{2}} |b_4-a_1|.
\end{align}
An example with $ D_n = 15 $, $ D_w = 12 $, $ \ell_{ne} = 5 $, $ \ell_{nw} = 6 $, $ \ell_{sw} = 4 $, $ \ell_{se} = 3 $  is given by
\vspace{-0.1cm}
\begin{center}
	\includegraphics[scale=0.3]{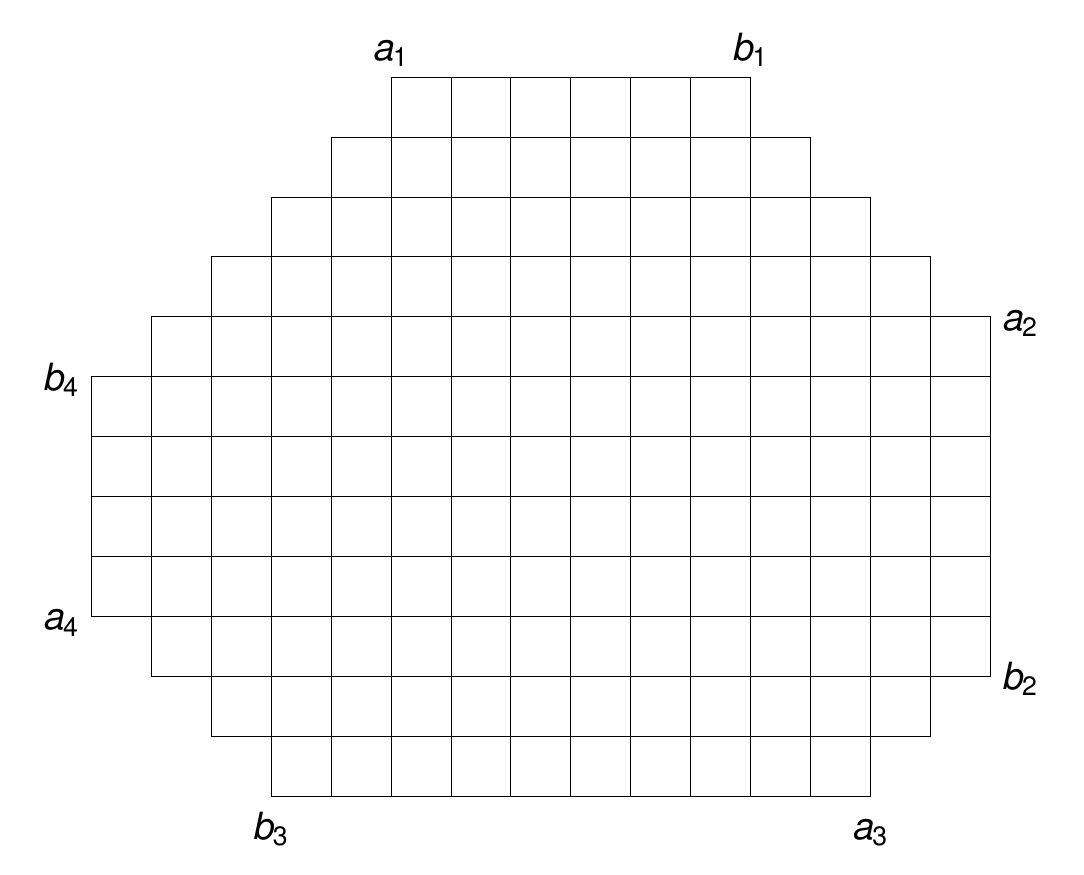}
\end{center}
\vspace{-0.5cm} 
 We often abuse the notation by  identifying   $ Q(D_n,D_w; \ell_{ne},\ell_{nw},\ell_{sw},\ell_{se}) $ with configurations from this set.
\item For $ Q \in Q(D_n,D_w; \ell_{ne},\ell_{nw},\ell_{sw},\ell_{se}) $, the upper right edge of length $ \ell_{ne} $ is called \emph{NE-edge}, the upper left edge of length $ \ell_{nw} $ \emph{NW-edge}, the down left edge of length $ \ell_{sw} $  \emph{SW-edge} and the down right edge of length $ \ell_{se} $  is called \emph{SE-edge}.
These four edges are also called \emph{oblique edges}.
The four remaining horizontal or vertical edges are called \emph{coordinate edges}. 
We call the upper coordinate edge \emph{N-edge}, the left one \emph{W-edge}, the bottom one \emph{S-edge} and the right coordinate edge \emph{E-edge}.

\vspace{-0.1cm}
\begin{center}
	\includegraphics[scale=0.3]{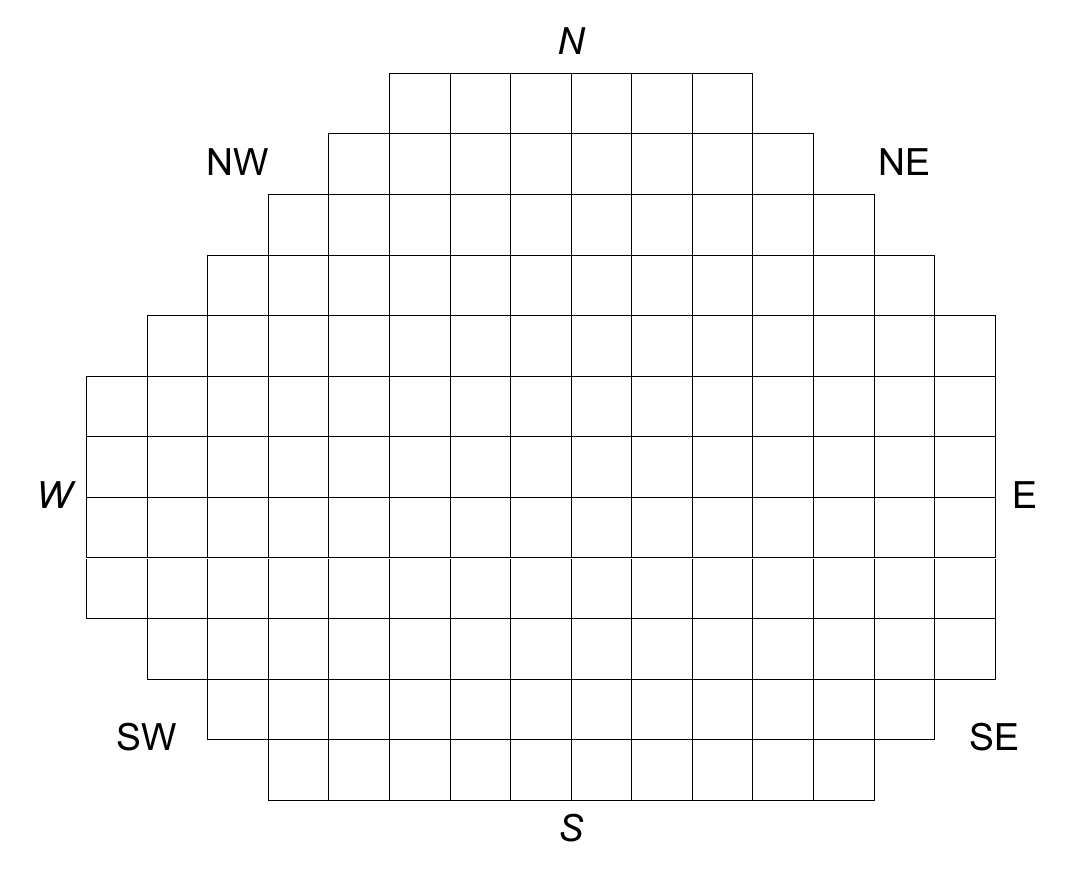}
\end{center}
\vspace{-0.5cm}


\item
$ Q(D_n,D_w; \ell_{ne},\ell_{nw},\ell_{sw},\ell_{se}) $ is called \emph{stable octagon} if each of his eight edges has length greater or equal to 2.
 


\item  We abbreviate $ Q(D_n,D_w; \ell,\ell,\ell,\ell) = Q(D_n,D_w;\ell) $ for all $ \ell \in \N $ and $ Q(D_n,D_w;\ell^\star ) = Q(D_n,D_w) $.  
 Moreover, we write $ Q(3\ell-2,3\ell-2;\ell ) = Q(\ell) $ for all $ \ell \in \N $, which corresponds to the case, where all eight edges have the same length given by $ \ell $. 

\item
$ Q(D_n,D_w;\ell)^{\mathrm{1pr}} $ denotes the set of all configurations that are obtained from a configuration in $ Q(D_n,D_w;\ell)$ by adding a protuberance somewhere at the \emph{interior} of one of its longest coordinate edges.
Here the interior of the coordinate edge contains every site of the edge except for the two sites at the end of the edge. 
The right droplet in Figure \ref{PandC_NN} provides an example.

\item
$Q(D_n,D_w;\ell)^\mathrm{2pr}$ denotes the set of all configurations that are obtained from a configuration in $Q(D_n,D_w;\ell)^\mathrm{1pr}$ by adding a second $ (+1) $--spin adjacent to the protuberance at the \emph{interior} of the coordinate edge. 
\end{itemize}
%
%
\ed

Note that the energy of an octagon $ Q \in Q(D_n,D_w;\ell_{ne},\ell_{nw},\ell_{sw},\ell_{se}) $ is given by
\begin{align}
\mathrm{H}_\mathrm{NN}(Q) = 
\mathrm{H}_\mathrm{NN}(\boxminus) 
-h D_n D_w
+2 J (D_n+ D_w)
+ \sum_{a  \in \{ne,nw,sw,se\}} F(\ell_a),
\end{align}
where $ F(\ell) = 	-K(2\ell-1) + \frac{1}{2}h(\ell-1)\ell	 $.
Now we can formulate the main result of this chapter.
\bt{ThmNN}
Under Assumption \ref{AssNN}, the pair $ (\boxminus,\boxplus) $ satisfies (H1) and (H2) so that Theorems \ref{T1}--\ref{T3} hold for the Ising model with next-nearest-neighbor attraction.
Moreover,
\bi
\item	$ \cP^ \star =Q(D^\star - 1,D^\star) $,
\item	$ \cC^ \star =Q(D^\star - 1,D^\star)^\mathrm{1pr} $,
\item $  \Phi(\boxminus,\boxplus) - \mathrm{H}_\mathrm{NN}(\boxminus) =\mathrm{H}_\mathrm{NN}(Q(D^\star - 1,D^\star))  +2J-4K-h=:\G^\star_{\mathrm{NN}} =: \mathrm{E}_{\mathrm{NN}}^\star  - \mathrm{H}_\mathrm{NN}(\boxminus) $,
\item $K^{-1}=   \frac{4(2L^\star-5)}{3}|\L|  $.  
\ei
\et
\bpr
The proof is divided into the Sections \ref{Step1NN}--\ref{Step6NN}. 
\epr

\vspace{-0.2cm}

\begin{remark}
	In the Lemmas \ref{Gate1NN} and  \ref{Gate2NN} we will see that, in order to  reach a configuration $ \s  $ with 
	$ P_HR(\s) , P_VR(\s) \geq  	D^\star   $, for any $ \gamma \in (\boxminus \ra \boxplus )_{\mathrm{opt}}$ it is essential  to pass from $Q(D^\star - 1,D^\star)$ to $Q(D^\star - 1,D^\star)^\mathrm{1pr}$.
		Hence,  if $ \gamma $ would pass through a configuration $ \eta  $ that is obtained from $Q(D^\star - 1,D^\star)$ by adding a protuberance at the shortest coordinate edge, $ \gamma $ would need to go back to $Q(D^\star - 1,D^\star)$ and pass through $Q(D^\star - 1,D^\star)^\mathrm{1pr}$ to be ``over the hill'', i.e.\ to be able to make the transition to $ \boxplus $ without reaching the energy level $ \G^\star_{\mathrm{NN}}  $ again.
	Consequently,  $ \eta \notin \cC^\star $, since such a configuration can not fulfill property \emph{3.)} of Definition \ref{Crit}, since $\Phi(\boxminus,Q(D^\star - 1,D^\star)) < \Phi(Q(D^\star - 1,D^\star),\boxplus)$ (as we will see in Section \ref{Step1NN}).
	\end{remark}

\subsection{Proof of $ \Phi(\boxminus,\boxplus) - \mathrm{H}_\mathrm{NN}(\boxminus) \leq  \G^\star_{\mathrm{NN}} $}

\label{Step1NN}
As in Section \ref{Step1A}, we need to construct a reference path $  \gamma_\mathrm{NN}: \boxminus \ra \boxplus$ such that 
\be{RefPathCondNN}
\max_{\eta \in  \gamma_\mathrm{NN} } \mathrm{H}_\mathrm{NN}(\eta) \leq  \mathrm{H}_\mathrm{NN}(\boxminus) +  \G^\star_{NN} = \mathrm{E}_{\mathrm{NN}}^\star .
\ee 

\emph{Construction of $  \gamma_\mathrm{NN} $.} 
We only sketch the construction of $  \gamma_\mathrm{NN} $, since we can rely on \cite{KO94} and we are mainly interested in the part of the path around the critical configuration.

\textbullet \ [From $ \boxminus $ to $ Q(2) $.]\\
See Scheme 5.1 of \cite{KO94}.

\textbullet \ [From $ Q(\ell) $ to $ Q(\ell+1) $ for all $ \ell= 2,\dots, \ell^\star -1 $.]\\
See Scheme 5.2 of \cite{KO94}.

\textbullet \ [From $ Q(D,D) $ to $ Q(D+1,D+1) $ for all $ D= \ell^\star ,\dots, \sqrt{|\L|}-2$.]\\
This transition is based on  Scheme 5.5 of \cite{KO94}, and it goes for example as follows. 
A $(+1)$--spin is added somewhere at the interior of the E-edge of $ Q(D,D) $. 
Afterwards, successively, adjacent $(-1)$--spins are flipped in this column  until $ Q(D+1,D; \ell^\star+1,\ell^\star,\ell^\star,\ell^\star+1) $ is reached.
Then a $(-1)$--spin is flipped at the upper end of the SE-edge. 
Now $(-1)$--spins are flipped until $ Q(D+1,D; \ell^\star+1,\ell^\star,\ell^\star,\ell^\star) $ is reached. 
Next, the same  is done at the NE-edge such that $ Q(D+1,D; \ell^\star,\ell^\star,\ell^\star,\ell^\star) = Q(D+1,D)$ is reached.
This procedure is  repeated below $ Q(D+1,D) $, i.e.\ first $(-1)$--spins are flipped at the S-edge until $ Q(D+1,D+1; \ell^\star,\ell^\star,\ell^\star+1,\ell^\star+1) $ is reached, then an oblique bar is added at the SW-edge  to reach $ Q(D+1,D+1; \ell^\star,\ell^\star,\ell^\star,\ell^\star+1) $, and finally, $(-1)$--spins are flipped at the SE-edge, until we arrive at $ Q(D+1,D+1) $.

\textbullet \ Lastly, flip all remaining $(-1)$--spins outside of $ Q(\sqrt{|\L|}-1, \sqrt{|\L|}-1) $ until $ \boxplus $ is reached.
\medskip

\emph{Inequality \eqref{RefPathCondNN} holds.}
The proof relies on the detailed computations made in (3.4a)--(3.4e) in \cite{KO94}.
Let $ k^\star  $ be such that $  \gamma_\mathrm{NN}(k^\star) \in Q(D^\star - 1,D^\star) $.
Then $ \mathrm{H}_\mathrm{NN}( \gamma_\mathrm{NN}(k^\star)) = \mathrm{E}_{\mathrm{NN}}^\star - 2 \tilde{J}  + h <  \mathrm{E}_{\mathrm{NN}}^\star$. 
If we go backwards in the path from that point on, then we will have to flip all $ (+1) $--spins on the NE-edge of $ Q(D^\star - 1,D^\star) $.
This is an increase of the energy in each step until only one $ (+1) $--spin remains on this edge (cf.\ \cite[(3.4a)]{KO94}).
At this point the energy equals to 
\begin{align} 
\mathrm{E}_{\mathrm{NN}}^\star - 2  \tilde{J}  + \ell^\star  h  < \mathrm{E}_{\mathrm{NN}}^\star, 
\end{align}
where we have used Assumption \ref{AssNN} b) and \eqref{InequNN}.
Flipping the last $ (+1) $--spin on this edge decreases the energy by $ 2K-h $ (cf.\ \cite[(3.4c)]{KO94}, but with here we flip a $ (+1) $--spin). 
Next, we do the same thing on the SE-edge, i.e.\ we flip all but one $ (+1) $--spins on this edge and arrive at the energy
\begin{align}
\mathrm{E}_{\mathrm{NN}}^\star - 2  \tilde{J} -2K + 2\ell^\star  h  < \mathrm{E}_{\mathrm{NN}}^\star.
\end{align}
Flipping the last $ (+1) $--spin on this edge, we arrive at $ \mathrm{E}_{\mathrm{NN}}^\star - 2  J +( 2\ell^\star +1) h $.
Finally, we need to flip all  but one $ (+1) $--spins on the E-edge, which leads to the energy level (cf.\ \cite[(3.4a)]{KO94})
\begin{align}
\mathrm{E}_{\mathrm{NN}}^\star - 2  J +( D^\star -1) h < \mathrm{E}_{\mathrm{NN}}^\star,
\end{align}
and flipping the last $ (+1) $--spin on this edge, we arrive at the energy $ \mathrm{E}_{\mathrm{NN}}^\star - 4  J +4K +D^\star  h $ (analogously to \cite[(3.4e)]{KO94}).
With the same reasoning, if we keep on going backwards in the path of $  \gamma_\mathrm{NN} $, we will always stay below $ \mathrm{E}_{\mathrm{NN}}^\star $, since the length of the edges of the circumscribing rectangles will be at most $ D^\star  -  1$.
Hence, we get that  
\be{RefPathFirstPartNN} 
\max_{i=1,\dots,k^\star}\mathrm{H}_\mathrm{NN}( \gamma_\mathrm{NN}(i)) < \mathrm{E}_{\mathrm{NN}}^\star .
\ee

We now analyze  the  path of $  \gamma_\mathrm{NN} $ after the step $ k^\star +2 $.
It holds that $ \mathrm{H}_\mathrm{NN}( \gamma_\mathrm{NN}(k^\star+2)) = \mathrm{E}_{\mathrm{NN}}^\star- h <  \mathrm{E}_{\mathrm{NN}}^\star$. 
First, $ L^\star -4$ $ (+1) $--spins are attached at the interior of the  S-edge.
The energy is decreased to $ \mathrm{E}_{\mathrm{NN}}^\star- (L^\star -3)h $.
Afterwards, a $ (+1) $--spin is added at the SW-edge, which leads to the energy (cf.\ \cite[(3.4c)]{KO94})
\begin{align} 
\mathrm{E}_{\mathrm{NN}}^\star+2K - (L^\star -2)h  <\mathrm{E}_{\mathrm{NN}}^\star,
\end{align}
where we have used the inequality $ L^\star \geq 2\ell^\star +1$, which follows immediately from Assumption \ref{AssNN} b).
Filling the SW-edge decreases the energy by $ (\ell^\star - 1)h $.
Then we do the same things for the SE-edge by attaching first a $ (+1) $--spin on this edge, which increases the energy to 
\begin{align}
\mathrm{E}_{\mathrm{NN}}^\star+4K - (L^\star +\ell^\star -2)h  <\mathrm{E}_{\mathrm{NN}}^\star,
\end{align}
and then filling up this edge, which decreases the energy to $ \mathrm{E}_{\mathrm{NN}}^\star+4K - (D^\star -1)h $.
Next, a protuberance is added at the interior of the E-edge.
We arrive at the energy level  (cf.\ \cite[(3.4e)]{KO94})
\begin{align} 
 \mathrm{E}_{\mathrm{NN}}^\star+4K +2 \tilde{J} - D^\star h = \mathrm{E}_{\mathrm{NN}}^\star+2J - D^\star h  <\mathrm{E}_{\mathrm{NN}}^\star.
\end{align}
If we keep  following the path of $  \gamma_\mathrm{NN} $, we will always stay below $ \G_{\mathrm{NN}}^\star $, since the length of the edges of the circumscribing rectangles will be at least $ D^\star $.
Combining this with \eqref{RefPathFirstPartNN} and the fact that $ \mathrm{H}_\mathrm{NN}( \gamma_\mathrm{NN}(k^\star+1)) = \mathrm{E}^\star_\mathrm{NN} $, we infer \eqref{RefPathCondNN}.

\subsection{Proof of $ \Phi(\boxminus,\boxplus) - \mathrm{H}_\mathrm{NN}(\boxminus) \geq  \G^\star _{\mathrm{NN}}$}

\label{Step2NN}
We first list a few observations taken from \cite{KO94}.
\bl{LocMinNN}
Let $ \s \in S $ be a local minimum of $ \mathrm{H}_\mathrm{NN} $. 
Then all clusters of  $ \s $ have distance at least $ \sqrt{2} $ from each other and each cluster is either a stable octagon or a rectangle that wraps around the torus. 
\el
\bpr
In   Lemma 2.1 in \cite{KO94} the following fact was proven.
Let $ \s \in S $ be a local minimum and let $ \s_1 $ be a cluster of $ \s $, then $ \s_1 = Q(\s_1)$, where  $ Q(\s_1) $ is the \emph{octagonal envelope} of $ \s_1 $, i.e.\
\begin{itemize}
	\item if $ \s_1 $ does not wind around the torus, then $ Q(\s_1) $ is the smallest octagon containing $ \s_1 $, and
	\item if $ \s_1 $  winds around the torus, then $ Q(\s_1) = R(\s_1) $.
\end{itemize}
See page 424 (before and after Scheme 3.2) in \cite{KO94} for the definition of the octagonal envelope. 
Moreover, it is shown, at the end of the proof of Lemma 2.1 in \cite{KO94}, that all clusters of  $ \s $ have distance at least $ \sqrt{2} $ and that if a cluster of $ \s $ is a octagon, it  must be a stable octagon.
\epr

\bl{MinRectNN}
Assume that $ \s \in S $ consists  of a unique cluster that does not wrap around the torus. 
Let $ R(\s) \in R( D_n,D_w) $ with $ D_n \geq D_w $.
\begin{itemize}
		\setlength\itemsep{-0.3em} 
	\item  If $ D_w \geq 2\ell^\star -1 $, then
\begin{align}
\mathrm{H}_\mathrm{NN}(\s) \geq  \mathrm{H}_\mathrm{NN}(Q(D_n,D_w)),
\end{align}
and equality holds if and only if $ \s=Q(D_n,D_w)	 $.

\item
If $ D_w < 2\ell^\star -1 \text{ and } D_w \text{ is odd}$, then
\begin{align}
\mathrm{H}_\mathrm{NN}(\s) \geq  \mathrm{H}_\mathrm{NN}(Q(D_n,D_w;\tfrac{1}{2}(D_w+1))),
\end{align}
and equality holds if and only if $ \s=Q(D_n,D_w;\tfrac{1}{2}(D_w+1))	 $.

\item
If $ D_w < 2\ell^\star -1 \text{ and } D_w \text{ is even}$, then
\begin{align}
\mathrm{H}_\mathrm{NN}(\s) \geq  \mathrm{H}_\mathrm{NN}(Q(D_n,D_w;\tfrac{1}{2}D_w,\tfrac{1}{2}D_w,\tfrac{1}{2}D_w+1,\tfrac{1}{2}D_w+1)),
\end{align}
and equality holds if and only if $ \s=Q(D_n,D_w;\tfrac{1}{2}D_w,\tfrac{1}{2}D_w,\tfrac{1}{2}D_w+1,\tfrac{1}{2}D_w+1)	 $.
\end{itemize}
%
\el
\bpr
See Lemma 3.2 and the proof of Lemma 4.1A in \cite{KO94}. 
The main step is to show that the function $ l \mapsto F(l)  $ is minimized in $ \ell^\star $.
\epr

In the following lemma we show that every optimal path has to cross $ Q(D^\star - 1,D^\star) $.
\bl{Gate1NN}
Let $ \gamma \in (\boxminus , \boxplus )_{\mathrm{opt}} $. 
Then $ \gamma $  has to cross $ Q(D^\star - 1,D^\star) $.
\el
\bpr
Assume the contrary, i.e.\ $ \gamma \cap Q(D^\star - 1,D^\star) = \emptyset$.
Using the same arguments as in the end of the proof of Lemma \ref{Gate1A}, we can restrict   to the case that throughout its whole path $ \gamma $ consists only of a unique cluster. 
On its way to $ \boxplus $, $ \gamma  $ has to cross a configuration, whose rectangular envelope has both  horizontal and vertical length greater or equal to $ D^\star $.
Let 
\begin{align} 
\bar{t}= \min\{ l\geq 0 \, | \, P_HR(\gamma(l)) , P_VR(\gamma(l)) \geq  	D^\star  	\} .
\end{align}
By the definition of $ \bar{t} $, we have that $ R(\gamma(\bar{t}-1))\in R(D^\star +m, D^\star -1 ) $ for some $ m\geq  0 $.

  \emph{Case 1} . [$ m=0$].\\
Obviously, $ D^\star \geq 2\ell^\star -1 $.
Hence, by 
Lemma \ref{MinRectNN} and since $ \gamma  $ does not cross $ Q(D^\star - 1,D^\star) $, we have that
\begin{align}
\mathrm{H}_\mathrm{NN}(\gamma(\bar{t}-1)) > \mathrm{H}_\mathrm{NN} (Q(D^\star  -1 , D^\star)) =  \mathrm{E}_{\mathrm{NN}}^\star -2\tilde{J} +h.
\end{align}
The minimal increase of energy to enlarge the rectangular envelope of a configuration is $ 2\tilde{J} -h  $.
Hence,
\begin{align}
\mathrm{H}_\mathrm{NN}(\gamma(\bar{t})) \geq \mathrm{H}_\mathrm{NN}(\gamma(\bar{t}-1)) + 2\tilde{J} -h>   \mathrm{E}_{\mathrm{NN}}^\star.
\end{align}
This contradicts the fact that $ \gamma \in (\boxminus , \boxplus )_{\mathrm{opt}} $, since 
$ \Phi(\boxminus,\boxplus) \leq \max_{\eta \in  \gamma_\mathrm{NN} } \mathrm{H}_\mathrm{NN}(\eta) \leq   \mathrm{E}_{\mathrm{NN}}^\star  $, where 
$ \gamma_\mathrm{NN} $ was constructed in  Section \ref{Step1NN}. 

 \emph{Case 2} . [$ m\in [1,\sqrt{|\L|} - D^\star)$].\\
Again, by Lemma \ref{MinRectNN}  we have that 
\begin{align}
\begin{split}
\mathrm{H}_\mathrm{NN}(\gamma(\bar{t}-1)) &\geq \mathrm{H}_\mathrm{NN} (Q( D^\star  +m, D^\star -1)) \\
&=\mathrm{H}_\mathrm{NN} (Q( D^\star   , D^\star-1))+ m(2J- h (D^\star -1))\\
&> \mathrm{E}_{\mathrm{NN}}^\star -2\tilde{J}+h.
\end{split}
\end{align}
As before, this leads to a contradiction, since
\begin{align}
\mathrm{H}_\mathrm{NN}(\gamma(\bar{t})) \geq \mathrm{H}_\mathrm{NN}(\gamma(\bar{t}-1)) + 2\tilde{J} -h>   \mathrm{E}_{\mathrm{NN}}^\star.
\end{align}

 \emph{Case 3} . [$ m=\sqrt{|\L|} - D^\star$].\\
In this case, $ \gamma(\bar{t}-1) $ wraps around the torus. 
One can easily observe that $ \mathrm{H}_\mathrm{NN}(\gamma(\bar{t}-1)) \geq \mathrm{H}_\mathrm{NN} ( R( \sqrt{|\L|} , D^\star  -1)) $.
We infer that 
\begin{align}
\mathrm{H}_\mathrm{NN}&(\gamma(\bar{t}-1)) \geq \mathrm{H}_\mathrm{NN} ( R( \sqrt{|\L|} , D^\star  -1)) \nonumber \\
&=\mathrm{H}_\mathrm{NN} (Q( D^\star, D^\star -1) ) +(\sqrt{|\L|} - D^\star ) (2J - h(D^\star -1)) - 2J(D^\star -1) - 4F(\ell^\star)\nonumber\\
&> \mathrm{H}_\mathrm{NN} (Q( D^\star , D^\star -1) )
= E_\mathrm{NN}^\star - 2\tilde{J} +h, 
\end{align}
where we have used  that $ F(\ell^\star) <0 $ and Assumption \ref{AssNN} d).
Finally, 
\begin{align}
\mathrm{H}_\mathrm{NN}(\gamma(\bar{t})) \geq \mathrm{H}_\mathrm{NN}(\gamma(\bar{t}-1))+ 2\tilde{J} -h>   E_\mathrm{NN}^\star.
\end{align}
This concludes the proof. 
\epr

Finally, the following lemma  concludes the proof of  $ \Phi(\boxminus,\boxplus) - \mathrm{H}_\mathrm{NN}(\boxminus) \geq  \G^\star _{\mathrm{NN}}$.

\bl{Gate2NN}
Let $ \gamma\in (\boxminus , \boxplus )_{\mathrm{opt}}  $.
In order  to cross a configuration whose rectangular envelope has both vertical and horizontal length greater or equal to $ D^\star $, $ \gamma $ has to pass through $ Q(D^\star - 1,D^\star) $ and $ Q(D^\star - 1,D^\star) ^\mathrm{1pr} $.
In particular, each optimal path between $  \boxminus $ and $ \boxplus $ has to cross $ Q(D^\star - 1,D^\star) ^\mathrm{1pr} $.
\el
\bpr
Consider the time step $ \bar{t} $ defined in the proof of Lemma \ref{Gate1NN}. 
It was shown there that necessarily $ \gamma(\bar{t}-1) $ needs to belong to $ Q(D^\star - 1,D^\star)$. 
Since $ P_VR(\gamma(\bar{t})),P_HR(\gamma(\bar{t})) \geq D^\star  $, 
$ \gamma(\bar{t}) $ must be obtained from $ \gamma(\bar{t}-1) $
by flipping a  $ (-1) $--spin at a site that is attached at the coordinate edge of a longer side of the droplet.
If it would not attach at the interior of the coordinate edge, then the energy level $ \mathrm{E}_{\mathrm{NN}}^\star +2K$ would be reached.
Hence, the protuberance must be added at the interior of the coordinate edge, which implies that $ \gamma(\bar{t}) $ needs to belong to $ Q(D^\star - 1,D^\star)^\mathrm{1pr} $.
\epr

\subsection{Identification of $ \cP^ \star $ and	$ \cC^ \star $}

\label{Step3NN}

In Section \ref{Step1NN} we have seen  that $ Q(D^\star - 1,D^\star) \subset  \cP^ \star $.
Now let $ \s \in \cP^ \star  $ and $ x \in \L $ be such that $ \s^x \in \cC^ \star $.
If follows from the definition of $ \cP^ \star $ and $ \cC^ \star $ that there exists $ \bar{\gamma} \in (\boxminus,\boxplus)_{\mathrm{opt}}$ and $ \ell \in \N $ such that 
\begin{enumerate}[(i)]
	\item $ \bar{\gamma}(\ell)= \s $ and $ \bar{\gamma}(\ell +1)=\s^x$,
	\item $ \mathrm{H}_\mathrm{NN}(\bar{\gamma}(k)) < \mathrm{E}_{\mathrm{NN}}^\star  $ for all $  k \in \{ 0,\dots,\ell	 	\}  $,
	\item $ \Phi(\boxminus,\bar{\gamma}(k) ) \geq  \Phi(\bar{\gamma}(k) , \boxplus )$  for all $  k \geq \ell + 1  $.
\end{enumerate}
By Lemma \ref{Gate2NN}, (ii) implies that  $\min(P_HR(\s),P_VR(\s))  \leq  D^\star  -1 $, since otherwise the energy level $ \mathrm{E}_{\mathrm{NN}}^\star $ would have been reached.
There are two possible cases.

  \emph{Case 1}. [$ P_HR(\s^x), P_VR(\s^x)  \geq  D^\star $].\\
Lemma \ref{Gate2NN} implies that necessarily $ \s \in Q(D^\star - 1,D^\star) $ and  $ \s^x \in Q(D^\star - 1,D^\star) ^\mathrm{1pr} $.

   \emph{Case 2}. [$ \min( P_HR(\s^x),P_VR(\s^x))  \leq  D^\star -1 $].\\
Also by  Lemma \ref{Gate2NN}, there must exist some $ k^\star \geq \ell + 2  $ such that $ \bar{\gamma}(k^\star ) \in Q(D^\star - 1,D^\star) $. 
But this contradicts (iii), since $ \Phi(\boxminus,\gamma(k^\star ) ) < \Phi(\bar{\gamma}(k^\star) , \boxplus )= \mathrm{E}_{\mathrm{NN}}^\star $.
%
Hence, only Case 1 can hold true.

We conclude that $ \cP^ \star=Q(D^\star - 1,D^\star) $ and	$ \cC^ \star= Q(D^\star - 1,D^\star) ^\mathrm{1pr} $.

\subsection{Verification of (H1)}

\label{Step4NN}
Obviously, $ S_{\mathrm{stab}}= \{ \boxplus	\} $, since $ \boxplus $ minimizes all three sums in \eqref{HNN}.
It remains to show that $ S_{\mathrm{meta}}= \{ \boxminus	\}$.

Let $ \s \in S\setminus\{\boxminus,\boxplus  \} $. 
As in Section \ref{Step4A}, we have to show that there exists $ \s' \in S $ such that $ \mathrm{H}_\mathrm{NN}(\s' ) < \mathrm{H}_\mathrm{NN}(\s)$ and $ \Phi ( \s,\s') - \mathrm{H}_\mathrm{NN}(\s) < \G_{\mathrm{NN}}^\star $. 
\medskip
 
  \emph{Case 1}. [$ \s $ contains a cluster, which is not a stable octagon and not a rectangle that wraps around the torus].\\
Lemma \ref{LocMinNN} implies that $ \s $ is not a local minimum, i.e.\ there exists $ x \in \L $ such that $ \mathrm{H}_\mathrm{NN}(\s^x ) < \mathrm{H}_\mathrm{NN}(\s)$ and  $ \Phi ( \s,\s^x) - \mathrm{H}_\mathrm{NN}(\s) =0 < \G_{\mathrm{NN}}^\star $.

 \emph{Case 2}. [$ \s $ contains a  cluster $ Q $, which is   a stable octagon  with $ D^\star\leq P_VR(Q) \leq \sqrt{\L}-1 $].\\
Let $ \s' $ be obtained from $ \s $ by attaching at $ Q $ an oblique bar at its NE-edge and its SE-edge respectively, and a vertical bar at its E-edge  in the same way that  was described in the third step of the construction of $ \gamma_\mathrm{NN} $ 
given in Section \ref{Step1NN}.
Then we obtain
\begin{align}
\begin{split}
\mathrm{H}_\mathrm{NN}(\s' ) -  \mathrm{H}_\mathrm{NN}(\s ) &\leq 2J - P_VR(Q) h \leq  2J - D^\star  h< 0, \quad \text{and} \\
\Phi ( \s,\s') - \mathrm{H}_\mathrm{NN}(\s) &\leq 2\tilde{J} - h < \G_{\mathrm{NN}}^\star .
\end{split}
\end{align}

  \emph{Case 3}. [$ \s $ contains a cluster $ Q $, which is a stable octagon with $ P_VR(Q) \leq D^\star -1 $].\\
Let $ \s' $ be obtained from $ \s $ as follows. 
First the uppermost $ (+1)$--spin at the NE-edge is flipped. 
Afterwards, successively,  adjacent $ (+1)$--spins are flipped until this oblique bar consist only of $ (- 1)$--spins. 
In the same way, the SE-edge and the E-edge of $ Q $  are detached by starting from the uppermost $ (+1)$--spin and then successively flipping all adjacent $ (+1)$--spins until the respective edge is detached from $ Q $.
Then
\begin{align}
\begin{split}
\mathrm{H}_\mathrm{NN}(\s' ) -  \mathrm{H}_\mathrm{NN}(\s ) &= - 2J + P_VR(Q) h 
\leq - 2J +  (D^\star-1) h 
< 0, \quad \text{and} \\
\Phi ( \s,\s') - \mathrm{H}_\mathrm{NN}(\s) &\leq (P_VR(Q) -1)h< \G_{\mathrm{NN}}^\star .
\end{split}
\end{align}

 \emph{Case 4}. [$ \s $ contains a cluster $ R $ that  is a rectangle that wraps around the torus.].\\
Let $ \s' $ be obtained from $ \s $ by attaching at $ R $ a bar that also wraps around the torus.
Then, by Assumption \ref{AssA} d), we have that 
\begin{align}
\begin{split}
\mathrm{H}_\mathrm{NN}(\s' ) -  \mathrm{H}_\mathrm{NN}(\s )&= 2\tilde{J} - \sqrt{|\L|} h < 0, \quad \text{and} \\
\Phi ( \s,\s') - \mathrm{H}_\mathrm{NN}(\s) &= 2\tilde{J} - h < \G_{\mathrm{NN}}^\star .
\end{split}
\end{align}

We conclude that $ S_{\mathrm{meta}}= \{ \boxminus	\}. $

\subsection{Verification of (H2)}

\label{Step5NN}

Obviously, $ |\{ \s \in \cP^\star \,  | \,  \s \sim \s' \}| = 1 $ for all $ \s' \in \cC^\star$.
Therefore, (H2) holds.

\subsection{Computation of $K$}

\label{Step6NN}

We proceed analogously to Section \ref{Step6A}. 

  \emph{Lower bound}. 
Note that $ \partial^+ \cC^\star \cap S^\star = Q(D^\star - 1,D^\star)  \cup Q(D^\star - 1,D^\star) ^\mathrm{2pr}  $, $  Q(D^\star - 1,D^\star)  \subset  S_\boxminus  $ and $ Q(D^\star - 1,D^\star) ^\mathrm{2pr} \subset S_\boxplus  $.
Hence,
\begin{align}
\begin{split}
\frac{1}{K} &\geq \min_{C_1,\dots,C_I \in [0,1]}
\min_{\substack{ h:S^\star \ra [0,1] \\
		\smallrestr{h}{S_\boxminus} =1, \smallrestr{h}{S_\boxplus} =0, \smallrestr{h}{S_i} =C_i \, \forall i }}
\frac{1}{2} \sum_{\eta,\eta' \in (\cC^\star)^+ } \mathbbm{1}_{\{	\eta \sim \eta'	 \}} [h(\eta) - h(\eta')]^2\\
&=
\min_{ h:\cC^\star \ra [0,1]}
\sum_{\eta \in \cC^\star } \left( \sum_{\eta' \in Q(D^\star - 1,D^\star) , \eta'\sim \eta  }	[1- h(\eta) ]^2	+ \sum_{\eta' \in Q(D^\star - 1,D^\star) ^\mathrm{2pr}, \eta'\sim \eta  }	h(\eta)^2			\right) \\
&=
\sum_{\eta \in \cC^\star }  \min_{ h \in [0,1]} \Big( |Q(D^\star - 1,D^\star) \sim \eta| \ 	[1- h ]^2	+ 
|Q(D^\star - 1,D^\star) ^\mathrm{2pr}\sim \eta| \  h^2			\Big) \\
&=
\sum_{\eta \in \cC^\star } \frac
{|Q(D^\star - 1,D^\star) \sim \eta| 	\cdot 
	|Q(D^\star - 1,D^\star) ^\mathrm{2pr}\sim \eta| }
{|Q(D^\star - 1,D^\star) \sim \eta| 	+ 
	|Q(D^\star - 1,D^\star) ^\mathrm{2pr}\sim \eta| }.
\end{split}
\end{align}
For all $ \eta \in \cC^\star  $ we have that $ |Q(D^\star - 1,D^\star) \sim \eta| = 1 $.
Moreover, there are four sites at the  longer coordinate edges of a critical droplet with $ |Q(D^\star - 1,D^\star) ^\mathrm{2pr}\sim \eta|=1 $, and $ 2(L^\star -4)  $ sites with  $ |Q(D^\star - 1,D^\star) ^\mathrm{2pr}\sim \eta|=2 $.
Further, there are $ |\L| $ possible locations for a configuration in $\cC^\star $, and  there are two analogue rotations for each critical droplet. 
Therefore, we obtain that
\begin{align}
\frac{1}{K} &\geq
\Big( 2(L^\star-4) \frac{2} {3 } + 4 \frac {1 } {2 } \Big)2 |\L|= \frac{4(2L^\star-5)}{3}|\L|.
\end{align}

  \emph{Upper bound}.
  The following proof uses the same arguments as in Section \ref{Step6A}. Hence, we shall only sketch the main arguments here. 
Define
\begin{align}
\begin{split}
\mathscr{S}^- =  \{\s\in S^\star\, &| \, 
\min(P_HR(\eta),P_VR(\eta))  \leq  D^\star  -1
 \text{ for all clusters  $ \eta $ of $ \s $} \} , \text{ and }\\
\mathscr{S}^+ =  \{\s\in S^\star\, &| \,  \text{ there exists a cluster  $ \eta $ of $ \s $ such that  }
P_HR(\eta),P_VR(\eta)  \geq  D^\star \  \}.
\end{split}
\end{align} 
Note that $ S^\star = \mathscr{S}^-\cup \mathscr{S}^+$, $ \mathscr{S}^- \cap \mathscr{S}^+ = \emptyset $, $ \cP^\star \subset \mathscr{S}^- $ and $ \cC^\star \subset \mathscr{S}^+ $.
\bl{UpperBoundNN}
Let $ \s \in \mathscr{S}^- $ and $ \s' \in \mathscr{S}^+ $.
Then $ \s \sim \s' $ if and only if  $ \s \in \cP^\star$ and $ \s' \in \cC^\star$.
\el
\bpr
The proof is a straightforward adaptation of the proof of Lemma \ref{Gate1NN}.
\epr

\noindent
The same arguments as in Section \ref{Step6A} yield that $ S_\boxminus \subset \mathscr{S}^- $, $ S_\boxplus \subset \mathscr{S}^+ $ and for all $ i =0,\dots,I $ we have that  either $ S_i \subset \mathscr{S}^- $  or $ S_i \subset \mathscr{S}^+ $.
Therefore, as in Section \ref{Step6A}, we can estimate the minimum in \eqref{EqPrefactor} from above by the minimum over all functions of the form \eqref{EqStep6FunctionA} and use Lemma \ref{UpperBoundNN} to infer that 
\begin{align}
\begin{split}
\frac{1}{K} &\leq 
\min_{\substack{ h:S^\star \ra [0,1] \\
		\smallrestr{h}{\mathscr{S}^-} =1, \smallrestr{h}{\mathscr{S}^+\setminus \cC^\star } =0}}
\frac{1}{2} \sum_{\eta,\eta' \in S^\star } \mathbbm{1}_{\{	\eta \sim \eta'	 \}} [h(\eta) - h(\eta')]^2 \\
&=
\min_{\substack{ h:(\cC^\star)^+ \ra [0,1] \\
		\smallrestr{h}{\mathscr{S}^-\cap \partial^+\cC^\star} =1, \smallrestr{h}{\mathscr{S}^+\cap \partial^+\cC^\star} =0}}
\frac{1}{2} \sum_{\eta,\eta' \in (\cC^\star)^+ } \mathbbm{1}_{\{	\eta \sim \eta'	 \}} [h(\eta) - h(\eta')]^2 \\
&=
\min_{ h:\cC^\star \ra [0,1]}
\sum_{\eta \in \cC^\star } \left( \sum_{\eta' \in Q(D^\star - 1,D^\star) , \eta'\sim \eta  }	[1- h(\eta) ]^2	+ \sum_{\eta' \in Q(D^\star - 1,D^\star) ^\mathrm{2pr}, \eta'\sim \eta  }	h(\eta)^2			\right) \\
&=
\frac{4(2L^\star-5)}{3}|\L|.
\end{split}
\end{align}

\section{Ising model with alternating magnetic field}
\label{pm}

We adapt the same strategy as in the Chapters 2 and 3 to a third modification of the Ising model (cf.\ Section \ref{resultsecpm}), where 
 the Hamiltonian is given by
\be{Hpm}
\mathrm{H}_{\pm}(\s) = - \frac{J}{2} \sum_{(x,y) \in \L^\star} \s(x)\s(y)
+ \frac{\ho }{2} \sum_{x \in \L_\mathrm{2}} \s(x)
- \frac{h_\mathrm{1}}{2} \sum_{x \in \L_\mathrm{1}} \s(x),
\ee
where $\s \in S$, $ J ,\ho,\he >0 $,  $ \Lo = \{	 (x_1,x_2)\in \L \, | \, x_2 \text{ is odd}	\} $ are the \emph{odd rows} in $ \L $,  $ \Le = \L\setminus\Lo $ are the \emph{even rows} and   $ \L^\star $ is the set of \textit{unordered nearest-neighbor bonds} in $ \L $.
%
One can rewrite $ \mathrm{H}_{\pm}(\s) $ geometrically as
\be{Hgeompm}
\mathrm{H}_{\pm}(\s) = 
\mathrm{H}_{\pm}(\boxminus) 
+\ho |\s \cap \Lo|	- \he |\s \cap \Le|
+J |\partial(\s)|.
\ee
Under the assumptions below, the critical lengths in this model are given by
\begin{align}\label{Critlengthspm}
l_b^\star = \left\lceil  \frac{\mu}{\e} \right\rceil  
\qquad \text{and} \qquad l_h^\star = 2 l_b^\star - 1,
\end{align}
where
\begin{align}
\begin{split}
& \e = \he - \ho , \quad \text{and} \\
&\mu = 2J - \ho.
\end{split}
\end{align}
$ l_b^\star $ will be the length of the basis of the critical droplet, and $ l_h^\star $ will be its height. 
The following assumptions will be made for this chapter.
\bass{Asspm}
\begin{enumerate}[a)]
	\item $ \he > \ho $,
	\item $ J > \he $,
	\item $ \frac {\mu	}{\e} \notin \N $, 
	\item $ |\L| > \left(  2 \left\lceil   \frac{2J (l_h^\star-1) +\ho}{4J-\e (l_b^\star-1)} \right\rceil   + l_h^\star \right)^2 $. 
\end{enumerate}
\eass
Assumption a) ensures that $ \boxplus $ is the \emph{stable configuration} in this system.
Assumptions b), c) and d) are made due to similar reasons as in the Chapters 2 and 3. 
Assumption b) can also be modified in various ways. 
E.g.\ one can take $ J<\he< 2J $. 
We refer to \cite{NO96}, page 10, where several other regimes are listed. 
In contrast to \cite{NO96}, in  this text, we only consider the regime given in Assumption \ref{Asspm}, since all other regimes can be handled in a similar way without using new ideas. 
It immediately follows from Assumption \ref{Asspm} c) that
\begin{align}\label{Inequpm}
(l_b^\star -1)\e < \mu  < l_b^\star  \e.
\end{align}
%

\begin{figure}[htbp]
	\centering
	\includegraphics[height= 5cm]{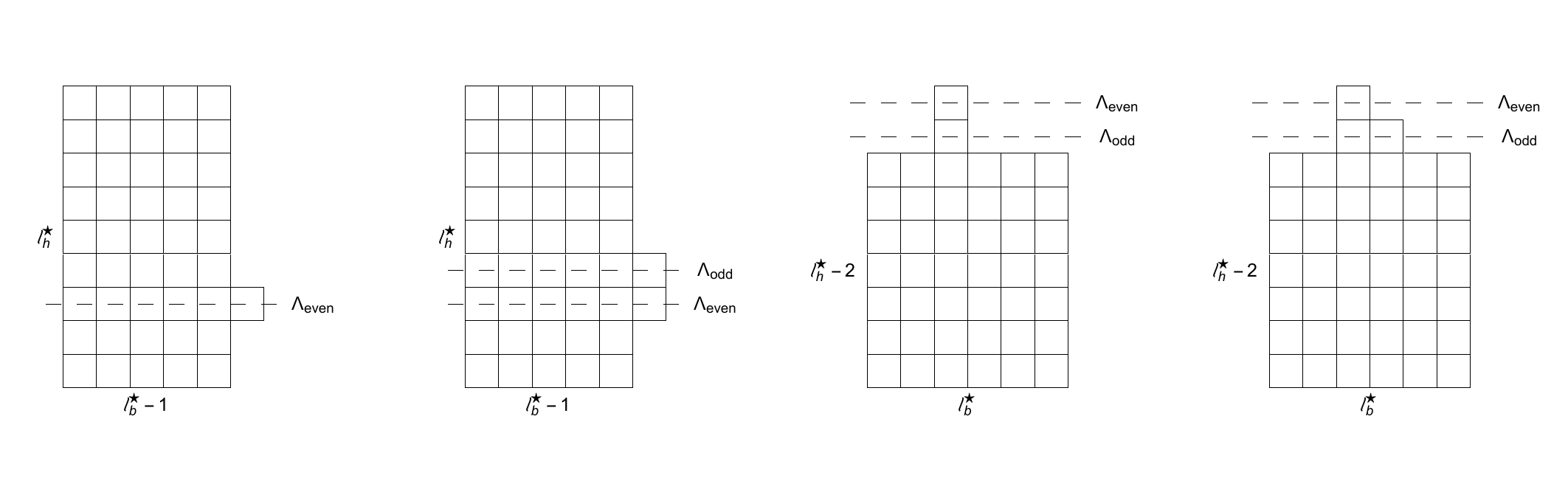}	
	\vspace{-1.5cm} \caption{{\footnotesize From left to right: An example of an element in $ \cP_1$, $ \cC_1$, $ \cP_2$ and $ \cC_2$}}
	\label{PandC_pm}
\end{figure}

In the following definition we define the protocritical and the critical configurations for this model.
Figure \ref{PandC_pm} below provides an example.

\bd{defipm} 
Let $ \s \in S $ consist of a unique cluster. 
$ l\in R(1\times 2) $ is called a \emph{2-protuberance attached at $ \s $} if there exists $ x \in l $ and $ \bar{y} \in \s $ such that $ |x-\bar{y}|=1 $ and $ \sum_{y \in \L: |y-x|=1} \s(y) = 0 $ and $ \sum_{y \in \L: |y-x'|=1} \s(y) = -2 $, where $ x' $ is the unique element in $ l \setminus x $.

We define the following subsets of $ S$.

$\cP_1 $ denotes the set of all configurations consisting only of a rectangle from $ R (( l_b^\star  -1) \times l_h^\star) $ that starts and ends in $ \Le $ (i.e.\ the bottom and the top row belong to $ \Le $)  and with an additional protuberance  attached at one of its vertical sides on a row in $ \Le $.

$\cC_1 $ denotes the set of all configurations that are obtained from a configuration in $\cP_1$ by adding a second $ (+1) $--spin  in $ \Lo $ adjacent to the protuberance and attached at the rectangle.

$\cP_2' $ denotes the set of all configurations consisting only of a rectangle from $ R (l_b^\star  \times (l_h^\star-2)) $ that starts and ends  in $ \Le $ and with an additional  horizontal bar of length $ 2 $ attached at one of the horizontal sides of the droplet.

$\cP_2'' $ denotes the set of all configurations consisting only of a rectangle from $ R (l_b^\star  \times (l_h^\star-2)) $ that starts and ends  in $ \Le $ and with an additional 2-protuberance attached at one of the horizontal sides of the droplet.

Define $ \cP_2 = \cP_2'\cup \cP_2'' $.

$\cC_2' $ denotes the set of all configurations that are obtained from a configuration in $\cP_2' $ by adding a  $ (+1) $--spin in $ \Le $ attached to the horizontal bar of length $ 2 $.

$\cC_2'' $ denotes the set of all configurations that are obtained from a configuration in $\cP_2'' $ by adding a  $ (+1) $--spin, which is both attached  to the 2-protuberance and to the rectangle. 

We easily observe that $  \cC_2'= \cC_2'' $. Define $ \cC_2 = \cC_2'= \cC_2'' $.
\ed

We now state the main result of this chapter.
\bt{Thmpm}
Under Assumption \ref{Asspm}, the pair $ (\boxminus,\boxplus) $ satisfies (H1) so that Theorem \ref{T1} a), Theorem \ref{T2} and Theorem \ref{T3} hold for the Ising model with alternating magnetic field.
Moreover,  
\bi
\item	$ \cP^ \star =\cP_1 \cup \cP_2$,
\item	$ \cC^ \star =\cC_1 \cup \cC_2 $,
\item $  \Phi(\boxminus,\boxplus) - \mathrm{H}_{\pm}(\boxminus) =  4J \, l_b^\star + \mu (l_b^\star - 1) - \e ( l_b^\star  (l_b^\star - 1) + 1 )  =:\G^\star_{\pm} =: \mathrm{E}_{\pm}^\star  - \mathrm{H}_{\pm}(\boxminus) $,
\item $K^{-1}=  
\frac{14\,( l_b^\star-1)}{3}|\L| $.
\ei
\et
\bpr
The proof is divided into the Sections \ref{Step1pm}--\ref{Step6pm}. 
\epr

\subsection{Proof of $ \Phi(\boxminus,\boxplus) - \mathrm{H}_{\pm}(\boxminus) \leq  \G^\star_{\pm} $}

\label{Step1pm}

As in Chapters 2 and 3, we construct a reference path $ \gamma_\pm: \boxminus \ra \boxplus$ such that 
\be{RefPathCondpm}
\max_{\eta \in \gamma_\pm } \mathrm{H}_{\pm}(\eta) \leq  \mathrm{H}_{\pm}(\boxminus) +  \G^\star_{\pm} = \mathrm{E}_{\pm}^\star .
\ee 

\emph{Construction of $ \gamma_\pm $.} $ \gamma_\pm $ is given through the following scheme. 

\textbullet \ Let $ \gamma_\pm(0)= \boxminus $.

\textbullet \ In the first step an arbitrary $(-1)$--spin  in $ \Le $ is flipped. 

 \textbullet \ [From $  R(l\times(2l-1)) $
   to $  R((l+1)\times(2l+1))   $ for  $ l \leq l_b^\star -1 $.]\\
 A protuberance is added to the right vertical side of the droplet at a row that belongs to $ \Le $.
Then successively  adjacent $(-1)$--spins are flipped  until the droplet belongs to $ R\left((l+1)\times(2l-1)\right) $.
Next,
a 
protuberance is added to the above horizontal side of the droplet, which is an odd row.
Afterwards, a second $(+1)$--spin is added above the protuberance on the even row. 
Hence, a
2-protuberance  attached  to the above horizontal side of the droplet was added.
Then, analogously as for this 2-protuberance, one adds successively adjacent $ 1\times 2  $ rectangles at the above horizontal side of the droplet until   $R((l+1)\times(2l+1))  $ is reached.

\textbullet \ [From $  R(l\times l_h^\star) $ 
to $  R((l+1)\times l_h^\star)   $ for  $ l \geq l_b^\star $.]\\
 A protuberance is added on the right  vertical side of the droplet at a row that belongs to $ \Le $, and successively  adjacent $(-1)$--spins are flipped  until the droplet belongs to $ R((l+1)\times l_h^\star) $.

\textbullet \ [From $  R(\sqrt{|\L|}\times l_h^\star)   $
 to $    \boxplus  $.]\\
As above, a  2-protuberance is added to the above horizontal side of the droplet, which is an odd row, 
and successively adjacent $ 1\times 2  $ rectangles  are added at the above horizontal side of the droplet, until  a configuration in $  R(\sqrt{|\L|}\times (l_h^\star+2)) $ is reached.  
This procedure is repeated until the configuration $ \boxplus $ appears.

\emph{Inequality \eqref{RefPathCondpm} holds.}
Let $ k^\star  $ be such that $ \gamma_\pm(k^\star) \in R ( l_b^\star \times (l_h^\star-2)) $.
Using \eqref{Hgeompm} and Assumption \ref{Asspm}, we observe that 
\begin{align}\label{EqInequPM}
\begin{split}
 \mathrm{H}_{\pm}(\gamma_\pm(k^\star)) &=
 \mathrm{H}_{\pm}(\boxminus) +  6J (  l_b^\star - 1 ) - \he  l_b^\star (l_b^\star - 1)  + \ho l_b^\star  (l_b^\star - 2)  \\
 &=
 \mathrm{H}_{\pm}(\boxminus) +  4J \,(  l_b^\star - 1 ) + \mu (l_b^\star - 1) - \e l_b^\star  (l_b^\star - 1)   - \ho  \\
 &=  \mathrm{E}_{\pm}^\star - 4 J +\e - \ho <  \mathrm{E}_{\pm}^\star.
\end{split}
\end{align} 
If we go backwards in the path from that point on, then we will have to cut the right vertical bar of the droplet.
While cutting this vertical bar, the highest energy level is reached when only two adjacent $(+1)$--spins remain, one in $ \Le $ and one in $ \Lo $. 
Indeed, at that point the energy in \eqref{EqInequPM} is increased by $ \e /2  (l_b^\star - 3)   + \ho $, so that it equals 
\begin{align}
\mathrm{E}_{\pm}^\star - 4J+  \e (l_b^\star -1)      < \mathrm{E}_{\pm}^\star,
\end{align}
where we have used \eqref{Inequpm}.
Cutting the last $(+1)$--spins, we reach $ R (( l_b^\star  -1) \times (l_h^\star-2)) $ and the energy decreases  to $ \mathrm{E}_{\pm}^\star - 6J+  \e l_b^\star  $. 
Next, we have to cut
the above two rows by successively cutting vertical bars of length $ 2 $ in these rows.
Doing that, the  highest energy point is the stage, where only  one vertical bar of length $ 2 $ and a single $(+1)$--spin in $ \Lo $ next to it have  remained. 
At this point the energy has increased by $ \e   (l_b^\star - 2)   + \he $
 and it 
equals to
\begin{align}
 \mathrm{E}_{\pm}^\star - 6J+  \e l_b^\star+ \e   (l_b^\star - 2)   + \he= \mathrm{E}_{\pm}^\star - 6J+ 2 \e (l_b^\star -1) + \he   < \mathrm{E}_{\pm}^\star.
\end{align}
Using the same arguments, if we keep on going backwards in the path of $ \gamma_{\pm} $, we will always stay below $ \mathrm{E}_{\pm}^\star $, since the sizes of the cut columns and rows further decrease. 
Hence, 
\be{RefPathFirstPartpm} 
\max_{i=1,\dots,k^\star}\mathrm{H}_{\pm}(\gamma_{\pm}(i)) < \mathrm{E}_{\pm}^\star .
\ee

We now consider the path of $ \gamma_{\pm} $ after the step $ k^\star +3 $.
We have that $\mathrm{H}_{\pm}(\gamma_\pm(k^\star +3) )=  \mathrm{E}_{\pm}^\star  $.
First,  the two rows above the droplet are filled.
This lowers the energy to 
$ \mathrm{E}_{\pm}^\star- \e (l_b^\star -1)- \ho $.
Afterwards, a protuberance is attached on the right vertical side of the droplet in a row that belongs to $ \Le $.
The energy is increased by $ 2J-\he $ and equals to
\begin{align} \mathrm{E}_{\pm}^\star +\mu  - \e l_b^\star -\ho < \mathrm{E}_{\pm}^\star.
\end{align}
Adding a second $ (+1) $--spin adjacent to the protuberance 
further increases the energy by $ \ho $. 
By  \eqref{Inequpm}, we still get 
\begin{align} \mathrm{E}_{\pm}^\star +\mu - \e l_b^\star < \mathrm{E}_{\pm}^\star .
\end{align}
If we fill this column, we further decrease the energy so that the energy still remains below $ \mathrm{E}_{\pm}^\star $.
In the following, analogously, columns are added successively on the right vertical side of the droplet and each column decreases the energy by $ \mu - \e l_b^\star  $.
This is repeated until the droplet wraps around the torus. 
It is easy to see that  the remaining part of $ \gamma $ also stays below $ \mathrm{E}_{\pm}^\star $.
Hence, 
\be{RefPathLastPartpm} 
\max_{i\geq k^\star+3}\mathrm{H}_{\pm}(\gamma_{\pm}(i)) \leq \mathrm{E}_{\pm}^\star .
\ee

Finally, we have that $ \mathrm{H}_{\pm}(\gamma_{\pm}(k^\star+1)) = \mathrm{E}_{\pm}^\star -2J + \he - \ho$ and $ \mathrm{H}_{\pm}(\gamma_{\pm}(k^\star+2)) = \mathrm{E}_{\pm}^\star - \ho $, which are clearly below $ \mathrm{E}_{\pm}^\star $. 
Hence, together with \eqref{RefPathFirstPartpm} and \eqref{RefPathLastPartpm}, we conclude \eqref{RefPathCondpm}.

\subsection{Proof of $ \Phi(\boxminus,\boxplus) - \mathrm{H}_{\pm}(\boxminus) \geq  \G^\star _{\pm}$}

\label{Step2pm}
Before we prove that $ \Phi(\boxminus,\boxplus) - \mathrm{H}_{\pm}(\boxminus) \geq  \G^\star _{\pm}$, we  need to collect some results that were established in \cite{NO96}.

\bd{DefiStablepm}
Let $ l_1 , l_2  \in \N $. 
We say that $  \s \in R(l_1 \times l_2) $ is a 
\emph{stable rectangle} if $ \s  $  starts and ends in $ \Le $ (i.e.\ its bottom and top row belong to $ \Le $), $ l_1 \geq 2 $ , $ l_2 \geq 3 $ and $ l_2  $ is odd.
Note that a stable rectangle  can possibly wrap around the torus. 
\ed

Recall \eqref{StabilityEq}.
Analogously to \cite{NO96}, we say that $ \s \in S $ is $ \ho $--stable if and only if $ \s \in S_{\ho} $.
 
\bl{LocMinpm}
$ \s \in S $ is $ \ho $--stable if and only if $ \s $ is a union of isolated stable rectangles.
\el
\bpr
This is the content of Proposition 3.1 in \cite{NO96} and the comment after it.
\epr

The following lemma is the analogue of  Corollary \ref{MinRectA}  for this model.  
\bl{MinRectpm}
Let $ \s \in S $ be such that $ R(\s) $ is a stable rectangle. 
Then
\begin{align}
\mathrm{H}_{\pm}(\s) \geq  \mathrm{H}_{\pm}(R(\s)),
\end{align}
and equality holds, if and only if $ \s=R(\s) $.
\el
\bpr
This is the content of Lemma 3.3 in \cite{NO96}.
\epr


%
Let $ l_1,l_2 \in \N $, and let $ R, R' \in R( l_1\times l_2) $.
Note that if $ l_2  $ is an odd number,   $ R $  starts  in $ \Lo $ and $ R' $  starts in $ \Le $, then 
$ \mathrm{H}_{\pm}(R)  > \mathrm{H}_{\pm}(R') $.
And if $ l_2  $ is even, then $ \mathrm{H}_{\pm}(R)  = \mathrm{H}_{\pm}(R') $.
Therefore, from now on, we set $ \mathrm{H}_{\pm}(l_1\times l_2)  = \mathrm{H}_{\pm}(R') $, which is the energetically more profitable choice. 
We will use this fact tacitly several times in the remaining part of this chapter.

\bl{Gate1pm}
Let $ \gamma \in (\boxminus , \boxplus )_{\mathrm{opt}} $. 
Then $ \gamma $  has to cross $ \cP_1 \cup \cP_2$.
\el

\bpr
Assume the contrary, i.e.\ $ \gamma \cap \{ \cP_1 \cup \cP_2 \}= \emptyset$.
Suppose first that throughout its whole path $ \gamma $ consists of a unique cluster.
At the end of this proof we treat the general case. 

Since $ \gamma $ leads to $ \boxplus $, there exists some time $ \bar{t} $ such that  $ P_VR(\gamma(j)) \geq l_h^\star $ and 
 $ P_HR(\gamma(j)) \geq l_b^\star $ for all $ j\geq \bar{t} $ and 
\begin{align}\label{Gate1pmEq00}
\bar{t}-1=  \max\big\{ j \geq 0 \ | \ P_VR(\gamma(j)) < l_h^\star \text{ or }  P_HR(\gamma(j)) < l_b^\star	\}.
\end{align} 
Note that $ \gamma(\bar{t}-1) $ has to   satisfy  either
\begin{enumerate}[1.)]
	\item $ P_HR(\gamma(\bar{t}-1)) = l_b^\star -1$ and $  P_VR(\gamma(\bar{t}-1)) =  l_h^\star +n $ for some $ n\geq 0 $, or
	\item $ P_VR(\gamma(\bar{t}-1)) = l_h^\star-1 $ and $  P_HR(\gamma(\bar{t}-1)) =  l_b^\star +m $ for some $ m\geq 0 $.
\end{enumerate}

  \emph{Case 1}. [$P_HR(\gamma(\bar{t}-1)) = l_b^\star -1 $ and $  P_VR(\gamma(\bar{t}-1)) =  l_h^\star +n $ for some $ n\geq 0 $].

 \emph{Case 1.1}. [$ n=0 $].\\
Let $ \tau $ be the first time that a second $ (+1) $--spin is added outside of $ R(\gamma(\bar{t}-1)) = R((l_b^\star -1 )\times l_h^\star )$, i.e.
\begin{align}
\tau= \min \big\{	j\geq \bar{t}+1 \ \big| \ |\gamma(j) \setminus R(\gamma(\bar{t}-1)) | = 	2	\big\}  	 
\end{align}
Note that $ |\gamma(\tau-1) \setminus R(\gamma(\bar{t}-1)) | = 	1 $ and that this protuberance is placed either at the right vertical side or at the left  vertical side of $ R(\gamma(\bar{t}-1)) $, since $ \gamma(\bar{t}-1) $ was the last configuration with the property $ P_HR(\gamma(\bar{t}-1)) = l_b^\star -1 $.
Analogously, $ P_VR(\gamma(\tau-1)) = l_h^\star $, otherwise, this would also contradict the definition of $ \bar{t}-1 $.
Now if $\gamma(\tau-1) \setminus R(\gamma(\bar{t}-1))  \in  \Lo $, we have that 
\begin{align}
\mathrm{H}_{\pm}(\gamma(\tau-1)) \geq  \mathrm{H}_{\pm}((l_b^\star-1)   \times  l_h^\star )+ 2J +\ho = \mathrm{E}_{\pm}^\star + \he> \mathrm{E}_{\pm}^\star.
\end{align}
This contradicts $ \gamma \in (\boxminus , \boxplus )_{\mathrm{opt}} $, since we already know from   Section \ref{Step1pm} that $ \Phi(\boxminus,\boxplus) \leq  E^\star _{\pm}$.
But if  $\gamma(\tau-1) \setminus R(\gamma(\bar{t}-1))  \in  \Le $, then,  since 
$ \gamma  $ does not cross $ \cP_1 $ and since the minimal increase of energy to enlarge the rectangular envelope is $ 2J -\he  $,  we have by Lemma \ref{MinRectpm} that 
\begin{align}
\mathrm{H}_{\pm}(\gamma(\tau-1)) >  \mathrm{H}_{\pm}((l_b^\star-1)   \times  l_h^\star ) + 2J -\he = \mathrm{E}_{\pm}^\star - \ho.
\end{align}
$ \gamma(\tau) $ is obtained from $ \gamma(\tau-1) $ by flipping a $ (-1) $--spin outside of $ R(\gamma(\bar{t}-1)) $.
One can easily see that the most profitable way is to flip a $ (-1) $--spin at a site that is adjacent to the protuberance of $ \gamma(\tau-1) $, which consequently must belong to $ \Lo $.
Hence,
\begin{align}
\mathrm{H}_{\pm}(\gamma(\tau)) \geq \mathrm{H}_{\pm}(\gamma(\tau-1)) + \ho > \mathrm{E}_{\pm}^\star ,
\end{align}
which leads to a contradiction.

\emph{Case 1.2}. [$ n=2k $ for some $ k>1 $].\\
According to Lemma \ref{MinRectpm}, we have that 
\begin{align}
\begin{split}
\mathrm{H}_{\pm}(\gamma(\bar{t})) &\geq \mathrm{H}_{\pm}(\gamma(\bar{t}-1)) + 2J-\he 
\geq \mathrm{H}_{\pm} (( l_b^\star  -1) \times (l_h^\star +2k)) + 2J - \he \\
&=\mathrm{H}_{\pm} (( l_b^\star  -1) \times l_h^\star) + k (4J- \e (l_b^\star -1) ) + 2J - \he
> \mathrm{E}_{\pm}^\star.
\end{split}
\end{align}
As before, this leads to a contradiction.

\emph{Case 1.3}. [$ n=2k+1 $ for some $ k\geq 0 $].\\
It holds  that either the top or bottom row of $ \gamma(\bar{t}) $ must belong to $ \Lo $. Similar to Case 1.2, we obtain a contradiction, since 
\begin{align}
\mathrm{H}_{\pm}(\gamma(\bar{t})) &\geq \mathrm{H}_{\pm}(\gamma(\bar{t}-1)) + 2J-\he
\geq \mathrm{H}_{\pm} (( l_b^\star  -1) \times (l_h^\star +2k)) + 4J + \ho - \he \nonumber \\
&\geq \mathrm{H}_{\pm} (( l_b^\star  -1) \times l_h^\star ) + 4J + \ho - \he>  \mathrm{E}_{\pm}^\star. \nonumber
\end{align}

\emph{Case 1.4}. [$P_VR(\gamma(\bar{t}-1)) =  \sqrt{|\L|}  $].\\ 
Using Assumption \ref{Asspm} d), we observe that 
\begin{align}
\begin{split}
\mathrm{H}_{\pm}(\gamma(\bar{t} - 1)) 
&\geq  \mathrm{H}_{\pm}((l_b^\star-1)   \times  \sqrt{|\L|} ) \\
&\geq \mathrm{H}_{\pm}((l_b^\star-1)   \times  l_h^\star ) + \lfloor (\sqrt{|\L|}-l_h^\star )/2 \rfloor (4J-\e (l_b^\star-1)) +\ho (l_b^\star-1)\\
&>  \mathrm{H}_{\pm}((l_b^\star-1)   \times  l_h^\star ) +\ho.
\end{split}
\end{align}
This leads to a contradiction, since
\begin{align}
\begin{split}
\mathrm{H}_{\pm}(\gamma(\bar{t})) 
&\geq \mathrm{H}_{\pm}(\gamma(\bar{t}-1)) + 2J-\he >  \mathrm{H}_{\pm}((l_b^\star-1)   \times  l_h^\star ) + \ho
+ 2J-\he = \mathrm{E}_{\pm}^\star .
\end{split}
\end{align}

\medskip
  \emph{Case 2}. [$ P_VR(\gamma(\bar{t}-1)) = l_h^\star-1 $ and $  P_HR(\gamma(\bar{t}-1)) =  l_b^\star +m $ for some $ m\geq 0 $].\\
Assume first that $ \gamma(\bar{t}) $ starts in $ \Lo $.
Hence, the top and the bottom row of $ \gamma(\bar{t}) $ belong to $ \Lo $.
Then, since $ \gamma(\bar{t}) $ is obtained from $ \gamma(\bar{t}-1) $ by adding a protuberance at a horizontal side of  $ R(\gamma(\bar{t}-1))  $, we have that 
\begin{align}\label{Gate1pmEq01}
\mathrm{H}_{\pm}(\gamma(\bar{t})) &\geq   \mathrm{H}_{\pm}(\gamma(\bar{t}-1)) + 2J+\ho.
\end{align}
Note that either the top or the bottom row of $ \gamma(\bar{t}-1) $ belongs to $ \Lo $.
By cutting this row, 
we can estimate the right-hand side of \eqref{Gate1pmEq01} from below by
\begin{align}\label{Gate1pmEq02}
\mathrm{H}_{\pm}((l_b^\star+m)   \times  (l_h^\star-2) ) + 4J + 2\ho.	
\end{align}
Moreover, \eqref{Gate1pmEq02} is bounded from below by 
\begin{align}\label{Gate1pmEq03}
	\mathrm{H}_{\pm}(l_b^\star   \times  (l_h^\star-2) ) + m( \mu-\e(l_b^\star - 1 ) )+ 4J + 2\ho, 
\end{align}
which is, obviously, strictly greater that $ \mathrm{E}_{\pm}^\star $.
This leads to a contradiction, and we can therefore, from now on, assume  that $ \gamma(\bar{t}) $ starts in $ \Le $. 

Note that
$ \gamma(\bar{t}) $ is obtained from $ \gamma(\bar{t}-1) $ either by  adding a protuberance at the  above horizontal side of $ R(\gamma(\bar{t}-1)) $ or the below one.
Without restriction, we suppose that a protuberance is added at the above horizontal side of $ R(\gamma(\bar{t}-1)) $.
Moreover, let $ P_HR(\gamma(\bar{t})) \times 	(l_h^\star  -2) $ denote the rectangle 
that is obtained from $ R(\gamma(\bar{t}-1)) $ by flipping all $(+1)$--spins from the top row of $ R(\gamma(\bar{t}-1)) $.
Note that $ P_HR(\gamma(\bar{t})) \times 	(l_h^\star  -2) $ starts from $ \Le $,  $  P_HR(\gamma(\bar{t}))= l_b^\star +m $ and that
$ \mathrm{H}_{\pm}(\gamma(\bar{t}))\geq\mathrm{H}_{\pm}(\gamma(\bar{t}-1))  + 2J - \he . $

\emph{Case 2.1}. [$ |\gamma(\bar{t}) \setminus \{P_HR(\gamma(\bar{t})) \times 	(l_h^\star  -2)\}|> 2 $].\\
In this case we necessarily have that $ \gamma(\bar{t}-1) $ has at least two $ (+1) $--spins in its uppermost row. 
If $ m=0 $, then,
since 
$ \gamma  $ does not cross $ \cP_2' $, we have by Lemma \ref{MinRectpm} that 
\begin{align}
\mathrm{H}_{\pm}(\gamma(\bar{t}-1)) >  \mathrm{H}_{\pm}(l_b^\star   \times  (l_h^\star-2) ) + 2J + 2\ho = \mathrm{E}_{\pm}^\star -2J + \he.
\end{align}
This leads to a contradiction, since 
\begin{align}
\mathrm{H}_{\pm}(\gamma(\bar{t}))\geq \mathrm{H}_{\pm}(\gamma(\bar{t}-1)) +2J- \he >  \mathrm{E}_{\pm}^\star.
\end{align}
If $ m>0 $ and $ P_HR( \gamma(\bar{t}-1) )< \sqrt{|\L|} $, then similarly, we observe
\begin{align}\label{Gate1pmEq04}
\begin{split}
\mathrm{H}_{\pm}(\gamma(\bar{t})) &\geq \mathrm{H}_{\pm}(\gamma(\bar{t}-1)) +2J- \he \geq  \mathrm{H}_{\pm}((l_b^\star+m)   \times  (l_h^\star-2) ) + 2J + 2\ho +2J- \he  \\
&= 
\mathrm{H}_{\pm}(l_b^\star   \times  (l_h^\star-2) ) + m(\mu-\e(l_b^\star - 1))+ 4J + 2\ho -\he > \mathrm{E}_{\pm}^\star,
\end{split}
\end{align}
which is a contradiction.
Finally, if $ P_HR( \gamma(\bar{t}-1) )= \sqrt{|\L|} $, we have that 
\begin{align}\label{Gate1pmEq05}
\mathrm{H}_{\pm}&(\gamma(\bar{t})) \geq \mathrm{H}_{\pm}(\gamma(\bar{t}-1)) +2J- \he  
\geq  \mathrm{H}_{\pm}( \sqrt{|\L|} \times  (l_h^\star-2) ) + 4J + 2\ho  - \he  \\
&= 
\mathrm{H}_{\pm}(l_b^\star   \times  (l_h^\star-2) ) + (\sqrt{|\L|} - l_b^\star )(\mu-\e(l_b^\star - 1)) - 2J (l_h^\star -1)+ 4J + 2\ho -\he\nonumber > \mathrm{E}_{\pm}^\star.\nonumber
\end{align}

\emph{Case 2.2}. [$ |\gamma(\bar{t}) \setminus \{P_HR(\gamma(\bar{t})) \times 	(l_h^\star  -2)\}|= 2 $].\\
Define
\begin{align} T= 
\max\Big\{ j \geq \bar{t}  \, \Big| \,  |\gamma(j) \setminus \{P_HR(\gamma(\bar{t})) \times 	(l_h^\star  -2)\}|\leq 2	\Big\}  
\end{align}
i.e.\ the last time that a configuration has only two $ (+1) $--spins outside of $ P_HR(\gamma(\bar{t})) \times 	(l_h^\star  -2) $.
From the maximality property of $ \bar{t} $, we have that $  P_VR(\gamma(T)) =  l_h^\star$
and $  P_HR(\gamma(T)) =  l_b^\star +m' $ for some $ m' \geq 0 $.
Moreover, we easily observe that $ \mathrm{H}_{\pm}(\gamma(T+1)) \geq \mathrm{H}_{\pm}(\gamma(T))+\ho  $.
As in  Case 2.1, we show that  every possible value of $ m' $ leads to a contradiction.
If $ m'=0 $, then,
since 
$ \gamma  $ does not cross $ \cP_2' $, Lemma \ref{MinRectpm} implies that 
\begin{align}
\mathrm{H}_{\pm}(\gamma(T)) >  \mathrm{H}_{\pm}(l_b^\star   \times  (l_h^\star-2) ) + 4J -\he + \ho = \mathrm{E}_{\pm}^\star - \ho,
\end{align}
and therefore
\begin{align}
\mathrm{H}_{\pm}(\gamma(T+1)\geq \mathrm{H}_{\pm}(\gamma(T)) + \ho >  \mathrm{E}_{\pm}^\star.
\end{align}
If $ m'>0 $ and $ P_HR( \gamma(T) )< \sqrt{|\L|} $, then 
\begin{align}
\begin{split}
\mathrm{H}_{\pm}(\gamma(T+1)) &\geq \mathrm{H}_{\pm}(\gamma(T)) + \ho \geq  \mathrm{H}_{\pm}((l_b^\star+m)   \times  (l_h^\star-2) ) + 4J  - \he + \ho +\ho  \\
&> \mathrm{E}_{\pm}^\star.
\end{split}
\end{align}
And if $ P_HR( \gamma(T) )= \sqrt{|\L|} $, we have that 
\begin{align}
\begin{split}
\mathrm{H}_{\pm}(\gamma(T+1)) &\geq \mathrm{H}_{\pm}(\gamma(T)) +\ho 
\geq  \mathrm{H}_{\pm}( \sqrt{|\L|} \times  (l_h^\star-2) ) + 4J - \he + 2\ho   > \mathrm{E}_{\pm}^\star.
\end{split}
\end{align}

Finally, we briefly sketch the proof for the case when $ \gamma $ can consist of several clusters. 
Recall the definitions of $ (n_j)_j, ((\gamma^k(j))_{k\leq n_j})_j, ((\ell_V^k(j))_{k\leq n_j})_j, ((\ell_H^k(j))_{k\leq n_j})_j, \ell_V   $ and $ \ell_H $ from the proof of Lemma \ref{Gate1A}.
Similarly as in \eqref{TwoCluster00Eq} and \eqref{TwoCluster01Eq}, we can show that 
for all $ j \in \N $, 
\begin{align}\label{TwoClusterpm00Eq}
\mathrm{H}_{\pm}&( \,\gamma(j) \,)\geq 
\sum_{k=1}^{n_j} \mathrm{H}_{\pm}( \,R(\gamma^k(j) )\,) - (n_j-1) \, 
\mathrm{H}_{\pm}(\boxminus)\\
&=\mathrm{H}_{\pm}(\boxminus) + 2 J \left(\sum_{k=1}^{n_j}   \ell_V^k(j)  +  \sum_{k=1}^{n_j}    \ell_H^k(j)\right) 
+ \ho \sum_{k=1}^{n_j}  \ell_H^k(j)  \lfloor \ell_V^k(j)/2 \rfloor - \he  \sum_{k=1}^{n_j}  \ell_H^k(j)  \lceil \ell_V^k(j)/2 \rceil.
\nonumber
\end{align}
Analogously to \eqref{TwoCluster0Eq} and \eqref{Gate1pmEq00}, 
define 
\begin{align} 
\tilde{t} -1 = \max\big\{ j \geq 0 \ | \ \ell_V(j) < l_h^\star \text{ or }  \ell_H (j)< l_b^\star	\}.
\end{align}
We have  that either $ \ell_H(\tilde{t}-1) = l_b^\star -1$ or $ \ell_V(\tilde{t}-1) = l_h^\star-1 $.
Proceeding as in  the first part of this proof and in the end of the proof of Lemma \ref{Gate1A}, we can now show  that, under the hypothesis that $ \gamma \cap \{ \cP_1 \cup \cP_2 \}= \emptyset$,   both cases lead to the fact that $ \mathrm{H}_{\pm}( \,\gamma(\tilde{t}) \,) > \mathrm{E}_{\pm}^\star $, which is a contradiction.
We omit the details and  conclude the proof of this lemma. 
\epr 

The following observation concludes the proof of $ \Phi(\boxminus,\boxplus) - \mathrm{H}_{\pm}(\boxminus) \geq  \G^\star _{\pm}$.

\bl{Gate2pm}
Let  $ \gamma\in (\boxminus , \boxplus )_{\mathrm{opt}}  $.
In order to cross at a time $ \bar{t} $ a configuration $ \gamma(\bar{t}) $ such that  $ P_VR(\gamma(j)) \geq l_h^\star $ and 
$ P_HR(\gamma(j)) \geq l_b^\star $ for all $ j\geq \bar{t} $, there must be some time $ t'\geq \bar{t}-1 $ such that   $  \gamma(t')  \in \cP_1 \cup \cP_ 2$ and $  \gamma(t'+1) \in \cC_1 \cup \cC_ 2$. 
In particular, every optimal path between $ \boxminus $ and $ \boxplus $ has to cross $ \cC_1 \cup \cC_2$.
\el

\bpr
Consider the  time step $ \bar{t} $ defined in the proof of Lemma \ref{Gate1pm}.
It was shown that there necessarily exists  a time $ t'\geq \bar{t} -1$ such that $ \gamma(t') \in \cP_1 \cup \cP_ 2 $. 
Note that 
\begin{align}\label{Gate2pmEq}
 P_VR(\gamma(j)) \geq l_h^\star \  \text{ and } \
 P_HR(\gamma(j)) \geq l_b^\star  \qquad \text{ for all } j \geq t'+1 .
\end{align}
In the following we show that $ \gamma(t'+1) \in \cC_1 \cup \cC_2$.

\emph{Case 1}. [	$ \gamma(t') \in \cP_1 \cup \cP_2''$	].\\
In this case,   $ \mathrm{H}_{\pm}(\gamma(t')) = \mathrm{E}_{\pm}^\star - \ho. $
Then it is easy to see  that $ \gamma(t'+1) $ must belong to $  \cC_1 \cup \cC_ 2''$.
Indeed,  for any other spin flip that fulfills the constraint  \eqref{Gate2pmEq}, the energy level of $ \gamma(t'+1) $ would exceed $ \mathrm{E}_{\pm}^\star $, and this violates the fact that $ \gamma\in (\boxminus , \boxplus )_{\mathrm{opt}}  $.
Note that we have tacitly used Assumption \ref{Asspm} a).

\emph{Case 2}. [	$ \gamma(t') \in \cP_2'$	].\\
By the definition of $ \bar{t} $, we have that $ t'= \bar{t} -1 $. 
And since $ P_VR(\gamma(\bar{t})) = l_h^\star $ and $ P_VR(\gamma(t')) = l_h^\star-1 $, we necessarily have that $ \gamma(t'+1) \in \cC_2'$.
This concludes the proof. 
\epr

\subsection{Identification of $ \cP^ \star $ and	$ \cC^ \star $}

\label{Step3pm}
Recall the definition of $ \cP^ \star $ and $ \cC^ \star $ from Definition \ref{Crit}.
Repeating similar computations as in Section \ref{Step1pm}, it is clear that $  \cP_1 \cup \cP_ 2 \subset  \cP^ \star $.
Now let $ \s \in \cP^ \star  $ and $ x \in \L $ be such that $ \s^x \in \cC^ \star $.
Then there exists $ \gamma \in (\boxminus,\boxplus)_{\mathrm{opt}}$ and $ \ell \in \N $ such that 
\begin{enumerate}[(i)]
	\item $ \gamma(\ell)= \s $ and $ \gamma(\ell +1)=\s^x$,
	\item $ \mathrm{H}_{\pm}(\gamma(k)) < \mathrm{E}_{\pm}^\star  $ for all $  k \in \{ 0,\dots,\ell	 	\}  $,
	\item $ \Phi(\boxminus,\gamma(k) ) \geq  \Phi(\gamma(k) , \boxplus )$  for all $  k \geq \ell + 1  $.
\end{enumerate}
As in the proof of Lemma \ref{Gate1pm} and in Lemma \ref{Gate2pm}, let 
\begin{align}
\bar{t}-1=  \max\big\{ j \geq 0 \ | \ P_VR(\gamma(j)) < l_h^\star \text{ or }  P_HR(\gamma(j)) < l_b^\star	\}.
\end{align}
We know from  Lemma \ref{Gate2pm}
that there exists $ t' \geq \bar{t} $ such that $ \gamma(t') \in  \cP_1 \cup \cP_ 2 $ and
$ \gamma(t'+1) \in  \cC_1 \cup \cC_ 2 $.
We get from fact (ii) that  $ \ell \leq t' $.

If $ \ell= t' $, then we have that  $ \s \in \cP_1 \cup \cP_ 2 $ and  $ \s^x \in \cC_1 \cup \cC_ 2$.

If $ \ell< t' $,
then fact  (iii) is violated, since $ \Phi(\boxminus,\gamma(t' ) ) < \Phi(\gamma(t') , \boxplus )= \mathrm{E}_{\pm}^\star $.
%
%
Hence, it must be the case that $ \ell= t' $.
We conclude that $ \cP^ \star=\cP_1 \cup \cP_ 2 $ and	$ \cC^ \star= \cC_1 \cup \cC_ 2$.

\subsection{Verification of (H1)}

\label{Step4pm}
Obviously, $ S_{\mathrm{stab}}= \{ \boxplus	\} $, since $ \he>\ho $.
It remains to show that $ S_{\mathrm{meta}}= \{ \boxminus	\}. $
Let $ \s \in S $. There are four cases.

\emph{Case 1}. [$ \s $ contains a cluster, which is not a stable rectangle].\\
Lemma \ref{LocMinpm} implies that $ \s $ is not $ \ho $--stable, i.e.\ there exists $ \s'  \in S $ such that $ \mathrm{H}_{\pm}(\s' ) < \mathrm{H}_{\pm}(\s)$ and
 $ \Phi ( \s,\s') - \mathrm{H}_{\pm}(\s) \leq \ho < \G_{\pm}^\star $.

 \emph{Case 2}. [$ \s $ contains a cluster $ R $, which is a stable rectangle  with $P_VR \geq l_h^\star $ and
 
 \qquad \qquad   $ P_HR<\sqrt{|\L|} $].\\
Let $ \s' $ be obtained from $ \s $ by attaching at the right vertical side  of $ R $ a column of length $ P_VR $.
We start to attach on an even row on the right  vertical side  of $ R $ and then successively flip adjacent spins until the column is filled. 
Then
\begin{align}
\begin{split}
\mathrm{H}_{\pm}(\s' ) &\leq     \mathrm{H}_{\pm}(\s ) + \mu - \frac{P_VR+1}{2} \e 
\leq \mathrm{H}_{\pm}(\s ) +  \mu - l_b^\star  \e< \mathrm{H}_{\pm}(\s), \text{ and }\\
\Phi ( \s,\s') - \mathrm{H}_{\pm}(\s) &\leq 2J - \he < \G_{\pm}^\star.
\end{split}
\end{align}

 \emph{Case 3}. [$ \s $ contains a cluster $ R $, which is a stable rectangle with $ P_VR \leq l_h^\star -2$ and  

 \qquad \qquad $ P_HR<\sqrt{|\L|} $].\\
Let $ \s' $ be obtained from $ \s $ by cutting the right column of $ R $. Then
\begin{align}
\begin{split}
\mathrm{H}_{\pm}(\s' ) &=    \mathrm{H}_{\pm}(\s )- \mu + \frac{P_VR+1}{2} \e 
\leq \mathrm{H}_{\pm}(\s ) -  \mu + (l_b^\star-1)  \e< \mathrm{H}_{\pm}(\s),  \text{ and }\\
\Phi ( \s,\s') - \mathrm{H}_{\pm}(\s) &\leq \frac{P_VR-1}{2} \e + \ho < \G_{\pm}^\star.
\end{split}
\end{align}

 \emph{Case 4}. [$ \s $ contains a cluster $ R $, which is a stable rectangle with $ P_HR=\sqrt{|\L|} $].\\
Let $ \s' $ be obtained from $ \s $ by attaching above $ R $ successively vertical bars of length $ 2 $ until the two rows above $ R $ wrap around the torus.
Then,
\begin{align}
\begin{split}
\mathrm{H}_{\pm}(\s' ) &=  \mathrm{H}_{\pm}(\s ) + 4J - l_1 \e < \mathrm{H}_{\pm}(\s), \quad \text{and} \\
\Phi ( \s,\s') - \mathrm{H}_{\pm}(\s) &\leq 4J - \e  < \G_\pm^\star.
\end{split}
\end{align}
This proves that $ S_{\mathrm{meta}}= \{ \boxminus	\}. $

\subsection{Computation of $K$}

\label{Step6pm}

Again, we proceed as in Section \ref{Step6A} and in Section \ref{Step6NN}.
Before estimating  $ K^{-1} $ from below and above, we define 
$ \bar{\cC}= \bar{\cC_1} \cup \bar{\cC_2}$, where
\begin{itemize}
	\item 
	$ \bar{\cC_1}$ is the set of all configurations $ \s $ that are obtained from a configuration  $\s' \in  \cC_1 $ as follows.
	There is a column in $ \s' $ that has length 2.
	$ \s $ is obtained from $ \s' $  by adding a third $ (+1) $--spin on the even row adjacent to this column, and
	\item 
	$ \bar{\cC_2}$ is the set of all configurations $ \s $ that are obtained from a configuration  $\s' \in \cC_2 $ as follows. 
	There is a component of three $ (+1) $--spins above or below the $ l_b^\star \times (l_h^\star -2) $-rectangle in $ \s' $.
	$ \s $ is obtained from $ \s' $ by adding a $ (+1) $--spin such that this component becomes a $ 2 \times 2 $-square.
\end{itemize}
\begin{figure}[htbp]
	\centering
	\includegraphics[height= 4cm]{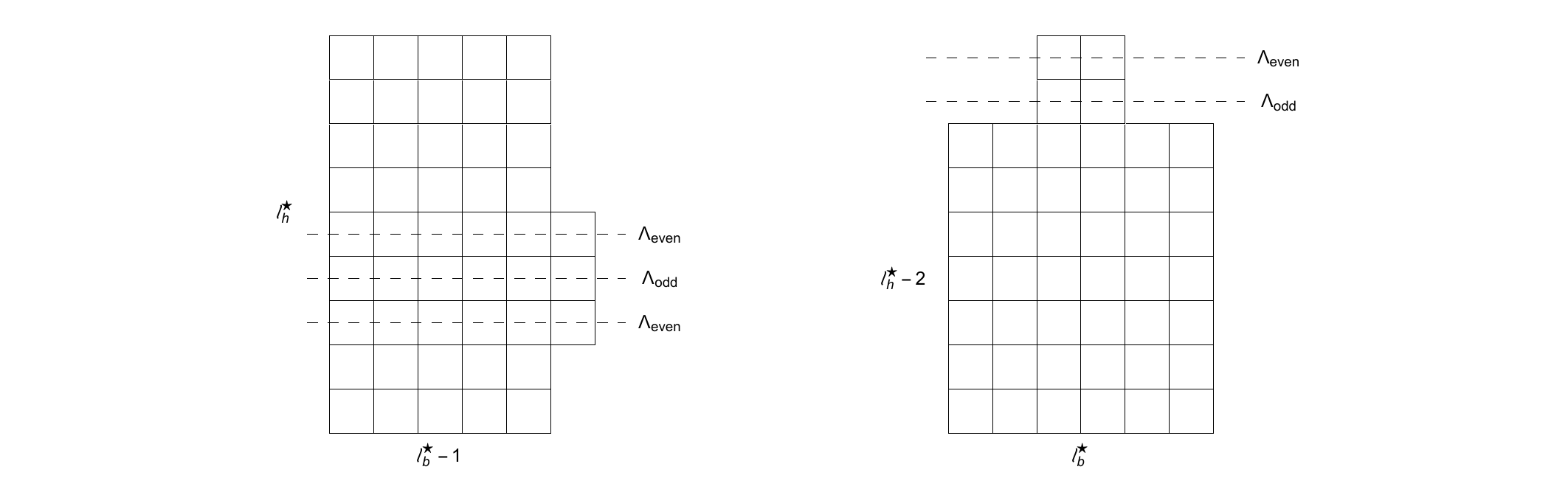}	
	\caption{{\footnotesize An example of an element in $ \bar{\cC_1}$ and $ \bar{\cC_2}$} }
\end{figure}
It is easy to see that $ \partial^+ \cC^\star \cap S^\star = \cP_1 \cup \cP_2 \cup \bar{\cC}_1 \cup \bar{\cC}_2= \cP^\star \cup \bar{\cC}$, \  $  \cP^\star \subset  S_\boxminus  $ \  and \  $  \bar{\cC} \subset S_\boxplus  $.

\emph{Lower bound}.
Using these definitions and facts, we can estimate $ K^{-1} $ as follows. 
\begin{align}
\begin{split}
\frac{1}{K} &\geq \min_{C_1,\dots,C_I \in [0,1]}
\min_{\substack{ h:S^\star \ra [0,1] \\
		\smallrestr{h}{S_\boxminus} =1, \smallrestr{h}{S_\boxplus} =0, \smallrestr{h}{S_i} =C_i \, \forall i }}
\frac{1}{2} \sum_{\eta,\eta' \in (\cC^\star)^+ } \mathbbm{1}_{\{	\eta \sim \eta'	 \}} [h(\eta) - h(\eta')]^2\\
 &= 
\min_{ h:\cC^\star \ra [0,1]}
\sum_{\eta \in \cC^\star } \left( \sum_{\eta' \in \cP^\star, \eta'\sim \eta  }	[1- h(\eta) ]^2	
+ \sum_{\eta' \in \bar{\cC} , \eta'\sim \eta  }	h(\eta)^2			\right) \\
&=
\sum_{\eta \in \cC_1 } \frac
{|\cP^\star\sim \eta| 	\cdot 
	|\bar{\cC} \sim \eta| }
{|\cP^\star\sim \eta| 	+ 
	|\bar{\cC} \sim \eta| }
+ \sum_{\eta \in \cC_2 } \frac
{|\cP^\star\sim \eta| 	\cdot 
	|\bar{\cC} \sim \eta| }
{|\cP^\star\sim \eta| 	+ 
	|\bar{\cC} \sim \eta| }.
\end{split}
\end{align}
For all $ \eta \in \cC_1  $ we have that $ |\cP^\star \sim \eta| = 1 $,
whereas for all $ \eta \in \cC_2  $ we have that $ |\cP^\star \sim \eta| = 2 $.
Moreover, $ |\bar{\cC} \sim \eta| = 1 $  for all $ \eta \in \cC^\star   $.
Finally, it can be seen easily that $ |\cC_1| = |\cC_2| = 4|\L| (l_b^\star-1) $.
Hence, 
\begin{align}
\frac{1}{K} \geq
|\cC_1| \frac{1} {2 } +  |\cC_2| \frac {2 } {3 }	
= 
\frac{14\,( l_b^\star-1)}{3}|\L|.
\end{align}

\emph{Upper bound}.
We say that a row or a column of a configuration is a \emph{singleton} if it consists only of a single $ (+1) $-- spin.
We define the following subsets of $ S^\star $. 
\begin{align}
\mathscr{S}^- =  \{\s\in S^\star\, &| \, \text{ for all clusters $ \eta  $ of $ \s $ we have that either } (P_VR(\eta) < l_h^\star) \nonumber \\
&\ \text{ or }  (P_HR(\eta) < l_b^\star) \nonumber\\
&\ \text{ or }  (P_HR(\eta) \geq l_b^\star ,	P_VR(\eta) = l_h^\star	\text{ and at least two rows of $ \eta $ are singletons}) \nonumber\\
&\ \text{ or }  (P_HR(\eta) = l_b^\star , P_VR(\eta) = l_h^\star 	\text{ and at least one column of $ \eta $ is a singleton}) \}, \nonumber\\
\mathscr{S}^+_1 = \{\s\in S^\star\, &| \, \text{ there exists a cluster $ \eta  $ of $ \s $ such that  } P_HR(\eta) = l_b^\star ,  P_VR(\eta) = l_h^\star \nonumber \\
&\ \text{ and no column of $ \eta $ is a singleton} \\
&\ \text{ and at most one row of $ \eta $ is a singleton} \},\nonumber \\
\mathscr{S}^+_2 = \{\s\in S^\star\, &| \, \text{ there exists a cluster $ \eta  $ of $ \s $ such that  } P_HR(\eta) > l_b^\star 	\text{ and } P_VR(\eta) = l_h^\star  \nonumber \\
&\ \text{ and at most one row of $ \eta $ is a singleton} \},\nonumber \\
\mathscr{S}^+_3 = \{\s\in S^\star\, &| \,  \text{ there exists a cluster $ \eta  $ of $ \s $ such that  } P_HR(\eta) = l_b^\star \text{ and } P_VR(\eta) >  l_h^\star \},\nonumber\\
\mathscr{S}^+_4= \{\s\in S^\star\, &| \, \text{ there exists a cluster $ \eta  $ of $ \s $ such that  } P_HR(\eta) > l_b^\star 	\text{ and } P_VR(\eta) > l_h^\star \}.\nonumber
\end{align} 
Set 
$ \mathscr{S}^+ =   \mathscr{S}^+_1 \cup \mathscr{S}^+_2 \cup \mathscr{S}^+_3\cup \mathscr{S}^+_4$.
Then, $ S^\star = \mathscr{S}^-\cup \mathscr{S}^+$, $ \mathscr{S}^- \cap \mathscr{S}^+ = \emptyset $, $ \cP^\star \subset \mathscr{S}^- $ and $ \cC^\star \subset \mathscr{S}^+ $.
\bl{UpperBoundpm}
Let $ \s \in \mathscr{S}^- $ and $ \s' \in \mathscr{S}^+ $.
Then $ \s \sim \s' $ if and only if $\s \in  \cP^\star  $ and $\s' \in  \cC^\star  $.
\el
\bpr
In the following we show   separately for  all different cases that the assumption that either $ \s \notin \cP^\star $ or $ \s' \notin \cC^\star $ leads to $ \s  \notin S^\star $ or $ \s'  \notin S^\star $, which is a contradiction.
Using the same arguments as in the proof of Lemma \ref{Gate1pm}, it is no restriction to assume that both $ \s'  $ and $ \s $ consist of a unique cluster and that   $ R(\s) $ starts in $ \Le $.
\medskip

 \emph{Case 1}. [	$ \s' \in \mathscr{S}^+_1 $	].
 
 \emph{Case 1.1}. [	$ P_VR(\s) < l_h^\star$	].\\
$ \s' $ is obtained from $ \s  $ by adding a row to $ \s $, which is a singleton.
Therefore, since $ \s' \in \mathscr{S}^+_1 $, each row of $ \s $ needs to have at least two $ (+1) $--spins. 
 Now the same computations as in Case 2.1 from the proof of Lemma \ref{Gate1pm} lead to a contradiction.  
Here  $ \gamma(\bar{t}) $ is replaced by $ \s' $ and $ \gamma(\bar{t}-1) $ is replaced by  $ \s $.
 
 \emph{Case 1.2}. [	$ P_HR(\s) < l_b^\star$	].\\
$ \s' $ is obtained from $ \s  $ by adding a column to $ \s $, which is a singleton.
Since no column of $ \s' $ is a singleton, $ \s \sim \s' $ can not hold true, which implies that this case is not possible.

 \emph{Case 1.3}. [	$ P_HR(\s) \geq l_b^\star ,	P_VR(\s) = l_h^\star$  and at least two rows of $ \s $ are singletons	].\\
$ \s' $ is obtained from $ \s $ by flipping  a $ (-1) $--spin in a row of $ \s $  that is a singleton.
The same computations as in the Case 2.2 from the proof of Lemma \ref{Gate1pm} lead to a contradiction. 
 Here  $ \gamma(T+1) $ is replaced by $ \s' $ and $ \gamma(T) $ is replaced by  $ \s $.
 
 \emph{Case 1.4}. [	$ P_HR(\s) = l_b^\star ,	P_VR(\s) = l_h^\star$  and at least one column  of $ \s $ is a  singleton	].\\
$ \s' $ is obtained from $ \s $ by flipping  a $ (-1) $--spin in a column of $ \s $  that is a singleton.
The same computations as in the Case 1.1 from the proof of Lemma \ref{Gate1pm}  lead to a contradiction.  
Here  $ \gamma(\tau) $ is replaced by $ \s' $ and $ \gamma(\tau-1) $ is replaced by  $ \s $.
\medskip

 \emph{Case 2}. [	$ \s' \in \mathscr{S}^+_2 $	].

 \emph{Case 2.1}. [	$ P_VR(\s) < l_h^\star$	].\\
 We necessarily have that $ P_HR(\s) = l_b^\star + m$ for some $ m>0 $.
  $ \s' $ is obtained from $ \s  $ by adding a row in $ \Le $ to $ \s $, which is a singleton.
  Since $ \s' \in \mathscr{S}^+_2 $, we have that  the odd row below the added row  contains at least two $ (+1) $--spins.
 Now the same computations as in the equations \eqref{Gate1pmEq04}--\eqref{Gate1pmEq05} lead to a contradiction.  
 Here  $ \gamma(\bar{t}) $ is replaced by $ \s' $ and $ \gamma(\bar{t}-1) $  by  $ \s $.
 
 \emph{Case 2.2}. [	$ P_HR(\s) < l_b^\star$	].\\
  It is easy to see that $ \s \sim \s' $ can not hold true in this case.

 \emph{Case 2.3}. [	$ P_HR(\s) \geq l_b^\star ,	P_VR(\s) = l_h^\star$  and at least two rows of $ \s $ are singletons	].\\
 See Case 1.3.
 
 \emph{Case 2.4}. [	$ P_HR(\s) = l_b^\star ,	P_VR(\s) = l_h^\star$  and at least one column  of $ \s $ is a  singleton	].\\
 $ \s' $ is obtained from $ \s  $ by adding a column to $ \s $, which is a singleton.
 Hence,  two columns  of $ \s' $  are singletons.
 This implies that 
 \begin{align}
 \begin{split}
 \mathrm{H}_{\pm}(\s' ) & 
 \geq \mathrm{H}_{\pm}((l_b^\star -1 )\times l_h^\star  ) + 4J -2\he =  \mathrm{E}_{\pm}^\star 	+2J-\ho - \he 
 > \mathrm{E}_{\pm}^\star.
 \end{split}
 \end{align}

 \emph{Case 3}. [	$ \s' \in \mathscr{S}^+_3 $	].
 
 \emph{Case 3.1}. [	$ P_VR(\s) < l_h^\star$	].\\
 It is easy to see that $ \s \sim \s' $ can not hold true in this case.
 
 \emph{Case 3.2}. [	$ P_HR(\s) < l_b^\star$	].\\
 $ \s' $ is obtained from $ \s  $ by adding a protuberance at the left vertical side or the right vertical side of $ R(\s) $.
 The same computations as in the Cases 1.2, 1.3 and 1.4 from the proof of Lemma \ref{Gate1pm} lead to a contradiction.
  Here  $ \gamma(\bar{t}) $ is replaced by $ \s' $ and $ \gamma(\bar{t}-1) $  by  $ \s $. 
 
 \emph{Case 3.3}. [	$ P_HR(\s) \geq l_b^\star ,	P_VR(\s) = l_h^\star$  and at least two rows of $ \s $ are singletons	].\\
 $ \s' $ is obtained from $ \s $ by adding  a protuberance at the top row or bottom row of $ R(\s) $.
 Let $ P_HR(\s) = l_b^\star + m $ for some $ m \geq 0  $.
 Obviously, $ \mathrm{H}_{\pm}(\s ) \geq \mathrm{H}_{\pm}((l_b^\star +m )\times (l_h^\star - 2 ) ) + 4J - \e $.
 This implies that 
 \begin{align}
 \begin{split}
 \mathrm{H}_{\pm}(\s' ) &\geq  \mathrm{H}_{\pm}(\s ) +2J -  \he
 \geq \mathrm{H}_{\pm}((l_b^\star +m )\times (l_h^\star - 2 ) ) + 6J -  \e -\he\\
 &=  \mathrm{E}_{\pm}^\star + m(\mu - \e (l_b^\star - 1))	+2J-\ho - \he 
 > \mathrm{E}_{\pm}^\star.
 \end{split}
 \end{align}
 
 \emph{Case 3.4}. [	$ P_HR(\s) = l_b^\star ,	P_VR(\s) = l_h^\star$  and at least one column  of $ \s $ is a  singleton	].\\
 $ \s' $ is obtained from $ \s $ by adding  a protuberance at the top row or bottom row of $ R(\s) $.
 Hence,  one column  and one row of $ \s' $  are singletons.
 Then, as in Case 2.4 above, 
 \begin{align}
 \begin{split}
 \mathrm{H}_{\pm}(\s' ) & 
 \geq \mathrm{H}_{\pm}((l_b^\star -1 )\times l_h^\star  ) + 4J -2\he 
 > \mathrm{E}_{\pm}^\star.
 \end{split}
 \end{align}

 \emph{Case 4}. [	$ \s' \in \mathscr{S}^+_4 $	].
 
 \emph{Case 4.1}. [	$ P_VR(\s) < l_h^\star$	].\\
 It is easy to see that $ \s \sim \s' $ can not hold true in this case.

\emph{Case 4.2}. [	$ P_HR(\s) < l_b^\star$	].\\
$ \s \sim \s' $ can not hold true in this case.

\emph{Case 4.3}. [	$ P_HR(\s) \geq l_b^\star ,	P_VR(\s) = l_h^\star$  and at least two rows of $ \s $ are singletons	].\\
 See Case 3.3.

\emph{Case 4.4}. [	$ P_HR(\s) = l_b^\star ,	P_VR(\s) = l_h^\star$  and at least one column  of $ \s $ is a  singleton	].\\
$ \s \sim \s' $ can not hold true in this case.
\epr

As in Section \ref{Step6A} and in Section \ref{Step6NN}, we have that $ S_\boxminus \subset \mathscr{S}^- $, $ S_\boxplus \subset \mathscr{S}^+ $.
Moreover, for all $ i =0,\dots,I $ either $ S_i \subset \mathscr{S}^- $  or $ S_i \subset \mathscr{S}^+ $ holds true.
Therefore, again as in Section \ref{Step6A} and in Section \ref{Step6NN}, we estimate the minimum in \eqref{EqPrefactor} from above by the minimum over all functions of the form \eqref{EqStep6FunctionA}.
Using Lemma \ref{UpperBoundpm} we infer that 
%
	\begin{align}
	\begin{split}
\frac{1}{K} &\leq 
\min_{\substack{ h:S^\star \ra [0,1] \\
		\smallrestr{h}{\mathscr{S}^-} =1, \smallrestr{h}{\mathscr{S}^+\setminus \cC^\star } =0}}
\frac{1}{2} \sum_{\eta,\eta' \in S^\star } \mathbbm{1}_{\{	\eta \sim \eta'	 \}} [h(\eta) - h(\eta')]^2 \\
&=
\min_{\substack{ h:(\cC^\star)^+ \ra [0,1] \\
		\smallrestr{h}{\mathscr{S}^-\cap \partial^+\cC^\star} =1, \smallrestr{h}{\mathscr{S}^+\cap \partial^+\cC^\star} =0}}
\frac{1}{2} \sum_{\eta,\eta' \in (\cC^\star)^+ } \mathbbm{1}_{\{	\eta \sim \eta'	 \}} [h(\eta) - h(\eta')]^2 \\
&= 
\min_{ h:\cC^\star \ra [0,1]}
\sum_{\eta \in \cC^\star } \left( \sum_{\eta' \in \cP^\star, \eta'\sim \eta  }	[1- h(\eta) ]^2	
+ \sum_{\eta' \in \bar{\cC} , \eta'\sim \eta  }	h(\eta)^2			\right) \\
&=
\frac{14\,( l_b^\star-1)}{3}|\L|.
\end{split}
\end{align}

\section*{Acknowledgment} 
The author would like to give many thanks to Anton Bovier and Muhittin Mungan for a lot of useful discussions and  suggestions.
Moreover,
he would like to thank the anonymous referee
for valuable comments.


\end{document}